\documentclass[a4paper,english]{amsart}

\usepackage[utf8]{inputenc}
\usepackage[T1]{fontenc}

\usepackage{amsfonts}	
\usepackage{amsmath}	
\usepackage{amssymb}	
\usepackage{amsthm}	
\usepackage{bbm}	
\usepackage{mathrsfs}
\usepackage{mathtools}	
\usepackage{stmaryrd}

\usepackage[pdftex, pdfstartview=FitW,colorlinks=true,linkcolor=blue,%
urlcolor=blue,citecolor=blue,
pdfauthor={Timothée Bénard and Weikun He},
pdftitle={Multislicing and effective equidistribution for random walks on some homogeneous spaces},
pdfkeywords={Random walk, Homogeneous space, Slicing theorem, Equidistribution}
]{hyperref}
\usepackage{cleveref}	
\usepackage{url}	
\usepackage{enumitem}	
\setenumerate[0]{label=\alph*)}   

\newtheorem{thm}{Theorem}[section]
 
\newtheorem{corollary}[thm]{Corollary} 
\newtheorem{lemma}[thm]{Lemma} 
\newtheorem{proposition}[thm]{Proposition}

\theoremstyle{definition}
\newtheorem{definition}[thm]{Definition}
\newtheorem{notation}[thm]{Setting}

\newcommand{\QQ}{\overline{\mathbb{Q}}}
\newcommand{\Q}{\mathbb{Q}}
\newcommand{\T}{\mathbb{T}}
\newcommand{\N}{\mathbb{N}}
\newcommand{\Z}{\mathbb{Z}}
\newcommand{\C}{\mathbb{C}}
\newcommand{\R}{\mathbb{R}}
\renewcommand{\L}{\mathbb{L}}
\newcommand{\K}{\mathbb{K}}
\newcommand{\F}{\mathbb{F}}

\newcommand{\eps}{\varepsilon}
\newcommand{\ut}{\mathbf{t}}
\newcommand{\ur}{\mathbf{r}}
\newcommand{\us}{\mathbf{s}}
\newcommand{\ka}{\mathfrak{a}}

\newcommand{\kg}{\mathfrak{g}}
\newcommand{\kh}{\mathfrak{h}}
\newcommand{\kP}{\mathfrak{P}}
\newcommand{\so}{\mathfrak{so}}

\newcommand{\wt}{\widetilde}

\newcommand{\ug}{\underline g}

\newcommand{\1}{\mathbbm{1}}

\renewcommand{\a}{\alpha}
\renewcommand{\b}{\beta}
\renewcommand{\d}{\delta}
\renewcommand{\k}{\kappa}

\renewcommand{\phi}{\varphi}
\newcommand{\cN} {{\mathcal N}}
\newcommand{\cF} {{\mathcal F}}

\newcommand{\cB} {{\mathcal B}}

\newcommand{\cV} {{\mathcal V}}

\newcommand{\cP} {{\mathcal P}}

\newcommand{\cD} {{\mathcal D}}
\newcommand{\cQ} {{\mathcal Q}}
\newcommand{\cE} {{\mathcal E}}

\newcommand{\cS} {{\mathcal S}}

\newcommand{\cR} {{\mathcal R}}
\newcommand{\cW} {{\mathcal W}}
\newcommand{\cO}{\mathcal{O}}

\newcommand{\sD} {{\mathscr D}}
\newcommand{\sC}{\mathscr{C}}

\newcommand{\bG}{\mathbf{G}}
\newcommand{\bH}{\mathbf{H}}
\newcommand{\bK}{\mathbf{K}}
\newcommand{\bR}{\mathbf{R}}
\newcommand{\bbF}{\mathbb{F}}

\DeclareMathOperator{\Ad}{Ad}

\DeclareMathOperator{\End}{End}

\DeclareMathOperator{\GL}{GL}
\DeclareMathOperator{\Gal}{Gal}

\DeclareMathOperator{\Lie}{Lie}
\DeclareMathOperator{\Lip}{Lip}
\DeclareMathOperator{\Mah}{Mah}

\DeclareMathOperator{\SL}{SL}
\DeclareMathOperator{\SO}{SO}

\DeclareMathOperator{\ad}{ad}
\DeclareMathOperator{\codim}{codim}

\DeclareMathOperator{\diam}{diam}
\DeclareMathOperator{\inj}{inj}

\DeclareMathOperator{\supp}{supp}
\DeclareMathOperator{\Span}{Span}
\renewcommand{\sl}{\mathfrak{sl}}
\DeclareMathOperator{\dist}{d}
\DeclareMathOperator{\M}{M}
\DeclareMathOperator{\den}{den}

\newcommand{\abs}[1]{\lvert#1\rvert}    
\newcommand{\abse}[1]{\left\lvert#1\right\rvert}

\newcommand{\norm}[1]{\lVert#1\rVert}   

\newcommand{\set}[1]{\{\, #1  \,\}}     
\newcommand{\setbig}[1]{\bigl\{\, #1 \,\bigr\}}

\newcommand{\dd}{\,\mathrm{d}}
\newcommand{\defeq}{\coloneqq}
\DeclareMathOperator{\Id}{Id}
\newcommand{\acts}{\curvearrowright}

\newcommand{\presigma}[1]{\prescript{\sigma}{}{#1}}

\DeclareMathOperator{\Gr}{Gr}
\DeclareMathOperator{\dang}{d_\measuredangle}
\DeclareMathOperator{\vol}{vol}

\renewcommand{\setminus}{\smallsetminus}
\renewcommand{\emptyset}{\varnothing}       
\renewcommand{\subset}{\subseteq}

\newcommand{\increasing}{\boxslash}         

\newcommand{\bfparagraph}[1]{\par\medskip \noindent{\bfseries #1}}

\newcommand{\overprec}[1]{\overset{#1}{\prec}}

\newenvironment{remark}{\par\medskip \noindent \textit{Remark.} \rmfamily}{\medskip}
   
\title[Multislicing and effective equidistribution for random walks]{Multislicing and effective equidistribution for random walks on some homogeneous spaces}

\author{Timoth\'ee B\'enard }
\address{CNRS – LAGA, Universit\'e Sorbonne Paris Nord, 99 avenue J.-B.
Cl\'ement, 93430 Villetaneuse}
\email{benard@math.univ-paris13.fr}

\author{Weikun He}
\address{State Key Laboratory of Mathematical Sciences, Academy of Mathematics and System Science, Chinese Academy of Sciences, Beijing 100190, P. R. of China}
\email{heweikun@amss.ac.cn}
\thanks{W.H. is supported by the National Key R\&D Program of China (No. 2022YFA1007500) and the National Natural Science Foundation of China (No. 12288201).}
\date{}

\begin{document}

\large


\begin{abstract}
We consider a random walk on a homogeneous space $G/\Lambda$ where $G$ is  $\SO(2,1)$ or $\SO(3,1)$ and  $\Lambda$ is a lattice. The walk is driven by a probability measure $\mu$ on $G$ whose support generates a Zariski-dense subgroup.

We show that for every starting point $x\in G/\Lambda$ which is not trapped in a finite $\mu$-invariant set, the $n$-step distribution $\mu^{*n}*\delta_{x}$ of the walk equidistributes toward the Haar measure.  Moreover, under arithmetic assumptions on the pair $(\Lambda, \mu)$, we show the convergence occurs at an exponential rate, tempered by the obstructions that $x$ may be high in a cusp or close to a finite orbit.

Our approach is substantially different from that of Benoist-Quint \cite{BQ2}, whose equidistribution statements only hold in Ces\`aro average and are not quantitative, that of Bourgain-Furman-Lindenstrauss-Mozes \cite{BFLM} concerning the torus case, and that of  Lindenstrauss-Mohammadi-Wang and Yang \cite{LMW,Yang,LMWY} about the analogous problem for unipotent flows.
A key new feature of our proof is the use of a new phenomenon which we call multislicing. The latter is a generalization of the discretized projection theorems à la Bourgain and we believe it presents independent interest.

\end{abstract}

\maketitle

\setcounter{tocdepth}{2}    
\tableofcontents

\section{Introduction} 

A major topic during the last half of the twentieth century has been to understand the orbits of unipotent flows on homogeneous spaces. The subject started with the work of Hedlund \cite{Hedlund1, Hedlund2, Hedlund3} and Furstenberg \cite{Furstenberg73} who treated the case where the homogeneous space is a compact quotient of $\SL_{2}(\R)$, and culminated with Ratner's theorems \cite{Rat1, Rat2, Rat3, Rat4} stating that unipotent flows on finite volume homogeneous spaces have homogeneous orbit closures and equidistribute inside them. A natural question, then formulated by Shah \cite{Shah98}, was to extend these conclusions to orbits of a group $\Gamma$ whose Zariski-closure is generated by unipotents, thus allowing for groups that may not even contain unipotents themselves. 
Progress on this question was made in the breakthrough work of Benoist-Quint.
Instead of a unipotent flow, they considered a random walk driven by a probability measure supported on $\Gamma$, for which a whole array of probabilistic tools are available.
Benoist and Quint extended Ratner's theorems to this probabilistic setting \cite{BQ1,BQ2,BQ3} under the assumption that $\Gamma$ has semisimple Zariski closure, thus answering positively Shah's question in this case.  Since then, Benoist-Quint's work has been adapted to study fractal measures \cite{SW19, ProSertShi23}, the dynamics of surface diffeomorphisms \cite{BrRod}, regularity of orbit of closures on moduli spaces \cite{EskMirz}.  Further generalizations include \cite{EskinLindenstrauss, BenDeS21}.


In parallel,  the question of making  Ratner's theorems effective emerged. Indeed Ratner's theorems give equidistribution statements of unipotent flows as the time parameter goes to infinity, but the proof is purely ergodic theoretic and yields no information as to when a given level of equidistribution is achieved. Effectivity is about giving an explicit rate of equidistribution, that would be computable in terms of the initial data. In the nilmanifold setting, the work of Green-Tao~\cite{GreenTao12} provides such estimates using Fourier analysis and the van der Corput inequality. The case of horospherical flows (as unipotent flows on a quotient of $\SL_{2}(\R)$) can be handled using a now standard thickening argument and exponential decay of matrix coefficients \cite{KM12, Tay19, Katz23}. More general homogeneous spaces raise serious difficulties. Major recent developments concern cases of small dimension, including the work of Lindenstrauss, Mohammadi, and Wang \cite{LM, LMW} dealing with homogeneous spaces modeled over $\SO(2, 1)\times \SO(2,1)$ or $\SO(3,1)$, and the work of Yang \cite{Yang} (and also \cite{LMWY}) tackling certain unipotent flows on $\SL_{3}(\R)/\SL_{3}(\Z)$. 

In the random walk setting, the first effective result was obtained in the breakthrough work of Bourgain-Furman-Lindenstrauss-Mozes \cite{BFLM} managing the case of the torus. The four authors showed that a random walk on $\T^d$ driven by a measure on $\SL_{d}(\Z)$ with support generating a big enough subgroup
equidistributes at an exponential rate toward the Lebesgue probability measure on $\T^d$, unless the starting point is close to a finite orbit of the walk. Since then, \cite{BFLM} has been generalized to semisimple linear random walks on the torus~\cite{He2020IJM,HS2022,HS2023} as well as to affine random walks on the torus~\cite{Boyer,HLL_Aff} and on some nilmanifolds~\cite{HLL_Nil}. 
All these works can be seen as probabilistic counterparts of Green-Tao's results for nilflows, and similarly, their approaches all rely crucially on Fourier analysis.

In this paper, we obtain the first effective equidistribution result for random walks with arbitrary starting point on spaces that are not nilmanifolds but instead are modeled on semisimple Lie groups, namely $\SO(2,1)$ or $\SO(3,1)$.  Effectivity aside, our equidistribution result also improves on that of Benoist-Quint because it holds without Ces\`aro average.
The overall strategy of proof is rather common:  first, establish that the random walk distribution acquires some positive dimension at a certain range of scales, then run the walk again to bootstrap the dimension arbitrarily close to that of the ambient space, and finally use a spectral gap argument to go from high dimension to equidistribution. However, the arguments we develop are significantly different from our predecessors, as we now explain.  

The method of \cite{BFLM} and subsequent works have not allowed progress on the problem due to the lack of a good analog for Fourier analysis on semisimple homogeneous spaces. Here we attack the problem from a different angle, \emph{we do not use Fourier analysis}. We develop new tools rooted in additive combinatorics
and apply them directly in the ambient space, rather than in the Fourier domain.
 The method also applies in the context of the torus, providing a new proof of \cite{BFLM} (Zariski-dense case).

 
Our approach differs also from \cite{LMW} and \cite{Yang} about unipotent flows. First in both \cite{LMW, Yang}, there is a well-defined notion of unstable manifold and the setting only requires bootstrapping dimension in a codimension $1$ subspace of the unstable direction, thus leading to considerations in dimension $1$ (as in \cite{LMW}) or 2 (in \cite{Yang}).  In our situation, there is \emph{no privileged unstable direction} (as it depends on the random instructions of the walk), so we must bootstrap dimension in all possible directions, thus forcing us to consider spaces of bigger dimension  ($3$ or $6$).
Second, as in \cite{LMW}, we use projection theorems (relying on additive combinatorics/incidence geometry) to gain dimension. 
However, \emph{the projection theorem we use is different and the way we use it is different}.
Indeed, in \cite{LMW} the unipotent flow projects as a non-degenerate curve in the direction where the bootstrap is performed, and this allows to apply a restricted Marstrand-type projection theorem, granting that dimension is preserved. Such an approach is not possible for us because it uses crucially that the expanded measure on the unipotent orbit is Lebesgue, while a fair analog in our context would be to expand a measure of positive but non-full Hausdorff dimension (actually a Furstenberg measure). Here non-full dimension forbids dimensional preservation. Instead, we develop a new approach based on a \emph{multislicing theorem} which generalizes the projection theorem of Bourgain \cite{Bourgain2010} (as well as its extensions by the second named author~\cite{He2020JFG} and Shmerkin \cite{Shmerkin}). 
More details on that theorem and how we use it are given in \Cref{Sec-ideas-more}. We believe our multislicing theorem will find many more applications. In a follow-up work with Han Zhang \cite{BHZ24, BHZ25}, we  use it to prove an analog of Khintchine's theorem for the middle-thirds Cantor set, hereby answering a long-standing question of Mahler \cite[Section 2]{Mahler84} also put forward by Kleinbock-Lindenstrauss-Weiss in \cite{KLW04}.  Works related to this question comprise \cite{KLW04, EFS11, KM12, SW19, Yu21, KhalilLuethi, DJ24}. 

\subsection{Statement of the main results}\label{ss:mainresults}
Let $G$ be a  connected real linear group which is isogenous to $\SO(2,1)$ or $\SO(3,1)$, let  $\Lambda\subseteq G$ be a lattice, set $X=G/\Lambda$.
The group $G$ acts on $X$ by left multiplication.
We write $m_X$ to denote the unique $G$-invariant Borel probability measure on $X$, commonly known as the \emph{Haar measure}.
Slightly abusively, we denote by $*$ both the convolution associated to the multiplication map $G \times G \to G$ and that associated to the action map $G \times X \to X$.
For a point $x \in X$, we write $\delta_x$ to denote the Dirac mass at $x$.
Given a probability measure $\mu$ on $G$,  we study the Markov chain (\emph{$\mu$-walk}) on $X$ with transition  distributions $(\mu*\delta_{x})_{x\in X}$.
Thus, the distribution at time $n \in \N$ of the $\mu$-walk starting at a point $x \in X$ is 
\[
\mu^{*n} * \delta_x = \underbrace{\mu * \dotsm * \mu}_{\text{$n$ times}} * \delta_x.
\]
Recall that we say $\mu$ has a \emph{finite exponential moment} if there is some $\kappa > 0$ such that
\[
\int \norm{\Ad(g)}^\kappa \dd \mu (g) < +\infty,
\]
where $\Ad \colon G \to \GL(\kg)$ is the adjoint representation of $G$ on its Lie algebra $\kg$ and $\norm{\cdot}$ is any norm on $\kg$.
We will write $\Gamma_{\!\mu}\subseteq G$ to denote the  subgroup  generated by $\supp\mu$, the support of $\mu$. 
We say $\Gamma_{\!\mu}$ is Zariski dense in $G$ if $\Ad (\Gamma_{\!\mu})$ is Zariski-dense in $\Ad(G)$;  we will often simply say that \emph{$\mu$ is Zariski-dense} (for short, as in \cite{BQ1}).

\bfparagraph{Equidistribution in law.} We start with the main qualitative output of the paper: the $n$-step distribution of a Zariski-dense random walk on $X$ equidistributes toward the Haar measure unless the starting point is trapped in a finite invariant set. 

\begin{thm}[Equidistribution in law]
\label{thm:equidis}
Let $X=G/\Lambda$ be a homogeneous space where $G$ is a connected real linear group with Lie algebra isomorphic to $\mathfrak{so}(2,1)$ or $\mathfrak{so}(3,1)$ and $\Lambda$ is a lattice in $G$. 
Let $\mu$ be a probability measure on $G$ having a finite exponential moment and whose support generates a Zariski-dense subgroup $\Gamma_{\!\mu}$ of $G$.

For every $x \in X$, we have
\begin{equation}\label{eqth11}
\mu^{*n}* \delta_x \rightharpoonup^* m_X
\end{equation}
unless the orbit $\Gamma_{\!\mu} x$ is finite.
\end{thm}

In the above statement, the symbol $\rightharpoonup^*$ refers to weak-$*$ convergence of measures, in other terms \eqref{eqth11} means that $\lim_{n\to +\infty}\mu^{*n}* \delta_x (f)= m_X(f)$ for every continuous bounded function $f:X\rightarrow \R$.

Clearly, \eqref{eqth11}  is a dichotomy: if $\Gamma_{\!\mu} x$ is finite then the random walk cannot equidistribute in $X$.

It is well-known that \Cref{thm:equidis} implies the following rigidity statements, originally due to Benoist-Quint~\cite{BQ1}.
\begin{description}
\item[Classification of stationary measures] Every ergodic $\mu$-stationary probability measure on $X$ is either $m_X$ or the normalized counting  measure on a finite $\Gamma_{\!\mu}$-orbit.
\item[Classification of orbit closures] Every $\Gamma_{\!\mu}$-orbit is either dense or finite.
\item[Stiffness] Every $\mu$-stationary measure on $X$ is $\Gamma_{\!\mu}$-invariant.
\end{description}

Our proof of \Cref{thm:equidis} does not rely on the work of Benoist-Quint. Hence our method provides a new proof of the above rigidity statements in the setting covered by \Cref{thm:equidis}. An alternative proof of rigidity of orbit closures has also been obtained by Lee-Oh \cite{LO19} in the context of compact homogeneous spaces modeled on a rank one simple Lie group.

Benoist-Quint  \cite{BQ3}  proved that the convergence \eqref{eqth11} holds \emph{in Cesàro average}. \Cref{thm:equidis} is new in that it shows that $\mu^{*n}*\delta_{x} \rightharpoonup^* m_{X}$ holds without Cesàro average.

The question of removing the Cesàro average was asked by Benoist-Quint in \cite[Question 7]{BQ-survey12}.
Since \cite{BQ-survey12}, this has been managed under the assumption that the driving measure has non-mutually singular convolution powers (\cite{Benard23-equidmass} showing that this case boils down to the Ces\`aro average case),  for random walks on some nilmanifolds \cite{BFLM, He2020IJM, HS2022, HS2023, HLL_Aff, HLL_Nil}, see also \cite{KhalilLuethi}.
\begin{remark}
Contrary to \cite{Benard23-equidmass} we do not assume that $\mu$ is aperiodic. Periodicity may prevent equidistribution within \emph{finite} orbits as noted in \cite{Prohaska}.  
\end{remark}

\bfparagraph{Effective equidistribution.}
We will derive \Cref{thm:equidis} (in a nontrivial way) from the next quantitative statement, \Cref{thm:small-dim-equid-0}. It roughly says that the $\mu$-walk on $X$ equidistributes toward $m_{X}$ at an exponential rate provided the initial distribution of the walk has a positive dimension and is not too concentrated in the cusps. 

We fix a right-invariant Riemannian metric on $G$, which amounts to fixing an Euclidean norm on its Lie algebra $\kg$.
This metric induces a Riemannian metric on $X$, which we refer to as a \emph{quotient right $G$-invariant} Riemannian metric on $X$.
We denote by $\inj : X \to \R_{>0}$ the injectivity radius on $X$.

The rate of equidistribution will be expressed with respect to Lipschitz test functions. 
Let $\Lip(X)$ denote the space of bounded Lipschitz functions from $X$ to $\R$.
Endow it with its usual norm, that is,
for $f \in \Lip(X)$,
\[\norm{f}_{\Lip} \defeq  \norm{f}_\infty + \sup_{x\neq y \in X} \frac{\abs{f(x)-f(y)}}{\dist(x,y)}.\]
\begin{thm}[Effective equidistribution I]
\label{thm:small-dim-equid-0}
Let $X = G/\Lambda$, $\mu$ be as in \Cref{thm:equidis}.
Equip the homogeneous space $X=G/\Lambda$   with a quotient right $G$-invariant Riemannian metric.

Given $\kappa> 0$, there is $\eps > 0$ 
such that the following holds for all $\delta \in(0, \eps]$.
Let $\nu$ be a Borel probability measure on $X$, satisfying 
\begin{equation} 
\label{eq:cond nu kappa}
\nu (B_\rho(x)) \leq \rho^\kappa \text{ for all } x \in X, \rho \in [\delta, \delta^\eps].
\end{equation}
Then for all $n\geq |\log \delta|$ and all $f \in \Lip(X)$ with $\norm{f}_{\Lip} \leq 1$, 
one has
\[\abs{\mu^{*n}*\nu(f) -  m_{X}(f)} \leq \delta^\eps+ 2\nu\{\inj\leq \delta^\eps\}.\]
\end{thm}

We can also allow Hölder-continuous functions as test functions and state the convergence in terms of Wasserstein distance.
For the sake of conciseness, we stick to Lipschitz test functions in this introduction.  The more general statement is presented later as \Cref{small-dim-equid}.

Condition \eqref{eq:cond nu kappa} on the initial distribution $\nu$ should be interpreted as saying that 
$\nu$ has some positive dimension $\kappa > 0$ at scales above $\delta$.
For example, if $\nu$ is compactly supported and $\kappa$-Frostman, then one gets immediately the convergence $\mu^{*n} *\nu \rightharpoonup^* m_X$ with an exponential rate with respect to a Wasserstein distance.


\bfparagraph{Effective equidistribution under arithmetic assumptions.}
What is unsatisfactory about \Cref{thm:small-dim-equid-0} is that it is not ready to be applied to a random walk starting from a deterministic point.
However, under extra arithmetic assumptions, we are able to quantify the convergence $\mu^{*n}*\delta_x \rightharpoonup^* m_X$ in \Cref{thm:equidis}.

More precisely, in the theorem below, we further assume  that $\Lambda$ is an \emph{arithmetic} lattice and that 
 \emph{$\mu$ is algebraic with respect to $\Lambda$}.
The second assumption means that all elements in $\Ad(\Gamma_{\!\mu})$ and $\Ad(\Lambda)$ have algebraic entries with respect to some fixed basis of $\kg$. See also \Cref{def:alg} for an alternative characterization. 
Those arithmetic conditions are natural given the state of the art on effective equidistribution (e.g. \cite{BenoistSaxce, BFLM, LMW}). We comment on them further below. 


Now we consider only random walks starting from a deterministic point $x \in X$.
How fast the $\mu$-walk on $X$ equidistributes should depend on the starting point $x$, taking into account the natural obstructions that $x$ may be far away in a cusp of $X$ or close to a finite $\Gamma_{\!\mu}$-orbit of small cardinality.
To quantify such obstructions, we introduce $x_{0} \defeq \Lambda/\Lambda\in X$ which we see as a basepoint for $X$, as well as
\[W_{\!\mu, R} \defeq \set{x \in X : \sharp (\Gamma_{\!\mu}x) \leq R}\]
the set of points whose $\Gamma_{\!\mu}$-orbit is finite of cardinality at most $R$.



\begin{thm}[Effective equidistribution II] \label{main-thm}
Let $G$ be a connected real linear group with Lie algebra isomorphic to $\mathfrak{so}(2,1)$ or $\mathfrak{so}(3,1)$. Let 
$\Lambda$ be an arithmetic lattice in $G$, set $X=G/\Lambda$ equipped with a quotient right $G$-invariant  Riemannian metric. Let $\mu$ be a Zariski-dense finitely supported probability measure on $G$ which is algebraic with respect to $\Lambda$. 

There exists a constant $A>0$  such that for all  $x\in X$,  $n\in\N$, $R \geq 2$ and  $f\in \Lip(X)$, we have 
\[|\mu^{*n}*\delta_{x}(f)-m_{X}(f)| \leq R^{-1} \|f\|_{\Lip} \]
as soon as $n \geq A\log R + A \max \{|\log \dist(x,W_{\!\mu, R^A})|,\, \dist(x,x_{0})\}$.
\end{thm}

\begin{remark}
\begin{enumerate}
\item 
In the case where $W_{\!\mu, R}=\emptyset$, we use the convention that 
\[\max(|\log \dist(x,W_{\!\mu, R})|, \dist(x,x_{0}))=\dist(x,x_{0}).\]
\item In view of \Cref{comparison2}, the term $\dist(x,x_0)$ can be replaced by $\log \inj(x)$ where $\inj(x)$ denotes the injectivity radius of the point $x$.
Both measure how high in the cusp the starting point $x$ is.
\item In the case where $\Gamma_{\!\mu}$ does have a finite orbit on $X$, the condition that $\mu$ is algebraic with respect to $\Lambda$ can be dropped. Indeed, a certain conjugate $\mu'=\delta_{h}*\mu*\delta_{h^{-1}}$ ($h\in G$) of $\mu$ is then algebraic with respect to $\Lambda$ (see \Cref{reduction-commensurable} and the remark following it), and   it is enough to establish \Cref{main-thm} for $\mu'$.

\item The result remains valid, if one considers more generally bounded $\beta$-Hölder continuous test functions, $\beta \in {(0,1]}$, and the corresponding norm (see \Cref{small-dim-equid} and the plan of proof that follows it). The constant $A$ can be taken uniformly for values $\beta$ avoiding a given neighborhood of $0$. 

\item In the lower bound for $n$, the term $A |\log \dist(x, W_{\!\mu, R^A})|$ reflects the time needed for the walk to escape from a small neighborhood of $W_{\!\mu, R^A}$.
The proof shows that it can be  specified to   $\lambda^{-1}|\log \dist(x,W_{\!\mu, R^A})|$ where  $\lambda$ is any number in $(0, \lambda_{\mu})$  provided we also allow the various occurrences of $A$ to depend on $\lambda$.
Here  $\lambda_{\mu}>0$ denotes the top Lyapunov exponent of the adjoint random walk associated to $\mu$, i.e. $\lambda_{\mu} \defeq \lim_{n\to +\infty} n^{-1}\int_G \log \|\Ad g\| \dd \mu^{*n}(g)$.
Similarly, the term $A\cdot \dist(x,x_{0})$ reflects the time needed for the walk to come back to a large but fixed compact set.
We could specify it to $\lambda^{-1}C \dist(x,x_{0})$ for some constant $C$ depending only on the metric. 
Here $C$ is important as the overall term should be invariant by dilation of the metric.

\item Similar conclusions are known for affine random walks on some nilmanifolds \cite{BFLM, He2020IJM, HS2022, HS2023, HLL_Aff, HLL_Nil}  and also for some specific blockwise upper triangular walks on (some) arithmetic quotients of $\SL_{2}(\R)^2$, $\SL_{2}(\C)$ \cite{LMW}, $\SL_{3}(\R)$ \cite{Yang,LMWY}, or $\SL_{d}(\R)$  \cite{KhalilLuethi} provided extra assumptions on the starting point. More on how our result relates to and differs from these works was discussed in the preamble preceding \Cref{ss:mainresults}.


\end{enumerate}
\end{remark}


We now record a few corollaries of the above theorem.

First, remark that under these arithmetic assumptions, the qualitative equidistribution of \Cref{thm:equidis} follows immediately, since for $x$ having infinite orbit, the lower bound on $n$ prescribed by \Cref{main-thm} is finite.



We can also describe the set of starting points for which one has equidistribution with exponential rate. Given $D>1$, say $x\in X$ is \emph{$(\mu, D)$-Diophantine}  if for all $R>1$ with  $W_{\!\mu, R}\neq \emptyset$, one has
\[\dist(x, W_{\!\mu, R}) \geq \frac{1}{D} R^{-D}.\]
Observe this condition gets weaker as $D\to +\infty$. Say $x$ is \emph{$\mu$-Diophantine generic}  if it is $(\mu, D)$-Diophantine for some $D$. The set of $\mu$-Diophantine generic points $x\in X$ has full $m_X$\nobreakdash-measure. It is equal to $X$ when $\Gamma_{\!\mu}$ has no finite orbit.

\begin{corollary}[Points with exponential rate of equidistribution]
\label{cr:dioph}
In the setting of \Cref{main-thm}, let $x\in X$. The following are equivalent:
\begin{enumerate}
\item \label{it:exprate1}The point $x$ is $\mu$-Diophantine generic.
\item \label{it:exprate2} There exists $C, \theta>0$ such that for every $n\geq 1$, $f\in \Lip(X)$, 
\begin{equation} \label{exp-eq}
|\mu^{*n}*\delta_{x}(f)-m_{X}(f)| \leq \|f\|_{\Lip} C e^{-\theta n}.
\end{equation}

\end{enumerate}
Moreover, the constants $(C,\theta)$ can be chosen uniformly when $x$ varies in a compact subset and is $(\mu,D)$-Diophantine for a fixed $D$. 
\end{corollary}

Finally, we can deduce from \Cref{main-thm} that finite orbits of $\Gamma_{\!\mu}$ equidistribute toward the Haar measure on $X$ with polynomial rate as the cardinality of the orbit goes to infinity.

\begin{corollary}[Polynomial equidistribution of finite orbits]  \label{cor-orbites}
In the setting of \Cref{main-thm},
let $Y\subseteq X$ be a finite  $\Gamma_{\!\mu}$-orbit of cardinality $R$.
Let $m_{Y}$ denote the uniform probability measure on $Y$.
Then for all $f\in \Lip(X)$, one has
\[|m_{Y}(f)-m_{X}(f)| \leq \|f\|_{\Lip} C R^{-c}\]
where $C,c>0$ depend only  on $X$ and $\mu$.
\end{corollary}

Observe that, because of Lemmas \ref{finite-orbits} and \ref{DQ-sep}, the rate of equidistribution of finite orbits cannot be faster than polynomial.
 
  In the context 
where $\Gamma_{\mu}$ is a lattice, \Cref{cor-orbites}  follows from Maucourant-Gorodnik-Oh \cite[Corollary 3.31]{GMO08}, see also \cite{Clozel-Oh-Ullmo}. 

 \Cref{cor-orbites}  is an effective version of \cite[Corollary 1.8]{BQ3}, and a probabilistic analog of the polynomial equidistribution of large periodic unipotent orbits (which follows from \cite{LMW}). 

\bfparagraph{Effectiveness.}
The constant $\eps$ in \Cref{thm:small-dim-equid-0} and the constant $A$ in \Cref{main-thm} are \emph{effective} in the sense that by following the proof, one should be able to obtain an explicit formula for $\eps$ and $A$ in terms of standard numerical properties of $G, \Lambda, X, \mu$, as well as the additional parameter $\kappa$ for $\eps$. Those properties include in a meaningful way  the escape rate of $\check\mu$-walk\footnote{We denote by $\check\mu$  the image of $\mu$ by $g\mapsto g^{-1}$.} on $\kg$, the non-concentration estimates of the limiting stationary measure on the relevant Grassmanian,  the spectral gap of the Markov operator associated to $\mu$ and acting on $L^2(X, m_{X})$, and also for \Cref{main-thm}, the algebraic complexity (Mahler measure) of elements in the support of $\mu$.  

\bfparagraph{Further discussions.}
One can ask whether the conclusion of \Cref{main-thm} still holds without any arithmetic assumption. 
To shed some light on that question, it is interesting to compare our effective equidistribution statements to those concerning random walks on compact simple Lie groups. 

Recall that  Benoist-Saxcé~\cite{BenoistSaxce}, building upon Bourgain-Gamburd~\cite{BG_SU2,BG_SUd}, show that an aperiodic random walk on a compact simple Lie group enjoys exponential equidistribution toward the Haar measure (also referred to as spectral gap in this context) if and only if  it satisfies an almost Diophantine property. Such a statement is similar in spirit to \Cref{thm:small-dim-equid-0} where we show that exponential equidistribution holds under an assumption of a positive dimension for the initial distribution. 

Our proof of \Cref{main-thm} is a combination of \Cref{thm:small-dim-equid-0} with  \Cref{thm-positive-dimension}, the latter establishing the required positive dimension for $\nu = \mu^{*m} * \delta_x$ with controlled time $m$, under arithmetic assumptions. Here again, the analogy with \cite{BG_SU2,BG_SUd, BenoistSaxce} continues: in those papers the authors manage to check the almost Diophantine condition when the driving measure is supported on matrices with algebraic entries. 

Removing this algebraic assumption in the works \cite{BG_SU2,BG_SUd, BenoistSaxce}  is a well-known long-standing open problem.
We believe eliminating the arithmetic conditions in our context (\Cref{main-thm}) could be equally challenging.

\subsection{Ideas of the proof} \label{Sec-ideas-more}
We explain the ideas involved in the paper.
It is  an opportunity to introduce our main tool, a new result in incidence geometry, which we refer to as a 
\emph{multislicing} theorem.

\bfparagraph{Three phases.} Starting from a point, the random walk takes three phases to equidistribute.
The reader can easily notice the similarity with the structure of the proof in \cite{BG_Fp} and that in \cite{LMW_CRAD, Yang}.

\smallskip
In \emph{phase I}, the random walk gains some initial dimension, some small but positive dimension.
More precisely, we show that for any small enough $\rho\in(0, 1)$ and $x\in X$ with infinite $\Gamma_{\mu}$-orbit, for large  enough $n$, 
we have 
\begin{equation} \label{posdim-intro-eq}
\sup_{y\in X} \mu^{*n}*\delta_{x}(B_{\rho}(y))<\rho^\kappa,
\end{equation}
where $\kappa > 0$ is a small constant depending only on the ambient metric space $X$ and the driving measure $\mu$.

Without arithmeticity of the lattice and algebraicity of the coefficients of the matrices in $\supp \mu$, the lower bound on the time $n$ that guarantees \eqref{posdim-intro-eq}  is not controlled, see the proof of \Cref{thm:equidis} in \Cref{Sec-proof-main-statements}.  However, under the arithmetic conditions of \Cref{main-thm}, we show that for some large constant $A=A(X, \mu)>0$, the estimate \eqref{posdim-intro-eq} holds as soon as  $n \geq |\log \rho| + A \max \{|\log \dist(x,W_{\!\mu, \rho^{-A}})|,\, \dist(x,x_{0})\}$ (and even if $x$ has finite $\Gamma_{\mu}$-orbit). This is the content of  \Cref{thm-positive-dimension}. 

Note that such results can be seen as a closing lemma: given a starting point $x$, if most trajectories land at nearby points at a large time $n$, then $x$ must be very close to a finite orbit of the walk. 
  
The argument we have for this phase is valid for more general situations.
Namely, $G$ is allowed to be any semisimple real linear group without compact factors.


\smallskip
In \emph{phase II}, we start with a distribution having some initial dimension and by running the random walk, we show that the dimension increases until it reaches any prescribed number smaller than $\dim X$.
In this phase, ideas of incidence geometry are used.
More precisely, we need a generalization of discretized projections theorems which we briefly explain below. 

Our argument for this phase does not require the lattice to be arithmetic.
However, we ask $G$ to be isogenous to $\SO(2,1)$ or $\SO(3,1)$.

\smallskip
In \emph{phase III}, we get equidistribution starting from a measure of high dimension on $X$.
This is done using the spectral gap of the action of $\Gamma_{\!\mu}$ on $X$.
The argument for this phase is very general: $G$ is allowed to be any semisimple real linear group without compact factor,  arithmeticity is not required. 

\medskip

By combining phases II and III, we obtain \Cref{thm:small-dim-equid-0}, while \Cref{thm:equidis} and \Cref{main-thm} are proven by adding phase I to the picture.

\bfparagraph{Dimension increment.}
Let us take a moment to explain the idea behind $\text{phase}\,\text{II}$, as this part of our argument is completely new compared to previous works.

First, consider the case of a random walk on $X = \R^2/\Z^2$ driven by a probability measure $\mu$ on $\SL_2(\Z)$.
Assume that the support of $\mu$ is finite and generates a Zariski-dense subgroup. Denote by $\lambda_\mu>0$ the Lyapunov exponent of the $\mu$-walk on $\R^2$. 
We consider a measure $\nu$ on $X$ having a dimension $\alpha \in [\kappa,2 - \kappa]$ at all scales greater than a certain $\delta>0$.
The goal is to show that the $n$-step distribution starting from $\nu$, namely $\mu^{*n}*\nu$, has dimension $\alpha + \eps$ at scale $\rho = \delta^{1/2}$ for $n = \frac{\abs{\log \rho}}{\lambda_\mu}$, after possibly removing from  $\mu^{*n}*\nu$ a part of mass at most $\delta^\eps$ (ultimately negligible).
To this end, we need to bound
\[
\mu^{*n}*\nu(B_{\rho}(x)) = \int_{\SL_2(\Z)} \nu(g^{-1}B_\rho(x)) \dd \mu^{*n}(g)
\]
for $x \in X$.
By the large deviation estimates, a typical element $g$ sampled according to $\mu^{*n}$ has norm roughly $e^{\lambda_\mu n} = \rho^{-1}$.
Thus, the set $g^{-1}B_\rho(x)$ is roughly an ellipsoid with a major axis of length $1$ and a minor axis of length $\rho^2=\delta$,  the direction of the major axis being the stable direction $\theta_g$ of $g$.
We know by a result due to Guivarc'h~\cite{Guivarch} that the distribution of $\theta_g$ has a Hölder regularity.
Hence we may apply Bourgain's discretized projection theorem~\cite{Bourgain2010} to the family of orthogonal projections parallel to $(\theta_g)_{g \sim \mu^{*n}}$. Seen as a slicing theorem, this tells us that for most directions $\theta_g$, after removing from $\nu$ a small part of mass $\delta^\eps$, every fiber above a $\delta$-ball in $\theta_g^\perp$ has $\nu$-mass at most $\delta^{\frac{\alpha}{2}+\eps}= \rho^{\alpha + 2 \eps}$. As such a fiber is a rectangle of side length $1 \times \delta$ with the major axis directed by $\theta_g$, this gives exactly the upper bound we need for $\nu(g^{-1}B_\rho(x))$.
This argument requires to convert Bourgain's projection theorem into a slicing theorem. Such conversion is standard, and the heuristic goes as follows. The measure $\nu$ has dimension at least $\alpha$ at scale $\delta$.
We can pretend it is the uniform measure on a $\delta$-separated set $A$ of cardinality $\delta^{-\alpha}$.
The projection theorem says that a typical projection of $A$ has $\delta$-covering number $\geq \delta^{-\frac{\alpha}{2} - \eps}$, i.e. we need this number of $1 \times \delta$\nobreakdash-fibers to cover $A$.
Thus, a typical fiber contains at most $\frac{\delta^{-\alpha}}{\delta^{-\frac{\alpha}{2} - \eps}} = \delta^{-\frac{\alpha}{2} +\eps}$ points of $A$.
So the $\nu$\nobreakdash-measure of a typical fiber is at most $\frac{\delta^{-\frac{\alpha}{2} +\eps}}{\delta^{-\alpha}} = \delta^{\frac{\alpha}{2} + \eps}$.

\smallskip
Next, consider the case where $X = \R^3/\Z^3$ and $\mu$ is a probability measure on $\SL_3(\Z)$ whose support is finite and generates a Zariski-dense subgroup in $\SL_3$.
For this discussion, let us focus on the case where the Lyapunov spectrum of $\mu$ acting on $\R^3$ is of the form $\lambda_\mu > 0 > - \lambda_\mu$.
We may attempt to repeat the argument that we used for $\R^2/\Z^2$ with parameters $\rho = \delta^{1/2}$ and $n = \frac{\abs{\log \rho}}{\lambda_\mu}$.
However, we notice that sets of the form $g^{-1}B_\rho(x)$ are roughly rotated rectangles\footnote{Such rectangles are also called thin tubes.} of side lengths $1 \times \delta^{1/2} \times \delta$.
Thus, a discretized projection theorem is insufficient for us to conclude, as rectangles of this form are not $\delta$-neighborhoods of fibers of a projection.
So we need to upgrade the projection theorem to a slicing type theorem allowing for rectangles with \emph{more than $2$ different side lengths} like these. We will call it a multislicing theorem.
Put differently, a multislicing theorem is roughly an estimate on the number of incidences between balls and thin tubes.

Before we move on, let us remark that our argument applied to the case of the torus leads to a new (and shorter) proof of the effective equidistribution of Bourgain-Furman-Lindenstrauss-Mozes~\cite{BFLM} in the case where the linear part of the acting group is Zariski-dense in $\SL_d$. 



\smallskip

Finally, let us turn to the setting where $X = \SO(2,1)/\Lambda$.
Compared to the $3$\nobreakdash-dimensional torus case above, the sets of the form $g^{-1}B_\rho(x)$ are now distorted.
They are ellipsoids only in some charts and these charts depend on $g$ and $x$. 
Thus, we need to allow  \emph{nonlinear} rectangles in our multislicing theorem.

\bfparagraph{Multislicing.} Our multislicing theorem will be stated later as \Cref{thm:slicing} and \Cref{coro:slicing}.
It generalizes Bougain's discretized projection theorem \cite{Bourgain2010}, as well as its extensions to higher rank projections \cite{He2020JFG} and to nonlinear projections \cite{Shmerkin}.
Besides being the key to the dimension increment in our argument, it is interesting on its own right and we believe it will find other applications.

To give a taste of the result, we state the special case that is used for random walks on $G/\Lambda$ with $G = \SL_2(\R)$ and $\Lambda < G$ an arbitrary lattice. 

Let $(E,F,H)$ be the standard $\sl_2$\nobreakdash-triple of $\kg = \sl_2(\R)$. Endow $\kg$ with the Euclidean structure for which $(E, F, H)$ is orthonormal, and endow $G$ with the associated right invariant metric. 
Let $K$ denote $\SO(2) < G$.
For $\theta \in K$, consider the map $\psi_\theta \colon \kg \to G$ defined by
\[
rE + sH + tF \mapsto \theta \exp(rE) \exp(sH) \exp(tF).
\]



\begin{thm}
\label{thm:slicingSL2}
Given $\kappa > 0$, there exists $\eps , \delta_0 > 0$ depending only on  the parameter $\kappa$  such that the following holds for every $\delta \in (0,\delta_0]$.
Let $\sigma$ be a probability measure on $K$, let $\nu$ be a Borel measure on $B_{1}^G$. Assume  the following non-concentration properties:
\[
\forall \rho \in [\delta,\delta^\eps], \quad \sup_{\theta \in K} \sigma \bigl(B^K_\rho(\theta)\bigr) \leq \rho^\kappa,
\]
and for some $\alpha \in [\kappa, 1- \kappa]$,
\[
\forall \rho \in [\delta,\delta^\eps], \quad \sup_{ g \in G} \nu( B^G_\rho(g)) \leq \rho^{3\alpha}.
\]
Then there exists $\cE \subseteq K$ such that $\sigma(\cE)\leq \delta^{\eps}$ and for all $\theta\in K \setminus \cE$, there is a set  $A_{\theta} \subseteq G$ with  $\nu( G \setminus A_{\theta}) \leq \delta^{\eps}$ and such that
\[ \sup_{v \in B^\kg_{1}}\, \nu\bigl(A_\theta \cap  \psi_\theta(v+R)\bigr) \leq \vol(R)^{\alpha + \eps},\]
where  
\(
R = \set{r E + s H + t F : r \in [0,1], \, s \in [0,\delta^{1/2}], t \in [0,\delta]}.
\)
\end{thm}


The proof of \Cref{thm:slicingSL2} roughly proceeds in two steps: (1) an estimate for the covering number of a set $A$ by nonlinear rectangles having only two side lengths but that are both non-macroscopic (i.e. given a positive power of $\delta$), (2) a submodular inequality for covering numbers that allows to combine the estimates of the first step and that of Shmerkin \cite{Shmerkin} to deduce the theorem.

\bfparagraph{About the restriction on $G$.} 
The condition that $G$  is isogenous to $\SO(2,1)$ or $\SO(3,1)$ is used only in the bootstrap phase of the argument (previously referred to as Phase II). This restriction enters for the following reason. Our multislicing theorems \ref{thm:slicing}, \ref{coro:slicing} require a Bourgain-type non-concentration property for the random partial flag associated with the involved random box; see condition \ref{it:NCsigma} in \Cref{thm:slicing} and the discussion immediately following its statement. In the setting of the  $\mu$-random walk on $X$, we observed earlier that the random box takes the form of a translate $g^{-1} B_{\rho^{1/2}}x$, where $g$ is distributed according to $\mu^{*n}$ for some $n=n(\mu, \rho)$.  
When $G$ is isogenous to $\SO(2,1)$ or $\SO(3,1)$,  we are able to verify the required non-concentration estimate for this random box. By contrast, for more general groups such as $G=\SL_{3}(\R)$ or $G=\SO(7,1)$ the non-concentration conditions needed to apply Theorems \ref{thm:slicing}, \ref{coro:slicing} fail; see the discussion in \cite[\S1.2]{BH25}. In \cite{BH25}, we relax the non-concentration hypothesis and ultimately extend  the results of the present paper to arbitrary simple Lie groups.


\subsection{Organization of the paper}
In \Cref{sec:slicing}, we establish our multislicing estimates. In \Cref{sec-intro-posdim} we show under arithmetic assumptions that the $\mu$-walk on $X$ acquires a positive dimension, and does so at an explicit rate. In \Cref{Sec-bootstrap}, we prove \Cref{thm:small-dim-equid-0}. The bootstrap phase corresponds to \Cref{dimension-bootstrap} and the endgame phase is encapsulated in \Cref{endgame}. The bootstrap relies on \Cref{sec:slicing} but not on   \Cref{sec-intro-posdim}, the endgame is self-contained. In \Cref{Sec-proof-main-statements}, we prove \Cref{thm:equidis}, \Cref{main-thm},  \Cref{cr:dioph} and \Cref{cor-orbites}. We also add \Cref{ss:Appendix} to complete some proofs involved in \Cref{sec:slicing}.

\subsection{Conventions and notations} \label{conventions-notations}
The cardinality of a set $A$ is denoted by $\sharp A$.
The neutral element of a group will be denoted by $\Id$.
The differential of a differentiable map $f$ at a point $x$ will be denoted by $D_x f$.

A ball of radius $\rho > 0$ and centered at $x$ in a metric space $X$ is denoted by $B^X_\rho(x)$.
When the ambient space is unambiguous from the context, we drop the $X$ and simply write  $B_\rho(x)$.
When the space has a distinguished point, for example, zero vector $0$ of a linear space of neutral element $\Id$ of a group, $B_\rho$ denotes the ball centered at the distinguished point.

We use the Landau notation $O(\,\cdot\,)$ and the Vinogradov symbol $\ll$. Given $a,b>0$, we also write $a\simeq b$ for $a\ll b\ll a$. We also say that a statement involving $a,b$  is valid under the condition $a\lll b$ if it holds provided  $a\leq \eps b$ where $\eps>0$ is a small enough constant. The asymptotic notations $O(\,\cdot\,)$, $\ll$, $\simeq$, $\lll$  refer to constants that can depend on other parameters (e.g. parameters needed to formulate the framework of the ongoing section, or introduced in the course of a proof). The dependence in the parameters considered as framework will be implicit, it will not appear in the notation. The other dependencies will be indicated as subscripts. For instance, $a \lll_p b$ means that the constant $\eps$ above depends on some ambient framework and an additional parameter $p$. What we mean by the ambient framework will be specified at the beginning of each section.

\addtocontents{toc}{\protect\setcounter{tocdepth}{0}}

\subsection*{Acknowledgements}
We are grateful to  P\'eter Varj\'u for pointing out that \Cref{thm:small-dim-equid-0} should imply \Cref{thm:equidis} and sharing his intuition on the proof.  
We also thank Yves Benoist, Emmanuel Breuillard, Yang Cao, Elon Lindenstrauss, Amir Mohammadi, Nicolas de Saxcé for enlightening discussions. 

\addtocontents{toc}{\protect\setcounter{tocdepth}{2}}

\section{Multislicing theorem}
\label{sec:slicing}

Let $d\geq 2$. 
We prove various lower bounds for the covering number of a set in $\R^d$ by nonlinear rectangles.
We deduce that if a measure is $\alpha$-Frostman with respect to balls for some $0<\alpha<d$, then it is $(\alpha+\eps)$-Frostman with respect to such rectangles.
The novelty of these results is that, contrary to previous work, we allow rectangles that can have more than $2$ different side lengths, with potentially none of them being macroscopic.
\bigskip

Consider a family of differentiable maps $\phi_\theta : B^{\R^d}_1 \rightarrow \R^d$, indexed by a measurable space $\Theta$. We assume that for every $x\in B_{1}$, the differential  $D_{x}\phi_\theta \in \End(\R^d)$ varies  measurably in $\theta\in\Theta$.



Let $0 \subsetneq V_1 \subsetneq \dots \subsetneq V_m \subsetneq \R^d$ be a flag of subspaces in $\R^d$ of length $m \geq 1$,  together with a sequence $(r_i)_{1\leq i \leq m}$ of real numbers satisfying 
\[
0 \leq r_1 < r_2 < \dots < r_m <1.
\]
For $\delta > 0$, we write $R =R ((V_{i})_{1\leq i\leq m}, (r_{i})_{1\leq i\leq m}, \delta)= B^{V_1}_{\delta^{r_1}} + \dotsb + B^{V_m}_{\delta^{r_m}} + B^{\R^d}_\delta$. If $A \subset \R^d$ is a subset, we let $\cN_{R}(A)$ denote the smallest number of translates of $R$ needed to cover $A$. We also write $\vol$ for the Lebesgue measure on $\R^d$.

The theorem below aims at giving a lower bound on the quantity $\cN_{R}(\varphi_{\theta}A)$, that is, the covering number of $A$ by sets of the form $\varphi_{\theta}^{-1}(x+R)$ ($x\in \R^d$), for most $\theta \in \Theta$ when $\theta$ varies according to some distribution satisfying certain non-distortion and non-concentration properties.
We state the non-concentration property in terms of an angle function $\dang$ which we now define.  Given $k=1, \dots, d-1$, let $\Gr(\R^d, k)$ denote the Grasmannian of subspaces of dimension $k$ of $\R^d$, and for $U\in \Gr(\R^d, k)$, $W\in \Gr(\R^d, d-k)$, set 
\[\dang(U, W) \defeq \abs{\det(u_{1}, \dotsc, u_{k}, w_{1}, \dotsc, w_{d-k})},\]
where $(u_{1}, \dotsc, u_{k})$ and $(w_{1}, \dotsc, w_{d-k})$ stand for any orthonormal basis of $U$ and $W$, and the determinant is taken with respect to any orthonormal basis of $\R^d$. Note that  $\dang(U, W)$ can be interpreted (up to bounded multiplicative constant and up to a power depending on $d$) 
as the smallest angle between two lines respectively included in $U$ and $W$. 
 
\begin{thm} 
\label{thm:slicing}
Given $(\phi_\theta)_{\theta \in \Theta}$, $(V_i)_{1\leq i \leq m}$ and $(r_i)_{1 \leq i \leq m}$ as above, given
$\kappa > 0$, there exist   $\eps_{0}$, $\delta_{0}>0$ depending only on $d$, $(r_{i})_{i}$ and $\kappa$ such that   the following holds for every $\eps\in (0, \eps_{0}]$ and  $\delta \in (0, \delta_{0}]$.

Let $\sigma$ be a probability measure on $\Theta$ and $A$ a subset of $B^{\R^d}_1$ satisfying the following properties:  
\begin{enumerate}[label=(\roman*), series=cond]
\item \label{it:distortion} For $\sigma$-almost every  $\theta \in \Theta$, 
\begin{gather}
\label{eq:phiLipschitz}
\forall x,y \in A, \quad \delta^{\eps} \norm{x-y}  \leq \norm{\phi_\theta(x) - \phi_{\theta}(y)} \leq \delta^{-\eps}\norm{x-y},\\
\label{eq:phiC2bounded}
\forall x,y \in A, \quad \norm{\phi_\theta(x) - \phi_{\theta}(y) - D_x\phi_\theta (x-y) }  \leq \delta^{-\eps}\norm{x-y}^2.
\end{gather}
\item \label{it:NCsigma} $\forall i\in \{1, \dotsc ,m\}$, $\forall x\in A$, $\forall \rho \geq \delta$, $\forall W \in \Gr(\R^d, \codim V_{i})$,
\begin{equation*}
\sigma \set{\theta\in\Theta : \dang((D_x \varphi_{\theta})^{-1}V_{i}, W)\leq \rho} \leq \delta^{-\eps} \rho^\kappa.
\end{equation*}
\item \label{it:NCA} There is $\alpha \in [\kappa, 1- \kappa]$ such that $\forall \rho \geq \delta$,
\[
\max_{x \in \R^d} \cN_\delta(A \cap B_\rho(x)) \leq \delta^{-\eps} \rho^{\alpha d} \cN_\delta(A).
\]
\end{enumerate}
Then the exceptional set  
\begin{equation}
\label{eq:Eset}
\begin{split}
\cE \defeq \{\, \theta\in\Theta : \exists A' \subseteq A \,\,&\text{ with }\,\,\cN_{\delta}(A')\geq \delta^\eps \cN_{\delta}(A)\\
& \text{ and } \,\,\cN_{R}(\phi_\theta A') < \vol(R)^{- \alpha - \eps}\}
\end{split}
\end{equation}
satisfies $\sigma(\cE)\leq \delta^\eps$.
\end{thm}

Assumption~\ref{it:distortion} controls the distortion of the random charts $\phi_{\theta}$. Specifically, \eqref{eq:phiLipschitz} asks that $\phi_\theta$ is $\delta^{-\eps}$-bi-Lipschitz, while \eqref{eq:phiC2bounded}  requires in spirit that the $C^2$-norm of $\phi_\theta$ is bounded by $\delta^{-\eps}$. Assumption~\ref{it:distortion} is automatically satisfied in the case where the maps $\phi_\theta$ are isometries of $\R^d$. Note this scenario is helpful to keep in mind, even though we will need to consider maps $\phi_\theta$ which are not isometries (but have bounded distortion) in the course of the paper.

Assumption~\ref{it:NCsigma} is a non-concentration condition for the random charts $\phi_{\theta}$. 
When the $\varphi_{\theta}$'s are isometries, $\varphi_{\theta}^{-1}(R)$ is just a Euclidean rectangle with the same shape as $R$ but placed in a certain way in space depending on $\theta$. 
The non-concentration assumption then asks that the face of  $\varphi_{\theta}^{-1}(R)$ spanned by sides of length greater than a constant does not accumulate too close to a subspace of complementary dimension as $\theta$ varies with law $\sigma$.
For general $\varphi_{\theta}$'s, the non-concentration requirement is similar but formulated at an infinitesimal scale. 

Assumption~\ref{it:NCA} is a Frostman-type non-concentation condition on $A$.
It rules out the possibility that $A$ is a ball, in which case the conclusion may fail. 

Finally, the conclusion of the theorem states that for most realizations of the random parameter $\theta$, every large subset $A'$ of $A$ cannot be covered by fewer than $\vol(R)^{- \alpha - \eps}$ (nonlinear) rectangles of the form $\varphi_{\theta}^{-1}(x+R)$, where $x\in \R^d$.




\medskip
In  \Cref{thm:slicing}, \emph{the case of  where $m=1$ and $r_{1}=0$} is  simply a reformulation of previously known results, due to Bourgain, Shmerkin,  and the second named author.
More precisely, Bourgain \cite[Theorem 5]{Bourgain2010} first proved the case where, in addition to $m=1$ and $r_{1}=0$, the maps $\varphi_{\theta}$ are isometries and $\codim V_{1}=1$. This means that one considers covering numbers by Euclidean rectangles with one small side and all other sides being macroscopic.
Note that the result of \cite{Bourgain2010} is expressed in terms of covering numbers by small balls of the image by orthogonal projection parallel to the macroscopic sides of the rectangle (alternative point of view).
In \cite{He2020JFG}, the second named author treated the case of higher rank projections, i.e. $\codim V_{1} \geq 1$, which means the rectangles are still Euclidean but may have several small sides of the same length.
Finally, Shmerkin \cite{Shmerkin} generalized these works to the nonlinear setting, i.e. the $\varphi_{\theta}$ are no longer isometries but only satisfy condition~\ref{it:distortion}. To be precise, Bourgain's theorem extends to \cite[Theorem 1.7]{Shmerkin} and, as explained in \cite[\S 6.4]{Shmerkin}, Shmerkin's argument combined with \cite[Theorem 1]{He2020JFG} gives the higher rank case.


\emph{All other cases of  \Cref{thm:slicing} (i.e $m\geq 2$, or $m=1$ and $r_{1}>0$) are new}. Note that $m\geq 2$ means that the rectangles that are considered may have $3$ different side lengths or more, while $r_{1}>0$ indicates that all sides are non-macroscopic. The case  $m\geq 2$ and $r_1 = 0$ will be used later in the paper to prove effective equidistribution of random walks. The case  $r_1 > 0$  presents independent interest (even for $m=1$) even though it will not play a further role in the paper.

\bigskip
\Cref{thm:slicing} will be applied to random walks in the following equivalent form.

\begin{corollary} 
\label{coro:slicing}
Given $(\phi_\theta)_{\theta \in \Theta}$, $(V_i)_{1\leq i \leq m}$ and $(r_i)_{1 \leq i \leq m}$ as above, 
$\kappa > 0$, there exist  $\eps_{0}$, $\delta_{0}>0$ depending only on $d$, $(r_{i})_{i}$ and $\kappa$ such that the following holds for every $\eps\in (0, \eps_{0}]$ and  $\delta \in (0, \delta_{0}]$.

Let $\sigma$ be a probability measure on $\Theta$ and $\nu$ a Borel measure on $B^{\R^d}_1$ satisfying conditions \ref{it:distortion} and \ref{it:NCsigma} from \Cref{thm:slicing} with $A\defeq \supp \nu$, as well as the following non-concentration property:
\begin{enumerate}[resume*=cond]
\item \label{it:NCnu} There is $\alpha \in [\kappa, 1- \kappa]$ such that $\forall x \in \supp(\nu)$, $\forall \rho \geq \delta$,
\[
\nu(B_\rho(x)) \leq \delta^{-\eps} \rho^{\alpha d}.
\]
\end{enumerate}
Then there exists $\cE \subseteq \Theta$ such that $\sigma(\cE)\leq \delta^{\eps}$ and for all $\theta\in\Theta \setminus \cE$, there is a set  $A_{\theta} \subseteq \R^d$ with  $\nu(\R^d \setminus A_{\theta}) \leq \delta^{\eps}$ and such that
\[ \sup_{x \in \R^d}\, ({\phi_\theta}_\star\nu_{|A_\theta}) (x+R) \leq \vol(R)^{\alpha + \eps},\]
where $R = B^{V_1}_{\delta^{r_1}} + \dotsb + B^{V_m}_{\delta^{r_m}} + B^{\R^d}_\delta$ as in \Cref{thm:slicing}.
\end{corollary}
Note that though we do not require $\nu$ to be a probability measure, this statement is essentially about probability measures.
Indeed, on the one hand, the condition \ref{it:NCnu} implies that $\nu(\R^d) \leq \delta^{-\eps}$.
On the other hand, the conclusion is trivially true if $\nu(\R^d) \leq \delta^\eps$ as we can take $A_\theta = \emptyset$.
Thus, normalizing $\nu$ to a probability measure only changes everything by at most a factor of $\delta^{-\eps}$.
Using this observation, we see that if the statement holds with $2\eps$ instead of $\eps$ for all probability measures $\nu$ then it holds with $\eps$ for any $\nu$.

\subsection{Preliminaries}
\label{Sec-prepa}
We set up some notation and collect some useful lemmata.

\bfparagraph{Asymptotic notations.} 
The notations $O(\,\cdot\,)$, $\ll$, $\simeq$, $\lll$ (see \ref{conventions-notations}) refer implicitly to  constants that may depend on the dimension $d$ 
but nothing else. 
Additional dependencies are indicated as subscripts, and sometimes we also recall the dependence on $d$ for clarity.

\bfparagraph{Partitions.} 
Let $\cP$ and $\cQ$ denote partitions of $\R^d$, let $A$ be a subset of $\R^d$.

We write $\cP(A)$ the set of cells of $\cP$ that meet $A$, that is,
\[
\cP(A) \defeq \set{ P \in \cP : P \cap A \neq \emptyset}.
\]
We further write $\cP(\nu) \defeq \cP(\supp(\nu))$ for a measure $\nu$.

We set $\cP_{|A}$  the partition of $A$ obtained by restricting  $\cP$-cells to $A$. Sometimes, we may see it as a partition of $\R^d$ by adding to it the complement of $A$. 

We say the partition  \emph{$\cQ$ refines $\cP$}, and write $\cP \prec \cQ$,  if  for every $Q \in \cQ$, we have $\sharp \cP(Q) = 1$. Note that $\cP$ and $\cQ$ admit a coarsest common refinement, written $\cP \vee \cQ$, which is obtained by taking the intersections of $\cP$-cells and $\cQ$-cells.

We say $\cQ$ \emph{roughly refines} $\cP$ with parameter $L \geq 1$, and write $\cP \overprec{L} \cQ$, if 
\[
\max_{Q \in \cQ}\, \sharp \cP(Q) \leq L.
\]
For example, if $\phi_\theta$ and $A\subseteq B^{\R^d}_{1}$ satisfy \eqref{eq:phiLipschitz} and if $\cD$ is a partition of $\R^d$ by translates of a given cube, then
\begin{equation}
\label{eq:phiDrho}
(\phi_\theta^{-1} \cD)_{|A}  \overprec{O(\delta^{-d \eps})} \cD_{|A} \text{ and vice versa}.
\end{equation}
Clearly, the relation $\overprec{L}$ is transitive in the sense that $\cP \overprec{L} \cP'$ and $\cP' \overprec{L'} \cP''$ implies $\cP \overprec{L L'} \cP''$.
Moreover, it is compatible with taking common refinements, that is, if $\cP \overprec{L} \cQ$ and $\cP' \overprec{L'} \cQ'$ for some partitions, then $\cP \vee \cP' \overprec{LL'} \cQ \vee \cQ'$.

\bfparagraph{Rectangles.}
Extend the flag $(V_i)_{1 \leq i \leq m}$ by setting $V_0 = \{0\}$ and $V_{m+1} = \R^d$.
Without loss of generality we can assume that $V_i = \Span(e_{1}, \dotsc, e_{\dim V_i})$ for each $i = 1,\dotsc,m$, where $(e_1,\dotsc,e_d)$ denotes the standard basis of $\R^d$.
We will use the shorthand $j_i \defeq \dim V_i - \dim V_{i-1}$ for $i \in \{1,\dotsc, m+1\}$.

Let $\increasing_m$ denote the set of real $(m+1)$\nobreakdash-tuples $(r_1, \dotsc, r_{m+1})$ such that
$0 \leq r_1 \leq \dots \leq r_{m+1} \leq 1$.
For $\delta > 0$ and $\ur = (r_1, \dotsc, r_{m+1}) \in \increasing_m$, 
we set $D_\delta^{\ur} \subset \R^d$ to be the rectangle
\[
D_\delta^{\ur} \defeq {[0,2^{k_1})}^{j_1} \times \dots \times {[0,2^{k_{m+1}})}^{j_{m+1}} \subset \R^d,
\]
where for each $i$, $k_i \in \Z$ is the unique integer such that $2^{k_i - 1} < \delta^{r_i} \leq 2^{k_i}$.
Note that the shape
\(
R_\delta^{\ur} \defeq B^{V_1}_{\delta^{r_1}} + \dotsb + B^{V_{m+1}}_{\delta^{r_{m+1}}},
\)
is covered by a bounded number of translates of $D_\delta^{\ur}$ and vice versa.

Let $\cD_\delta^\ur$ denote the partition corresponding to the tiling of $\R^d$ by $D_\delta^\ur$ and its translates.
Set also $\cD_\delta \defeq \cD_\delta^{(1,\dotsc,1)}$ to be the dyadic cube partition of side length $\simeq \delta$.

For $\ur,\us \in \increasing_m$, set
\[\ur \wedge \us= (\min\{s_{i}, r_{i}\})_{1\leq i\leq m+1},\qquad \ur \vee \us = (\max\{s_{i}, r_{i}\})_{1\leq i\leq m+1}.\]
We write $\ur \leq \us$ if $r_i \leq s_i$ coordinatewise. 
Using this notation, we have that $\cD_\delta^\ur \prec \cD_\delta^\us$ whenever $\ur \leq \us$ and that $\cD_\delta^\ur \vee \cD_\delta^\us = \cD_\delta^{\ur \vee \us}$.

\bfparagraph{Covering numbers.}
For a partition $\cP$ and a subset $A$ we write $\cN_\cP(A) \defeq \sharp \cP(A)$.
For tilings by rectangles, we simply write  $\cN_\delta^{\ur}(A) \defeq \sharp\cD_\delta^{\ur}(A)$, and   $\cN_\delta(A)$ when $\ur=(1, \dots, 1)$. 
Note that up to a bounded factor, it is the number of translates of $R_\delta^\ur$ needed to cover $A$.

Obviously, if $\cP \overprec{L} \cQ$ for some parameter $L \geq 1$ then $\cN_\cP(A)\leq L \cN_\cQ(A)$.
For example, \eqref{eq:phiDrho} implies that for any subset $A \subset B_1^{\R^d}$ and any $\phi_\theta$ satisfying condition~\eqref{eq:phiLipschitz}, we have
\begin{equation}
\label{eq:Nrhophi}
\forall \rho > 0,\quad \delta^{d \eps} \cN_\rho(A) \ll  \cN_\rho(\phi_\theta A) \ll \delta^{- d \eps} \cN_\rho(A).
\end{equation}

The following lemma will allow us to restrict to $\cD_{\delta^{r_{1}}}$-cells when estimating covering numbers by non-linear rectangles.  
\begin{lemma}\label{lm:sumAcapQ}
Let $\delta, \eps > 0$, let $A \subseteq B_1^{\R^d}$ and $\theta \in \Theta$ satisfying \eqref{eq:phiLipschitz}, let $\ur = (r_1,\dotsc,r_{m+1}) \in \increasing_m$. We have
\[
\cN_\delta^\ur(\phi_\theta A) \gg \delta^{d \eps} \sum_{Q \in \cD_{\delta^{r_1}}} \cN_\delta^\ur(\phi_\theta (A \cap Q)).
\]
\end{lemma}
\begin{proof}
The left-hand side is the covering number of $A$ by the partition $\varphi_{\theta}^{-1}\cD_\delta^{\ur}$. The sum in the right-hand side is the covering number of $A$ by the partition $(\varphi_{\theta}^{-1}\cD_\delta^{\ur}) \vee \cD_{\delta^{r_1}}$. From \eqref{eq:phiDrho} and $(r_1,\dotsc,r_1) \leq \ur$, we have, restricted to $A$,
\[ \cD_{\delta^{r_1}} \overprec{O(\delta^{-d\eps})} \varphi_\theta^{-1} \cD_{\delta^{r_1}} \prec \varphi_{\theta}^{-1}\cD_\delta^{\ur},\]
 whence $(\varphi_{\theta}^{-1}\cD_\delta^{\ur}) \vee  \cD_{\delta^{r_1}} \overprec{O(\delta^{-d\eps})} \varphi_{\theta}^{-1}\cD_\delta^{\ur}$, and the claim follows. 
\end{proof}

\bfparagraph{Regularity.}
Let $\cP$ and $\cQ$ be two partitions of $\R^d$, with $\cP \prec \cQ$. 
We say a set $A\subseteq \R^d$ is \emph{regular} between $\cP \prec \cQ$ if
\[
\forall P \in \cP(A),\quad \cN_\cQ(A \cap P) = \frac{\cN_\cQ(A)}{\cN_\cP(A)}.
\]
Note that some authors use the word \emph{uniform} instead. Regularity means that the number of $\cQ$-cells meeting  $A$  inside a given $\cP$-cell does not depend on the $\cP$-cell (provided the latter does intersect $A$). 

We say $A$ is \emph{regular} with respect to a filtration of partitions $\cP_1 \prec \dots \prec \cP_n$ if for each $i \in \{1,\dotsc, n-1\}$, $A$ is regular between $\cP_i \prec \cP_{i+1}$.

The following observation will be used to guarantee that given a set $A$ that is regular between $\cP \prec \cQ$,  if a subset $A'$ meets many  $\cQ$-cells  of $\cQ(A)$ then it meets many $\cP$-cells of $\cP(A)$. 

\begin{lemma}\label{lm:reg and subset}
If $A$ is regular between $\cP \prec \cQ$ and $A' \subset A$ is a subset, then
\[
\frac{\cN_\cP(A')}{\cN_\cP(A)} \geq \frac{\cN_\cQ(A')}{\cN_\cQ(A)}.
\]
\end{lemma}
\begin{proof}
We have
\[
\cN_\cQ(A') = \sum_{P \in \cP(A')} \cN_\cQ(A' \cap P) \leq \sum_{P \in \cP(A')} \cN_\cQ(A \cap P) = \cN_\cP(A')\frac{\cN_\cQ(A)}{\cN_\cP(A)}.\qedhere
\]
\end{proof}

The following regularization process is due to Bourgain.
\begin{lemma}[Regularization] \label{tree-structure}
Let $n\geq 2$ and $\cP_1 \prec \dots \prec \cP_n$ be a filtration of partitions of $\R^d$.
Let $A\subseteq \R^d$. 
There exists a subset $A' \subset A$ which is regular with respect to $\cP_1 \prec \dots \prec \cP_n$ and satisfies
\[
\cN_{\cP_n}(A') \geq \frac{\cN_{\cP_n}(A)}{\prod_{i=2}^{n} 2(1+ \log_2 \max_{P \in \cP_{i-1}} \cN_{\cP_{i}}(P))}.
\]
Moreover, $A'$ can be chosen to be the intersection of $A$ with a union of $\cP_n$\nobreakdash-cells.
\end{lemma}
This result is well-known, we indicate a brief proof for the reader's convenience.
More can be found in \cite[\S 2]{Bourgain2010} or the survey by Shmerkin for proceedings of the ICM \cite[Lemma 2.2]{shmerkin21}.
\begin{proof}
Let $M_n \defeq \log_2 \max_{P \in \cP_{n-1}} \cN_{\cP_n}(P)$.
Partition $\cP_{n-1}(A)$ by putting in the same class $\mathcal{C}_{j}$ the cells $Q \in \cP_{n-1}(A)$ such that $\cN_{\cP_n}(A \cap Q) \in [2^j, 2^{j+1})$ where $0 \leq j < M_n$. 
By the pigeonhole principle, there is some $j_{0}$ such that $\cN_{\cP_n}(\cup_{Q \in \mathcal{C}_{j_{0}}} A\cap Q) \geq \cN_{\cP_n}(A) / (1+M_n)$.
Remove from $A$ all the $\cP_{n-1}$-cells outside of $\mathcal{C}_{j_{0}}$, and if needed, remove as well at most half of the $\cP_n$-cells  of each $A \cap Q$ for $Q\in \mathcal{C}_{j_{0}}$ to obtain a subset $A'\subseteq A$ that is regular from $\cP_{n-1}$ to $\cP_n$ and satisfies $2\cN_{\cP_n}(A')\geq \cN_{\cP_n}(A)/(1+M_n)$.
Trim $A'$ in the same manner to get regularity from  $\cP_{n-2}$ to $\cP_{n-1}$, and so on.

The ``moreover'' part follows by the construction of $A'$.
\end{proof}

\subsection{Submodularity for covering numbers}
The proof strategy of \Cref{thm:slicing} is an induction on the integer $m$ ruling the number of side lengths.
The key to performing the induction step is a submodular inequality for covering numbers. We present it in the following lemma. 


\begin{lemma}[Submodular inequality] \label{submod-cn}
Let $\cP,\cQ,\cR,\cS$ be partitions and $A$ a subset (of some ambient space).
Assume that $\cR= \cP \vee \cQ$, $\cS \prec \cP$ and $\cS \prec \cQ$.
Then for every $c >0$, there is a subset $A' \subset A$ such that $\cN_{\cR}(A') \geq (1-c) \cN_{\cR}(A)$ and
\begin{align}\label{submod-coveringnumber}
\cN_\cP(A) \cN_\cQ(A)\geq \frac{c^2}{4}  \cN_\cR(A) \cN_{\cS}(A').
\end{align}
\end{lemma}

Later we will apply this to tilings by rectangles with parallel sides: $\cP=\cD_\delta^{\ur}$, $\cQ=\cD_\delta^{\us}$, $\cR=\cD_\delta^{\ur\vee \us}$, $\cS=\cD_\delta^{\ur \wedge \us}$. 

\bigskip

\begin{remark}
\begin{enumerate} 
\item An analogous estimate for the entropy is well-known.
Namely, with the above assumptions, for any probability measure $\nu$, we have
\begin{equation*}
H(\nu,\cP) + H(\nu,\cQ) \geq H(\nu,\cR) + H(\nu,\cS). 
\end{equation*}

\item It is not possible to replace $A'$ by $A$ in \eqref{submod-coveringnumber}. To see why, consider the example where $d=3$ and  $A$ is the intersection of $(\Z e_{1}\oplus \Z e_{2}) \cup \Z e_{3}$ with $B(0, R)$ where $R>0$ is large. 
Let $\cP$ be the tiling by the rectangle $[0,R) \times [0,1) \times [0,1)$ and its translates.
Let $\cQ$ be that induced similarly by $[0,1) \times [0,R) \times [0,1)$.
Then $\cR=\cP \vee \cQ$ is the tiling induced by  $[0,1) \times [0,1) \times [0,1)$.
Finally let $\cS$ be the tiling by $[0,R) \times [0,R) \times [0,1)$ and its translates. These partitions satisfy the assumptions of \Cref{submod-cn}, but $\cN_\cP(A)\simeq R$, $\cN_\cQ(A) \simeq R$, $\cN_\cR(A) \simeq R^2$, $\cN_\cS(A)\simeq R$, whence the inequality \eqref{submod-coveringnumber} fails for $A'=A$.
\end{enumerate} 
\end{remark}

\begin{proof}[Proof of \Cref{submod-cn}]
Without loss of generality, we may assume $A$ has at most one element in each cell of $\cR$.
Let $\eta$ be the uniform probability measure on $A$.
As $\cP$ refines $\cS$, we have the equality $\sum_{C \in \cS(A)} \cN_\cP(A \cap C) = \cN_\cP(A)$, which can be rewritten as
\[
\sum_{C \in \cS(A)} \eta(C) \frac{\cN_\cP(A \cap C)}{\eta(C)} = \cN_\cP(A).
\]
Applying the Markov inequality, we obtain 
\[
\eta\Bigl( \bigcup \set{C \in \cS : \cN_\cP(A\cap C) > 2c^{-1} \eta(C) \cN_\cP(A)} \Bigr) < c/2. 
\]
The same holds for $\cQ$ in place of $\cP$.
Define $A'$ to be the intersection of $A$ with the union of $C \in \cS(A)$ such that
\[
\cN_\cP(A\cap C) \leq 2c^{-1} \eta(C) \cN_\cP(A) \text{ and } \cN_\cQ(A\cap C) \leq 2c^{-1} \eta(C) \cN_\cQ(A).
\]
Then $\eta(A') \geq 1 - c$ and for every $C \in \cS(A')$, using that $\cR= \cP \vee \cQ$,
\begin{align*}
\cN_{\cR}(A \cap C) &\leq \cN_\cP(A\cap C) \cN_\cQ(A\cap C) \\
& \leq 4c^{-2} \eta(C)^2 \cN_\cP(A) \cN_\cQ(A).
\end{align*}
Noting that $\eta(C) = \frac{\cN_{\cR}(A \cap C) }{ \cN_{\cR}(A) }$, we find
\[
\cN_{\cR}(A) \leq 4c^{-2} \eta(C) \cN_\cP(A) \cN_\cQ(A).
\]
Summing over $C\in \cS(A')$, we obtain the desired inequality.
\end{proof}

The following corollary will not be used in this paper but it presents independent interest. In a discrete context, we bound the cardinality of a finite set in terms of the cardinality of its projections. Given a subspace $V\subseteq \R^d$, we denote by $\pi_{V}$  the orthogonal projector to $V$. 

\begin{corollary}[Submodularity for projections] \label{cover-prop}
Let $Z\subseteq \R^d$ be a finite set. Let $V, W\subseteq \R^d$ be two subspaces. Let $c>0$. Then there exists $Z'\subseteq Z$ such that $\sharp \pi_{V+ W}(Z') \geq (1-c) \sharp \pi_{V+ W}(Z)$ and 
\[\sharp \pi_{V}(Z) \, \sharp \pi_{W}(Z)  \geq \frac{c^{2}}{4} \sharp \pi_{V+ W} (Z) \, \sharp \pi_{V\cap W}(Z').\]
\end{corollary}

The interpretation is that the left-hand side counts twice the coordinates in $V\cap W$ whence they need to appear twice in the lower bound.
This result complements the Bollobas-Thomason uniform cover theorem~\cite{BollobasThomason}, obtaining similar bounds when all coordinates are counted the same number of times.
It can also be used to simplify the proof of \cite[Proposition 34]{He2020JFG}. 

\begin{proof}[Proof of \Cref{cover-prop}] Write $V^\perp= E\oplus (V^\perp \cap W^\perp)$, $W^\perp=F\oplus (V^\perp \cap W^\perp)$ for some $E,F\in \Gr(\R^d)$. If $E$, $F$,  and $V^\perp \cap W^\perp$  are mutually orthogonal, then the result  follows by applying \Cref{submod-cn}  to the partitions $\pi_{S}^{-1}\cD^{S}_{\delta}$ where $S\in \{V,W,V+W,V\cap W\}$ and $\cD^{S}_{\delta}$ refers to a suitable partition of $S$ by $\delta$-cubes (with small $\delta$). The general case can then be reduced to the previous orthogonal orthogonal case by 
introducing $\varphi \in \GL_{d}(\R)$ such that $\varphi E$, $\varphi F$  and $\varphi(V^\perp \cap W^\perp)$ are mutually orthogonal, and observing that $\sharp \pi_{S}A=\sharp \pi_{S_{\varphi}}\varphi A$ where $S_{\varphi}:=(\varphi S^\perp)^\perp$.  
\end{proof}

\subsection{Subcritical estimates}

The next two subsections are dedicated to proving lower bounds on the covering number of a set $A$ seen in the chart $\varphi_{\theta}$ by rectangles of the form $R_\delta^{\ur} \defeq B^{V_1}_{\delta^{r_1}} + \dotsb + B^{V_{m+1}}_{\delta^{r_{m+1}}}$ where  $\ur \in \increasing_m$ and $\delta>0$.
The heuristic is that $\cN^{\ur}_{\delta}(\varphi_{\theta}A)$ should be at least greater than the geometric average $\left(\cN_{\delta^{r_{1}}}(A)^{ j_{1}}\dotsm \cN_{\delta^{r_{m+1}}}(A)^{ j_{m+1}}\right)^{1/d}$ where $j_{i}\defeq \dim V_{i}-\dim V_{i-1}$,  such lower bound  corresponding to the worst case scenario where $A$ is a ball.
The next proposition shows that this lower bound holds indeed up to a small loss, provided  $\theta$ does not belong to a small exceptional subset. The term \emph{subcritical} refers to this small loss.

\begin{proposition}\label{pr:subcritical}
Given $(\phi_\theta)_{\theta \in \Theta}$, $(V_i)_{1\leq i \leq m}$, $\ur = (r_1,\dotsc,r_{m+1}) \in \increasing_m$ and $\kappa > 0$, 
there exists $C > 1$ depending only on $d$, $\ur$, $\kappa$ such that the following holds for all $\eps \lll_{d, \ur} 1 $ and all $\delta \lll_{d,\ur,\kappa,\eps} 1$.

Let $\sigma$ be a probability measure on $\Theta$  and $A\subseteq B_1^{\R^d}$ a subset satisfying \ref{it:distortion},  \ref{it:NCsigma} and
\begin{enumerate}[resume*=cond]
\item \label{it:regular} The set $A$ is regular with respect to the filtration $\cD_{\delta^{r_1}} \prec \dots \prec \cD_{\delta^{r_{m+1}}}$.
\end{enumerate} 
Then the exceptional set  
\begin{equation}
\label{eq:excep-}
\begin{split}
\cE \defeq \Bigl\{\, \theta  \in \Theta : \exists A' \subseteq A \,\,&\text{ with }\,\,\cN_{\delta^{r_{m+1}}}(A')\geq \delta^{\eps} \cN_{\delta^{r_{m+1}}}(A)\\
& \text{ and } \,\,\cN_\delta^{\ur}(\phi_\theta A') < \delta^{C \eps \abs{\log \eps}} \prod_{i=1}^{m+1} \cN_{\delta^{r_i}}(A)^{j_i/d} \,\Bigr\}
\end{split}
\end{equation}
has measure $\sigma(\cE) \leq \delta^{\eps}$.
\end{proposition}

\begin{remark}
In the case where $r_{i}\in \{0,r_{m+1}\}$ for all $i$ (equivalently $\ur=(0, 1)$ or is constant) the regularity assumption \ref{it:regular} on $A$ can be removed, see \Cref{pr:CoarseNonProj}. In all other cases, it is \emph{necessary}. 
To see why, consider the example where  $d=2$,  $\phi_{\theta}$ is the rotation of angle $\theta$, and $\ur =(1/2, 1)$. 
Take $A=A_{1}\sqcup A_{2}$ where $A_{1}$ is $\delta^{1/2}$-separated with $\sharp A_{1} \simeq \delta^{-1/2}$, while $A_{2}$ is the intersection between $\delta \Z^2$ and a ball of radius $\delta^{1/2}$, in particular $\sharp A_{2} \simeq \delta^{-1}$.
Observe that $A$ is not regular between $\cD_{\delta^{1/2}} \prec \cD_{\delta}$. Moreover $\cN_{\delta^{1/2}}(A) \simeq \delta^{-1/2}$, $\cN_{\delta}(A) \simeq \delta^{-1}$, and $\cN_{\delta}^\ur (\phi_\theta A) \simeq \delta^{-1/2}$ uniformly in $\theta$. 
In view of these estimates, for any angle $\theta$ and  $\delta\lll 1$, one has
\[\cN_{\delta}^{\ur} (\phi_\theta A) \leq  \delta^{1/5} \left(\cN_{\delta^{1/2}}(A)\cN_{\delta}(A)\right)^{1/2} \]
thus forbidding the conclusion of \Cref{pr:subcritical}, regardless of $\sigma$. What makes it possible to remove the regularity assumption on $A$ when $\ur=(0, 1)$ is the regularization procedure of \Cref{tree-structure}. This procedure does not apply for other types of rectangles because it only preserves $\cN_{\delta}$ but not other $\cN_{\delta^r}$
\end{remark}

The case where $r_1 = \dots = r_{m+1}$ are all equal is trivial in view of \eqref{eq:phiLipschitz}.
Working with $\delta^{r_{m+1}}$ at the place of $\delta$, we can always assume $r_{m+1} = 1$.

The special case of $m = 1$ and $r_1 = 0$ is the subcritical counterpart of Shmerkin's nonlinear projection theorem. 
It can be shown by combining \cite[Proposition 29]{He2020JFG} and the linearization techniques of \cite{Shmerkin}.
Since its proof does not contain any new idea, we postpone it to \Cref{ss:Appendix}.
We will deduce the general case from this special case.

The next case is $m=1$ and $r_1 > 0$.

\begin{proof}[Proof of \Cref{pr:subcritical} when $m=1$ and $0 < r_1 < r_2 = 1$] 

Fix $\eps\in (0, 1/2]$, let $\delta>0$ be a parameter and write for a shorthand $\rho = \delta^{r_1}$, let  $\sigma$ and $A\subseteq B_{1}^{\R^d}$ as in the statement. We may assume that \eqref{eq:phiLipschitz} holds for \emph{every} $\theta\in \Theta$.

We first decompose into $ \cD_\rho$-cells. For each $\theta\in \Theta$, for $A'\subseteq A$,  \Cref{lm:sumAcapQ} implies
\begin{equation}\label{proj-0}
\cN_\delta^{\ur}(\phi_\theta A') \gg \delta^{d\eps} \sum_{Q\in \cD_\rho(A)} \cN_\delta^{\ur}(\phi_\theta(A'\cap Q)).
 \end{equation}
 
We then consider each cell individually. 
Apply the proposition with $(\varphi_{\theta})_{\theta\in \Theta}$,  $\{V_{1}\}$,  $(r_{1},r_{2})=(0,1)$ and write $C' = C'(d,\kappa) > 1$ the associated constant, assume $\delta\lll_{d,\kappa, 4\eps}1$ accordingly.
For $Q\in \cD_\rho(A)$, writing $A_{Q} \defeq A\cap Q$, we get some event $\cE_{Q}\subseteq \Theta$ with $\sigma(\cE_{Q})\leq \delta^{4\eps}$ and such that for $\theta \notin \cE_{Q}$,  any subset $A'_{Q}\subseteq A_{Q}$ with $\cN_{\delta} (A'_{Q})\geq \delta^{4\eps} \cN_{\delta} ( A_{Q})$ satisfies 
\[
 \cN^{\ur}_{\delta}(\varphi_{\theta}A'_{Q}) \geq  \cN^{(0, 1)}_{\delta}(\varphi_{\theta}A'_{Q}) \geq  \delta^{C' \eps \abs{\log \eps}} \cN_\delta(A_{Q})^{j_2/d}
\]
and hence with the regularity of $A$,
\begin{equation}\label{proj-1}
 \cN^{\ur}_{\delta}(\varphi_{\theta}A'_{Q}) \geq \delta^{C' \eps \abs{\log \eps}} \cN_{\rho}(A)^{-j_{2}/d}\cN_{\delta}(A)^{j_{2}/d}.
\end{equation}

In view of \eqref{proj-0} and \eqref{proj-1}, we define for $\theta\in \Theta$,
\[
\cQ_{\mathrm{bad}}(\theta) \defeq \set{Q \in \cD_{\rho}(A) : \theta \in \cE_{Q} }
\]
and for $A'\subseteq A$,
\[
\cQ_{\mathrm{large}}(A') \defeq \set{ Q \in \cD_\rho(A) : \cN_\delta(A'\cap Q) \geq \delta^{4\eps} \cN_\delta(A_{Q})}.
\]
It follows that
\begin{align*}
\cN_\delta^{\ur}(\phi_\theta A') 
 &\gg \delta^{d \eps}\sum_{Q\in \cQ_{\mathrm{large}}(A')\smallsetminus \cQ_{\mathrm{bad}}(\theta)} \cN_\delta^{\ur}(\phi_\theta(A'\cap Q))\\
 &\geq \sharp(\cQ_{\mathrm{large}}(A')\smallsetminus \cQ_{\mathrm{bad}}(\theta)) \, \delta^{(C' + d) \eps  \abs{\log \eps}} \cN_{\rho}(A)^{-j_{2}/d}\cN_{\delta}(A)^{j_{2}/d}. 
\end{align*}

To conclude, we just need to show that for many $\theta$ and every large subset $A' \subset A$, the set of cells $\cQ_{\mathrm{large}}(A')\smallsetminus \cQ_{\mathrm{bad}}(\theta)$ represents a large proportion of $\cD_{\rho}(A)$. 
More precisely, we claim that there is a set $\cE' \subset \Theta$ with $\sigma(\cE') \leq \delta^\eps$ such that for every $\theta \notin \cE'$ and every $A'\subseteq A$ with $\cN_{\delta}(A')\geq \delta^{\eps} \cN_{\delta}(A)$, we have 
\[
\sharp (\cQ_{\mathrm{large}}(A') \setminus \cQ_{\mathrm{bad}}(\theta)) \geq \delta^{3\eps} \cN_\rho(A).
\]
 This would imply for such $\theta$ and $A'$,
\[
\cN_\delta^{\ur}(\phi_\theta A') \geq \delta^{(C'+ d + 3)\eps  \abs{\log \eps}} \cN_\rho(A)^{j_1/d} \cN_\delta(A)^{j_2/d},
\]
thus finishing the proof of this special case.

We now show the claim. On the one hand, we set
\[
\cE' \defeq \set{\theta \in \Theta \,:\, \sharp  \cQ_{\mathrm{bad}}(\theta) \geq \delta^{3\eps} \cN_\rho(A)}.
\] 
 Fubini's theorem and the bound on $\max_{Q \in \cD_\rho(A)}\sigma(\cE_Q)\leq \delta^{4\eps}$ implies $\sigma(\cE')\leq \delta^\eps$.
On the other hand, let  $A'\subseteq A$ be a subset with $\cN_{\delta}(A')\geq \delta^{\eps} \cN_{\delta}(A)$.
Then noting that $\cN_\delta(A' \cap Q) \leq \cN_\delta(A \cap Q) = \cN_\delta(A) \cN_\rho(A)^{-1}$ for every $Q \in \cD_\rho(A)$, we deduce from the definition of $\cQ_{\mathrm{large}}(A')$ that
\begin{align*}
\delta^\eps \cN_\delta(A) \leq \cN_\delta(A') &= \sum\nolimits_{Q \in \cD_\rho(A)} \cN_\delta(A' \cap Q) \\
&\leq \sharp \cQ_{\mathrm{large}}(A') \cN_\delta(A) \cN_\rho(A)^{-1} + \delta^{4\eps} \cN_\delta(A),
\end{align*}
whence 
\(
\sharp \cQ_{\mathrm{large}}(A') \geq \delta^{2\eps} \cN_\rho(A)
\).
This justifies the claim and concludes the proof.
\end{proof}

The case $m=1$ being solved, it serves as the initialization in the induction argument below.
\begin{proof}[Proof of \Cref{pr:subcritical}, case $m \geq 2$.]
We proceed by induction on the number $m$, or equivalently the number of different values among $(r_i)_{1\leq i \leq m+1}$ since equality among $(r_i)$ amounts to remove a term in the flag $(V_i)$.
Thus, assume that $0 \leq r_1 < \dots < r_{m+1} = 1$. We may also suppose \eqref{eq:phiLipschitz} holds for \emph{every} parameter $\theta\in \Theta$.

Let $\theta\in \Theta$ and $A' \subseteq A$ satisfying $\cN_\delta (A') \geq \delta^{\eps} \cN_\delta (A)$.
The goal is to bound from below the quantity $\cN_\delta^\ur(\phi_\theta A')$ for $\theta$ outside of an exceptional set of small measure (and independent of $A'$). 
Throughout the proof, $\delta > 0$ is assumed to be small enough depending on $(d, \ur, \kappa, \eps)$, in particular $\delta^{-\eps}$ is larger than any quantity of the form $ O(\abs{\log \delta}^{O(1)})$.

Write $\rho = \delta^{r_2}$ as a shorthand.
Set
\begin{gather*}
\ut \defeq \ur \vee (r_2,\dotsc,r_2) = (r_2,r_2,r_3,\dotsc,r_{m+1}),\\
\us \defeq \ur \wedge (r_2,\dotsc,r_2) = (r_1,r_2,r_2,\dotsc,r_2)
\end{gather*}
so that $\cD_\delta^\ut \defeq \cD_\delta^\ur \vee \cD_\rho$ and both $\cD_\delta^\ur$ and $\cD_\rho$ refine $\cD_\delta^\us$.

By \Cref{tree-structure}, there exists a subset $A_1 \subseteq A'$ which is regular with respect to the filtration
$\phi_\theta^{-1} \cD_\rho \prec \phi_\theta^{-1}\cD_\delta^{\ut} \prec \phi_\theta^{-1}\cD_{\delta}$ while (recall \eqref{eq:Nrhophi})
\[
\delta^{-d \eps}\cN_\delta(A_1) \gg \cN_\delta(\phi_\theta A_1)  \geq \delta^\eps \cN_\delta(\phi_\theta A') \gg \delta^{(d + 1)\eps}  \cN_\delta(A'),
\]
hence 
\begin{equation}
\label{eq:A1geqA}
\cN_\delta(A_1) \geq \delta^{3 d \eps} \cN_\delta(A).
\end{equation}


By \eqref{eq:Nrhophi} and the submodular inequality from \Cref{submod-cn} applied with $\cP = \cD_\delta^\ur$, $\cQ = \cD_\rho$ , one has
\begin{equation}
\label{eq:submodA1A2}
\delta^{-d \eps}\cN_\delta^{\ur}(\phi_\theta A') \cN_\rho(A)
\gg
\cN_\delta^{\ur}(\phi_\theta A_1)  \cN_\rho (\phi_\theta A_1)
\gg
\cN_\delta^{\ut}(\phi_\theta A_1) \cN_\delta^{\us}(\phi_\theta A_2)
\end{equation}
for some subset $A_2 \subset A_1$ satisfying  
\(
\cN_\delta^{\ut}(\phi_\theta A_2) \gg \cN_\delta^{\ut}(\phi_\theta A_1)
\). We are going to apply the induction hypothesis to bound below each term in the right-hand side of \eqref{eq:submodA1A2}.  

Note that $A_{2}$ is also large in $A$ at scale $\rho$.
Indeed, by \Cref{lm:reg and subset} the regularity of $\phi_\theta A_1$ between $\cD_\rho \prec \cD_\delta^\ut$ implies that 
\[
\cN_\rho( \phi_\theta A_2) \geq \frac{\cN_\rho(\phi_\theta A_1)}{\cN_\delta^\ut(\phi_\theta A_1)} \cN_\delta^\ut( \phi_\theta A_2) \gg \cN_\rho(\phi_\theta A_1)
\]
and then $\cN_\rho(A_2) \gg \delta^{2 d \eps}\cN_\rho(A_1)$ in view of \eqref{eq:Nrhophi}.
Similarly, the regularity of $A$ between $\cD_\rho \prec \cD_\delta$ and \eqref{eq:A1geqA} imply $\cN_\rho(A_1) \geq \delta^{3d\eps} \cN_\rho(A)$.
Put together,
\begin{equation}
\label{eq:A2geqA}
\cN_\rho(A_2) \geq \delta^{6d\eps} \cN_\rho(A).
\end{equation}

By the induction hypothesis applied to $\ut$ and $3d\eps$ at the place of $\eps$, and recalling \eqref{eq:A1geqA}, there is $C_1 = C_1(d,\kappa, \ut) > 1$ such that
\[
\cN_\delta^{\ut}(\phi_\theta A_1)   
\geq 
\delta^{C_1\eps \abs{\log \eps}} \cN_{\delta^{r_{2}}}(A)^{(j_{1}+j_{2})/d} \cN_{\delta^{r_{3}}}(A)^{j_{3}/d} \dotsm \cN_{\delta^{r_{m+1}}}(A)^{j_{m+1}/d}
\]
whenever $\theta$ is outside a set independent of $A'$ and of $\sigma$-measure at most $\delta^{3d\eps}$.

By the induction hypothesis applied to $\us$ and $6d\eps$ at the place of $\eps$ and recalling \eqref{eq:A2geqA}, there is $C_2 = C_2(d,\kappa, \us) > 1$ such that
\[
\cN_\delta^{\us}(\phi_\theta A_2)
\geq \delta^{C_2\eps  \abs{\log \eps}}
\cN_{\delta^{r_{1}}}(A)^{j_{1}/d}  \cN_{\delta^{r_{2}}}(A)^{1-j_{1}/d}
\]
whenever $\theta$ is outside a set independent of $A'$ and of $\sigma$-measure at most $\delta^{6d\eps}$.

Plugging these back to \eqref{eq:submodA1A2} and simplifying by $\cN_\rho(A)$, we obtain
\[
\cN_\delta^{\ur}(\phi_\theta A') \gg \delta^{(C_1 +C_2 + d)\eps  \abs{\log \eps}} \prod_{i = 1}^{m+1} \cN_{\delta^{r_i}}(A)^{j_i/d}
\]
for $\theta$ outside a set independent of $A'$ and of $\sigma$-measure at most $\delta^{3d\eps} + \delta^{6d\eps} \leq \delta^\eps$.
\end{proof}

\subsection{Supercritical estimates}
Our next result is the \emph{supercritical} counterpart to \Cref{pr:subcritical}.
This means the output is a small gain on the heuristic 
$\cN_\delta^{\ur}(\phi_\theta A) \geq \prod_{i=1}^{m+1} \cN_{\delta^{r_i}}(A)^{j_i/d}$ instead of a small loss.
The price to pay is that we need a non-concentration assumption on the set $A$.

\begin{proposition}\label{pr:supercritical}
Given $(\phi_\theta)_{\theta\in\Theta}$, $(V_i)_{1\leq i \leq m}$, $\ur = (r_1,\dotsc,r_{m+1}) \in \increasing_m$, and $\kappa > 0$, 
the following holds for $\eps, \delta \lll_{d,\ur,\kappa} 1$. 

Let $\sigma$ be a probability measure on $\Theta$ and $A\subseteq B_1^{\R^d}$ a subset satisfying in addition to \ref{it:distortion}, \ref{it:NCsigma} and \ref{it:regular} the following single scale non-concentration condition.
\begin{enumerate}[resume*=cond]
\item \label{it:NCAr2} there exists $j \in \{1, \dots, m\}$ such that, writing 
\[\rho = \delta^{r_{j+1}} \cN_{\delta^{r_{j+1}}}(A)^{1/d} \cN_{\delta^{r_j}}(A)^{-1/d},\] we have
\[
\max_{x\in \R^d} \cN_{\delta^{r_{j+1}}}(A \cap B_{\rho}(x)) \leq \delta^\kappa \cN_{\delta^{r_{j+1}}}(A)\cN_{\delta^{r_j}}(A)^{-1}.
\]
\end{enumerate} 
Then the exceptional set 
\begin{equation*}
\begin{split}
\cE \defeq \Bigl\{\, \theta  \in \Theta : \exists A' \subseteq A \,\,&\text{ with }\,\,\cN_{\delta^{r_{m+1}}}(A')\geq \delta^{\eps} \cN_{\delta^{r_{m+1}}}(A)\\
& \text{ and } \,\,\cN_\delta^{\ur}(\phi_\theta A') < \delta^{-\eps} \prod_{i=1}^{m+1} \cN_{\delta^{r_i}}(A)^{j_i/d} \,\Bigr\}
\end{split}
\end{equation*}
has measure $\sigma(\cE) \leq \delta^{\eps}$.
\end{proposition}

\begin{remark}
The condition \ref{it:NCAr2} requires some explanation.
\begin{enumerate}
\item \label{comment-a} By the regularity assumption \ref{it:regular}, the ratio $\cN_{\delta^{r_{j+1}}}(A)\cN_{\delta^{r_{j}}}(A)^{-1}$ corresponds to the covering number by $\cD_{\delta^{r_{j+1}}}$ of $A\cap Q$ for an arbitrary  cell  $Q\in \cD_{\delta^{r_{j}}}(A)$. Hence $B_{\rho}(x)$ represents a ball of volume close to the volume of the $\delta^{r_{j+1}}$-neighborhood of $A\cap Q$. In sum \ref{it:NCAr2} asks that $A\cap Q$ is not concentrated in a small number of balls, and $\kappa$ rules the quality of the non-concentration. 

\item The left hand side being at least $1$, the condition~\ref{it:NCAr2} implies $r_{j+1} - r_{j} \gg \kappa$.
In particular $r_{j} < r_{j+1}$.
This is reasonable since the conclusion clearly fails when $r_1 = \dots = r_{m+1}$.
\item The condition~\ref{it:NCAr2} also implies
a conditional non-degeneracy assumption (for maybe a different $\kappa > 0$)
\[
\delta^{-\kappa(r_{j+1}-r_{j})} \leq \cN_{\delta^{r_{j+1}}}(A)\cN_{\delta^{r_{j}}}(A)^{-1} \leq \delta^{-(d-\kappa)(r_{j+1}-r_{j})}.
\]
\item The condition~\ref{it:NCAr2} is weaker than the combination of the non-degeneracy condition above with the following Frostman-type non-concentration.
\[
\forall \rho \geq \delta^{r_{j+1}-r_{j}},\quad \max_{x\in \R^d} \cN_{\delta^{r_{j+1}}}(A \cap B_{\delta^{r_{j}}\rho}(x)) \leq \rho^\kappa \cN_{\delta^{r_{j+1}}}(A)\cN_{\delta^{r_{j}}}(A)^{-1}.
\]
The projection theorems in \cite{Bourgain2010} and in \cite{He2020JFG} are stated with such combination (and $m=1$, $r_{1}=0$).
In~\cite{Shmerkin}, Shmerkin showed that a single-scale non-concentration condition is enough (and necessary). See \cite[Remark 1.8]{Shmerkin} for a discussion.
For the discretized sum-product theorem, the single-scale non-concentration condition is put forward by Guth-Katz-Zahl~\cite{GKZ}.


\item 
It is important that the non-concentration assumption \ref{it:NCAr2} on $A$ is conditional in the sense explained in \ref{comment-a}. For instance, assume $d=2$ and $m=1$. A non-concentration condition with respect to $\cD_{\delta^{r_{2}}}$ that would only be global and not conditional to each cell of $\cD_{\delta^{r_{1}}}$ would allow for $\cN_{\delta^{r_2}}(A\cap Q) = 1$ uniformly in $Q \in \cD_{\delta^{r_1}}(A)$.  But then $\cN_\delta^{\ur}(\phi_\theta A) \simeq \cN_{\delta^{r_{1}}}(A)$ and $\cN_{\delta^{r_2}}(A) = \cN_{\delta^{r_1}}(A)$, so there is no hope of dimensional gain. Similar problems arise in case of relative high dimension $\cN_{\delta^{r_2}}(A\cap Q)\simeq\delta^{-2(r_{2}-r_{1})}$. 
\end{enumerate} 

\end{remark} 

The proof of this proposition follows the same pattern as that of \Cref{pr:subcritical}.
The base case where $m = 1$ and $r_1 = 0$ is Shmerkin's nonlinear projection theorem, see\footnote{Shmerkin's statement is in fact slightly more restrictive as it requires $C^2$-differentiability for the random charts and bounded distortion, the latter meaning the scalars $\delta^{-\eps}, \delta^{-\eps}$ in \ref{it:distortion} are replaced by fixed arbitrary constants depending on which $\delta$ must be small enough. However,  Shmerkin's proof does yield the case where $m = 1$ and $r_1 = 0$ of \Cref{pr:supercritical}. See also \Cref{ss:Appendix}  where we use Shmerkin's argument to justify the base case of the subcritical estimate from \Cref{pr:subcritical} in our more general context.} \cite[Theorem 1.7 and \S 6.4]{Shmerkin}.
Next, we show the case of $m= 1$ and $r_1 > 0$.

\begin{proof}[Proof of \Cref{pr:supercritical} assuming that $m=1$]  We can assume $0 < r_1 < r_2 = 1$ and also suppose that \eqref{eq:phiLipschitz} holds for \emph{every} parameter $\theta\in \Theta$.

We will use the case of $m=1$ and $r_1 = 0$ as a black box.
To this end, consider $\ur' \defeq (0, r_{2})=(0,1) \in \increasing_1$ and assume $\eps, \delta \lll_{d, \kappa}1$   so that the conclusion of the proposition holds for the parameter $\ur'$ and with $(d+4) \eps$ at the place of $\eps$.
For each $Q\in \cD_{\delta^{r_1}}(A)$, we check that the proposition with parameter $\ur'$ can be applied to the set $A \cap Q$.
Indeed, the regularity condition~\ref{it:regular} with $\ur'$ holds trivially and \ref{it:NCAr2} with $\ur$ for $A$ combined with the regularity of $A$ implies \ref{it:NCAr2} with $\ur'$ for $A \cap Q$ (see comment \ref{comment-a}).
Thus, by the proposition with $m=1$ and parameter $\ur'$, 
we have $\sigma(\cE_Q) \leq \delta^{4\eps}$ where
\begin{equation*}
\begin{split}
\cE_Q \defeq \bigl\{\, \theta  \in \Theta : \exists A'_Q \subseteq A \cap Q \,\,&\text{ with }\,\,\cN_\delta(A'_Q)\geq \delta^{4\eps} \cN_\delta(A \cap Q)\\
& \text{ and } \,\,\cN_\delta^{\ur}(\phi_\theta A'_Q) < \delta^{-(d + 4) \eps} \cN_{\delta^{r_2}}(A \cap Q)^{j_2/d} \,\bigr\}.
\end{split}
\end{equation*}


The rest of the proof follows essentially the same lines as that of \Cref{pr:subcritical} case $m=1$.
Define $\cQ_{\mathrm{bad}}(\theta) \subset \cD_{\delta^{r_1}}(A) $, $\cE' \subset \Theta$ and then $\cQ_{\mathrm{large}}(A') \subset \cD_{\delta^{r_1}}(A)$ in the same way.
Then $\sigma(\cE') \leq \delta^\eps$ and for any $\theta \in \Theta \setminus \cE'$ and any $A' \subset A$ with $\cN_{\delta^{r_2}}(A') \geq \delta^\eps \cN_{\delta^{r_2}}(A)$, we have
\(
\sharp \bigl(\cQ_{\mathrm{large}}(A') \setminus \cQ_{\mathrm{bad}}(\theta) \bigr) \geq \delta^{3\eps} \cN_{\delta^{r_1}}(A)
\).
Then, by \Cref{lm:sumAcapQ} and the regularity of $A$,
\begin{align*}
\cN_\delta^{\ur}(\phi_\theta A') &\gg \delta^{d\eps}  \sum_{Q\in \cQ_{\mathrm{large}}(A') \setminus \cQ_{\mathrm{bad}}(\theta)} \cN_\delta^{\ur}(\phi_\theta(A'\cap Q))\\
&\geq \delta^{-4 \eps} \sharp \bigl(\cQ_{\mathrm{large}}(A') \setminus \cQ_{\mathrm{bad}}(\theta) \bigr)  \cN_{\delta^{r_2}}(A)^{j_2/d} \cN_{\delta^{r_1}}(A)^{j_2/d}\\
&\geq \delta^{-\eps}\cN_{\delta^{r_1}}(A)^{j_1/d}\cN_{\delta^{r_2}}(A)^{j_2/d}.
\end{align*}
This finishes the proof of the case $m=1$.
\end{proof}

For $m\geq 2$ and $j=1$ we use the submodular inequality again to reduce to a $m=1$ situation and a $m - 1$ situation. 
This time we do not rely on an induction, as we just use \Cref{pr:subcritical} in the $m-1$ situation.
This justifies why our assumption \ref{it:NCAr2} only concerns one scale transition (from $\delta^{r_{1}}$ to $\delta^{r_{2}}$) and not all of them. 

\begin{proof}[Proof of \Cref{pr:supercritical} for $m \geq 2$ and $j=1$] 
We may  suppose \eqref{eq:phiLipschitz} holds for \emph{every} parameter in $\Theta$.
Let $\theta \in \Theta$ and $A' \subseteq A$ a subset satisfying $\cN_{\delta^{r_{m+1}}}(A') \geq \delta^{\eps} \cN_{\delta^{r_{m+1}}}(A)$.

The same argument based on \Cref{tree-structure} and \Cref{submod-cn} as in the proof of \Cref{pr:subcritical} shows that
\begin{equation}
\label{eq:submodA1A2again}
\cN_\delta^{\ur}(\phi_\theta A') \cN_{\delta^{r_2}}(A)
\gg
\delta^{d\eps}\cN_\delta^{\ut}(\phi_\theta A_1) \cN_\delta^{\us}(\phi_\theta A_2)
\end{equation}
for 
\(\ut \defeq \ur \vee (r_2,\dotsc,r_2) = (r_2,r_2,r_3,\dotsc,r_{m+1})\),
\(\us \defeq \ur \wedge (r_2,\dotsc,r_2) = (r_1,r_2,r_2,\dotsc,r_2)\)
and some subsets $A_2 \subset A_1 \subset A'$ satisfying 
\begin{equation*}
\cN_{\delta^{r_{m+1}}}(A_1) \geq \delta^{3d\eps} \cN_{\delta^{r_{m+1}}}(A) \text{ and } \cN_{\delta^{r_2}}(A_2) \geq \delta^{6d\eps} \cN_{\delta^{r_2}}(A).
\end{equation*}

By \Cref{pr:subcritical} applied to $\ut$ and $3d\eps$ at the place of $\eps$, there is $C = C(d,\kappa, \ut) \geq 2$ such that for $\delta\lll_{d,\ur, \kappa, \eps}1$, 
\[
\cN_\delta^{\ut}(\phi_\theta A_1)   
\geq 
\delta^{C \eps \abs{\log \eps}} \cN_{\delta^{r_{2}}}(A)^{(j_{1}+j_{2})/d} \cN_{\delta^{r_{3}}}(A)^{j_{3}/d} \dotsm \cN_{\delta^{r_{m+1}}}(A)^{j_{m+1}/d}
\]
whenever $\theta$ is outside a set $\cE_1$ of measure $\sigma(\cE_1) \leq \delta^{3d\eps}$.

Note that the vector $\us$ consists of only two distinct values $r_1 < r_2$. This situation reduces to $m=1$ with the flag $(V_1)$. 
Applying the special case proved above with $\max\{2C \eps \abs{\log \eps},6d\eps\}$ at the place of $\eps$, we have for some $\eps =\eps(d,\kappa, C, \us) \in (0, 1/2]$ and for all $\delta\lll_{d,\ur, \kappa}1$, 
\[
\cN_\delta^{\us}(\phi_\theta A_2) 
\geq \delta^{- 2C \eps \abs{\log \eps}}
\cN_{\delta^{r_{1}}}(A)^{j_{1}/d}  \cN_{\delta^{r_{2}}}(A)^{1-j_1/d}
\]
whenever $\theta$ is outside a set $\cE_2$ of measure $\sigma(\cE_2) \leq \delta^{6d\eps}$.

Plugging these back to \eqref{eq:submodA1A2again} and simplifying by $\cN_{\delta^{r_2}}(A)$ shows the desired lower bound for 
$\cN_\delta^\ur(\phi_\theta A')$ for all $\theta \notin \cE_1 \cup \cE_2$ and all large $A' \subset A$. 
\end{proof}

\begin{proof}[Proof of \Cref{pr:supercritical} for $m \geq 2$ and $j\geq 2$]
This final case follows from the previous one, with the same method of proof: cut the sequence $\ur$ at $r_{j}$, i.e. consider
\(\ut \defeq \ur \vee (r_j,\dotsc,r_j) = (r_j, \dots,r_j,r_{j+1}\dotsc,r_{m+1})\),
\(\us \defeq \ur \wedge (r_j,\dotsc,r_j) = (r_1,r_2,\dots, r_j,\dotsc,r_j)\), and apply the case $m\geq 2$, $j=1$  of  \Cref{pr:supercritical}  to $\ut$, and  \Cref{pr:subcritical} to $\us$ (the order is reversed this time).  
\end{proof}

\subsection{Proof of Theorem \ref{thm:slicing}}
In this subsection, we proceed to the proof of \Cref{thm:slicing}.
Extend the vector $(r_i)_{1 \leq i \leq m}$ by $r_{m+1} = 1$ to form $\ur = (r_i)_{1 \leq i \leq m+1} \in \increasing_m$. 
Note that for the set $R$ in the statements, $\vol(R)$ can be expressed explicitly: 
\[\vol(R) \simeq \delta^{\sum_{i=1}^{m+1} r_i j_i}.\]
Recall also that $\cN_R$ and $\cN_\delta^\ur$ are comparable, i.e. $\cN_R\simeq \cN_\delta^\ur$.
Finally, observe that if the statement in  \Cref{thm:slicing} holds for some $\eps > 0$ and $\delta > 0$ then it holds for $\eps'$ and $\delta$ for all $\eps' \in (0, \eps)$ since making $\eps$ smaller only strengthens the assumptions and weakens the conclusion.
Therefore, to prove \Cref{thm:slicing}, it is enough to show the statement holds for some $\eps = \eps(d, \kappa, \ur)> 0$ and all $\delta>0$ smaller than a constant depending on $d$, $\ur$, $\kappa$ and $\eps$.

We start the proof with the special case where $A$ satisfies an extra regularity assumption. 

\begin{proof}[Proof of \Cref{thm:slicing} when $A$ is regular.]
Let $\eps, \delta \in (0, 1/2]$ be parameters. Let $\sigma, A$ be as in the statement of the theorem.
Assume that $A$ satisfies the condition~\ref{it:regular}, that is, $A$ is regular with respect to $\cD_{\delta^{r_1}} \prec \dots \prec \cD_{\delta^{r_{m+1}}}$.

Note that the condition \ref{it:NCA} applied with $\rho = \delta^{r_i}$ gives 
\begin{equation}\label{eq-dadte}
\forall i \in \{ 1,\dotsc, m+1 \},\quad 
\cN_{\delta^{r_i}}(A) \geq  \delta^{-r_i \alpha d+\eps}.
\end{equation}

Let $C = C(d,\kappa, \ur) > 0$ be the constant given by \Cref{pr:subcritical}.
If there is $i \in \{1,\dotsc,m+1\}$ such that
\[
\cN_{\delta^{r_i}}(A) \geq  \delta^{-r_i \alpha d - (C + 1)\frac{d}{j_{i}}\eps \abs{\log \eps}}.
\]
Then combined with \eqref{eq-dadte} we see immediately that the exceptional set $\cE$ of \eqref{eq:Eset} is contained in the exceptional set of \eqref{eq:excep-} (at least if $\eps$ is small enough in terms of $d$, and  $\delta\lll_{\eps}1$ as we may suppose).
This concludes the proof of $\sigma(\cE) \leq \delta^\eps$ in this case.

Therefore, we can assume that for every  $i \in \{1,\dotsc,m+1\}$,
\begin{equation}
\label{eq:Aleq}
\cN_{\delta^{r_i}}(A) \leq  \delta^{-r_i \alpha d - (C + 1)\frac{d}{j_{i}}\eps \abs{\log \eps}}.
\end{equation}
In this case, we use \Cref{pr:supercritical}.
To do so, we need to check the condition \ref{it:NCAr2}.  Let $j\in \{1, \dots, m\}$. 
The condition \ref{it:NCA} combined with \ref{it:regular} implies
for all $\rho \geq \delta^{r_{j+1}}$, 
\[\max_{x \in \R^d} \cN_{\delta^{r_{j+1}}}(A \cap B_\rho(x)) \ll \delta^{-\eps} \rho^{\alpha d} \cN_{\delta^{r_{j+1}}}(A).\]
In particular, specializing to the specific $\rho$ of \ref{it:NCAr2}, we get
\[
\max_{x \in \R^d} \cN_{\delta^{r_{j+1}}}(A \cap B_\rho(x)) \ll \bigl(\delta^{-\eps + r_{j+1}\alpha d} \cN_{\delta^{r_{j+1}}}(A)^{\alpha} \cN_{\delta^{r_{j}}}(A)^{1-\alpha}\bigr)  \cN_{\delta^{r_{j+1}}}(A) \cN_{\delta^{r_{j}}}(A)^{-1}.
\]
In view of \eqref{eq:Aleq},
\[\delta^{-\eps + r_{j+1}\alpha d} \cN_{\delta^{r_{j+1}}}(A)^{\alpha} \cN_{\delta^{r_{j}}}(A)^{1-\alpha} \leq \delta^{(r_{j+1} - r_{j})\alpha(1-\alpha)d - (C + 2)d\eps \abs{\log \eps}}.\]
The right hand side is $\leq \delta^{\kappa'}$ with $\kappa' \defeq \frac{1}{2}\min\{\kappa, (r_{j+1} - r_{j})\}^3$ provided  $\eps$ is small enough in terms of $C, d,\kappa'$.
Thus \ref{it:NCAr2} holds for $\kappa'$ at the place of $\kappa$ and \Cref{pr:supercritical} can be applied with $\kappa'$.
We obtain that for $\eps, \delta\lll_{d,\ur,\kappa}1$, the following exceptional set
\begin{equation*}
\begin{split}
\Bigl\{\, \theta  \in \Theta : \exists A' \subseteq A \,\,&\text{ with }\,\,\cN_{\delta}(A')\geq \delta^{\eps} \cN_{\delta}(A)\\
& \text{ and } \,\,\cN_\delta^{\ur}(\phi_\theta A') < \delta^{-(d + 2) \eps} \prod_{i=1}^{m+1} \cN_{\delta^{r_i}}(A)^{j_i/d} \,\Bigr\}
\end{split}
\end{equation*}
has $\sigma$-measure less than $\delta^\eps$.
But this set clearly contains $\cE$ of \eqref{eq:Eset} thanks to \eqref{eq-dadte}.
This finishes the proof of \Cref{thm:slicing} under extra regularity.
\end{proof}
To reduce to this special case, we use an exhaustion procedure similar to the one used in \cite[Proposition 25]{He2020JFG}.
\begin{proof}[Proof of \Cref{thm:slicing}]
Throughout this proof, $\delta > 0$ is assumed to be small enough in terms of $d,\ur,\kappa$ and also $\eps$ (as we may). In particular,  $\delta^{-\eps}$ is larger than all the quantity of the form $ O(\abs{\log \delta}^{O(1)})$ that will appear.

In this paragraph, by regular we mean regular with respect to the filtration $\cD_{\delta^{r_1}} \prec \dots \prec \cD_{\delta^{r_{m+1}}} = \cD_\delta$.
Using \Cref{tree-structure}, we can find a regular subset of $A_1 \subset A$ of size $\cN_\delta(A_1) \geq \delta^{\eps} \cN_\delta(A)$.
Then repeat to find a large regular subset $A_2 \subset A \setminus \bigcup \cD_\delta(A_1)$ and so on.
Repeat this until $\cN_\delta(A \setminus \bigcup_{k \in I} A_k) \leq \delta^{2\eps}\cN_{\delta}(A)$.
The family of regular subsets $(A_k)_{k \in I}$ then satisfies
that $\cD_\delta(A_k)$ are pairwise disjoint and 
\[
\forall k \in I,\quad \cN_\delta(A_k) \geq \delta^{3\eps}\cN_\delta(A).
\]

This combined with the condition~\ref{it:NCA} implies
\[
\max_{x \in \R^d} \cN_\delta(A_k \cap B_\rho(x)) \leq \delta^{-4\eps} \rho^{\alpha d} \cN_\delta(A_k).
\]
Thus, the special case proved above can be applied to the set $A_k$.
We conclude that for $\eps>0$ small enough in terms $d$, $\ur$, $\kappa$,
we have
\(\sigma(\cE_k) \leq \delta^{4\eps}\)
for 
\begin{equation*}
\begin{split}
\cE_k \defeq \{\, \theta\in\Theta : \exists A' \subseteq A_k \,\,&\text{ with }\,\,\cN_{\delta}(A')\geq \delta^{4\eps} \cN_{\delta}(A)\\
& \text{ and } \,\,\cN_{R}(\phi_\theta A') < \vol(R)^{- \alpha - 4\eps}\}.
\end{split}
\end{equation*}

We claim that the exceptional set $\cE$ of \eqref{eq:Eset} satisfies
\[
\cE \subset \bigcup_{J} \bigcap_{k \in J} \cE_k
\]
where the union is taken over all $J \subset I$ such that $\sum_{k \in J} \cN_\delta(A_k) \geq \delta^{2\eps} \cN_\delta(A)$.
This would finish the proof since by the Markov inequality (see \cite[Lemma 20]{He2020JFG}), we see that 
\[
\sigma(\cE) \leq \delta^{-2\eps} \max_{k \in I} \sigma(\cE_k) \leq \delta^\eps.
\]

Indeed, let $\theta \in \cE$.
By definition there exists $A' \subset A$ such that $\cN_\delta(A') \geq \delta^\eps \cN(A)$ and $\cN_\delta^\ur(\phi_\theta A') < \vol(R)^{- \alpha - \eps}$.
Consider
\[
I_{\mathrm{large}}(A') \defeq \set{ k \in I : \cN_\delta(A' \cap A_k) \geq \delta^{2\eps} \cN_\delta(A_k) },
\]
so that for every $k \in I_{\mathrm{large}}(A')$, the subset $A' \cap A_k \subset A_k$ witnesses $\theta \in \cE_k$.
Hence $\theta \in \bigcap_{k \in I_{\mathrm{large}}(A')} \cE_k$.
Moreover,
\begin{align*}
\cN_\delta(A') &\leq \cN_\delta\bigl(A \setminus \bigcup\nolimits_{k \in I} A_k\bigr) + \sum_{k \in I} \cN_\delta(A' \cap A_k) \\
&\leq \delta^{2\eps}\cN_\delta(A) + \sum_{k \in I_{\mathrm{large}}(A')} \cN_\delta(A_k) + \sum_{k \notin I_{\mathrm{large}}(A')}  \delta^{2\eps} \cN_\delta(A_k)\\
&\leq \sum_{k \in I_{\mathrm{large}}(A')} \cN_\delta(A_k) + 2\delta^{2\eps}\cN_\delta(A).
\end{align*}
Then $\sum_{k \in I_{\mathrm{large}}(A')} \cN_\delta(A_k) \geq \delta^{2\eps} \cN_\delta(A)$.
This finishes the proof of the claim.
\end{proof}

\subsection{Proof of Corollary \ref{coro:slicing}}
We proceed to the proof of \Cref{coro:slicing}. 
We established above a  \emph{lower} bound on the covering number of a set $A$ by some rectangles (as in a projection theorem) and we deduce an \emph{upper} bound on how much $A$ fills those rectangles (as in a slicing theorem). Such conversion, at least in spirit, dates back to Marstrand's seminal paper~\cite{Marstrand54}.
In a discretized setting, such a procedure is well-known. 
We include the details for the reader's convenience.

\begin{proof}[Proof of \Cref{coro:slicing}, special case]
Let $\eps, \delta >0$ be parameters. Let $\sigma, \nu$ be as in the statement of the corollary. 
We first deal with the case where, for $A \defeq \supp(\nu)$,
\[
\forall Q \in \cD_{\delta}(A),\quad \frac{\nu(\R^d)}{2\cN_\delta(A)} \leq  \nu(Q) \leq \frac{2\nu(\R^d)}{\cN_\delta(A)},
\]
so $\nu$ is essentially equally distributed among the $\cD_{\delta}$-cells. This assumption allows to convert non-concentration properties of the measure $\nu$ into similar properties for the support $A$.

As we remarked after the statement of \Cref{coro:slicing} and as a renormalization does not affect this extra assumption, we can assume without loss of generality that $\nu$ is a probability measure.

Then observe that for any union $B$ of $\cD_{\delta}$\nobreakdash-cells,
\[
\nu(B)\simeq {\cN_{\delta}(B)}/{\cN_{\delta}(A)}.
\]
Combined with the non-concentration assumption~\ref{it:NCnu} on $\nu$, this implies that for all $\rho\geq \delta$ and all $x\in \R^d$, 
\[
\cN_{\delta}(A\cap B_\rho(x)) \ll \nu(B_{2\rho}(x)) \cN_{\delta}(A) \ll \delta^{-\eps}\rho^{\alpha d}\cN_{\delta}(A).
\]
Here and below, assume $\delta > 0$ is small enough so that $\delta^{-\eps}$ is larger than any quantity of the form $ O(1)$.
Thus, the non-concentration assumption~\ref{it:NCA} holds for $A$ with $2\eps$ at the place of $\eps$.

Take $\eps>0$ small enough in terms of $d$, $(r_{i})_{i}$, $\kappa$ so that \Cref{thm:slicing} holds with $2 \eps$ at the place of $\eps$.
Thus, provided that $\delta \lll_{d, (r_{i})_{i},\kappa}1$, we have $\sigma(\cE) \leq \delta^{2\eps}$ for the exceptional set 
\begin{equation*}
\begin{split}
\cE \defeq \{\, \theta\in\Theta : \exists A' \subseteq A \,\,&\text{ with }\,\,\cN_{\delta}(A')\geq \delta^{2\eps} \cN_{\delta}(A)\\
& \text{ and } \,\,\cN_{R}(\phi_\theta A') < \vol(R)^{- \alpha - 2\eps}\}
\end{split}
\end{equation*}

Let $\theta\in\Theta \setminus \cE$.
Consider the part of $\varphi_{\theta}A$ covered by rectangles of large measure with respect to ${\phi_\theta}_\star \nu$. We will show that this part is small. 
Write $\ur = (r_1,\dotsc,r_m,1) \in \increasing_m$ and consider
\[
\cQ_{\mathrm{bad}}(\theta) \defeq \set{ P \in \cD_\delta^\ur : {\phi_\theta}_\star \nu(P) > \vol(R)^{\alpha +2\eps}}. 
\]
Since $\nu$ is a probability measure, $\sharp \cQ_{\mathrm{bad}}(\theta) < \vol(R)^{- \alpha - 2\eps}$.
Then letting
\[
A_\theta' = A \cap \bigcup_{P \in \cQ_{\mathrm{bad}}(\theta)} \phi_\theta^{-1} P,
\]
we have $\cN_\delta^\ur(\phi_\theta A_\theta') \leq \sharp \cQ_{\mathrm{bad}}(\theta) < \vol(R)^{- \alpha - 2\eps}$.
Recalling that $\theta \notin \cE$, we deduce that $\cN_\delta(A'_\theta) < \delta^{2\eps} \cN_\delta(A)$.

Now set $A_\theta = A \setminus A_\theta'$.
Then, on the one hand,
\[
\nu(\R^d \setminus A_\theta) = \nu(A_\theta') \ll \frac{\cN_\delta(A_\theta')}{\cN_\delta(A)} \leq \delta^{\eps}.
\]
On the other hand, for any $P \in \cD_\delta^\ur$, either $P \in \cQ_{\mathrm{bad}}(\theta)$, then $A_\theta \cap \phi_\theta^{-1}P = \emptyset$, or $P \notin \cQ_{\mathrm{bad}}(\theta)$, hence in any case, 
\[
\nu(A_\theta \cap \phi_\theta^{-1}P) \leq \vol(R)^{\alpha +2\eps}.
\]
Consequently, for any $x \in \R^d$,
\[
\nu(A_\theta \cap \phi_\theta^{-1}(x+R)) \leq \vol(R)^{\alpha +\eps},
\]
since $\cN_\delta^\ur(x + R) \ll 1$.
This finishes the proof of this special case.
\end{proof}

We now explain how to deduce the general case.
\begin{proof}[Proof of \Cref{coro:slicing}, general case]
Write  $M = 2 d \abs{\log \delta}$.
For $k \in \{1,\dotsc,M\}$, let $A_{k}$ be the union of cubes $Q\in \cD_{\delta}$ such that $2^{-k}<\nu(Q) \leq 2^{-(k-1)}$ and let $\nu_k = \nu_{|A_k}$.
Hence $\nu=\sum_{k \geq 1}^M \nu_{k} + \nu_{> M}$ where $\nu_{> M}$ is the restriction of $\nu$ to the union of cubes $Q \in \cD_\delta$ such that $\nu(Q) \leq 2^{-M}$.
As $\nu$ is supported on $B_1^{\R^d}$, we must have $\nu_{> M}(\R^d) \ll \delta^{-d} 2^{-M}$, hence $\nu_{> M}(\R^d) \leq \delta$ provided that $\delta \lll_d 1$.

For each $k \in \{1,\dotsc, M\}$, by its construction, $\nu_k$ satisfies the assumption of the special case above.
It also inherits from $\nu$ the non-concentration condition~\ref{it:NCnu} as $\nu_k \leq \nu$.
By the special case already shown, up to assuming $\eps, \delta \lll_{d,(r_{i}), \kappa }1$, there is a subset $\cE_k \subset \Theta$ of measure $\sigma(\cE_k) \leq \delta^{2\eps}$ and such that for every $\theta\in\Theta \setminus \cE_k$, 
there is a set  $A_{k,\theta}\subseteq A_{k}$ satisfying  $\nu_{k}(A_k \setminus A_{k,\theta})\leq \delta^{2\eps}$ and
\[
\sup_{x\in \R^d} ({\phi_\theta}_\star\nu_{k|A_{k,\theta}}) (x+R)\leq \vol(R)^{\alpha + 2\eps}.
\]

Set $\cE \defeq \bigcup_{k=1}^M \cE_{k}$ and for $\theta \notin \cE$, $A_{\theta} \defeq \bigcup_{k=1}^M A_{k,\theta}$.
Then $\sigma(\cE) \leq M \delta^{2\eps} \leq \delta^\eps$,
\[
\nu(\R^d \setminus A_{\theta}) \leq \sum_{k=1}^M \nu_k( A_k\setminus A_{k,\theta})  + \nu_{>M}(\R^d) \leq M\delta^{2\eps} + \delta \leq \delta^\eps.
\] 
and finally, for any $x \in \R^d$, using that the various $\nu_k$ are mutually singular,
\begin{align*}
({\phi_\theta}_\star\nu_{|A_{\theta}}) (x+R) &\leq \sum_{k=1}^M ({\phi_\theta}_\star\nu_{k|A_{k,\theta}}) (x+R)  \\
& \leq M \vol(R)^{\alpha + 2\eps} \\
& \leq \vol(R)^{\alpha + \eps},
\end{align*}
provided that $\delta \lll_{d,\eps} 1$.
As the dependence on $\eps$ for the upper bound on $\delta$ can be removed, this concludes the proof. 
\end{proof}

\section{Effective positive dimension}
\label{sec:positive} \label{sec-intro-posdim}
The goal of this section is to show that a random walk on an arithmetic homogeneous space reaches exponentially fast a positive dimension provided the driving measure is algebraic with respect to the ambient arithmetic structure and the starting point does not satisfy some natural obstructions. 

\bigskip

Let  $G$ be a connected semisimple real linear group with no compact factor.

\begin{definition}[Arithmetic lattice] A subgroup $\Lambda\subseteq G$ is \emph{arithmetic} if there exist 
a semisimple $\Q$-group  $\bG$, a  $\R$-anisotropic normal $\R$-subgroup $\bK\subseteq \bG$, and a Lie group isomorphism $\varphi : (\bG_{\R}/\bK_{\R})^\circ\rightarrow G/Z(G)$ such that $\varphi(\bG_{\Z}\mod \bK_{\R})$ is commensurable to $\Lambda \mod Z(G)$.   
\end{definition}

Here the superscript $\circ$ refers to the identity component (for the smooth topology). By the Borel-Harish Chandra Theorem, an arithmetic subgroup is a lattice in $G$. For more on arithmetic lattices, see  \cite[Chapter IX]{Margulis91}.

\begin{definition}[Algebraic driving measure] 
\label{def:alg}
We say a probability measure $\mu$  on $G$ is \emph{algebraic} with respect to an arithmetic lattice $\Lambda$
if for some (hence any) triple $(\bG, \bK, \varphi)$ as above, the measure
$\phi^{-1}_\star(\mu \mod Z(G))$ is concentrated on $\bG_{\QQ} \mod \bK_{\R}$.
\end{definition}

An alternative characterization, which will not be used, is to ask that there exists a basis of $\kg\defeq \Lie(G)$ in which all elements in $\Ad \Lambda$ and $\Ad_{\star}\mu$-almost all elements in $\Ad G$ are represented by matrices with algebraic entries.

\bigskip
The section is dedicated to the proof of the following statement.

\begin{thm}[Effective positive dimension] \label{thm-positive-dimension}
Let   $G$ be a connected semisimple real linear group with no compact factor, $\Lambda$ an arithmetic lattice in $G$, set $X=G/\Lambda$ equipped with a quotient right $G$-invariant  Riemannian metric, write  $x_{0}=\Lambda/\Lambda\in X$ the basepoint. Let $\mu$ be a finitely supported probability measure on $G$ which is Zariski-dense and algebraic with respect to $\Lambda$. 

 Given parameters $x\in X$, $n\geq 1$, $A,C, \kappa, \rho_0>0$, one has
\begin{equation}\label{eq-posdim-c}
\forall \rho\geq \rho_0,\, \forall y\in X,\quad \mu^{*n}*\delta_{x}(B_{\rho}(y)) \leq C\rho^\kappa
\end{equation}
provided that
$n \geq  \abs{\log \rho_0} +A\max\bigl\{ \abs{\log \dist(x,W_{\mu, \rho_0^{-1}})}, \dist(x,x_{0})\bigr\}$ and the conditions $A,C \ggg_{X, \mu} 1\ggg_{X, \mu}   \kappa$. 
 \end{thm}

 The assumption of algebraicity on $\mu$ can be removed if $\Gamma_{\mu}$ admits a finite orbit on $X$ (see the proof of \Cref{dimpos-fo}).

The most simple example where \Cref{thm-positive-dimension} applies  is the case where $G=\SL_{2}(\R)$, $\Lambda=\SL_{2}(\Z)$, and $\mu$ is supported on $\SL_{2}(\R\cap \QQ)$ with finite support containing two hyperbolic matrices with pairwise distinct eigenlines. This context should provide insight to follow the proof. A classification of arithmetic lattices of most simple Lie groups can be found in \cite[Chapters 6 \& 18]{Morris15}. See also \cite[p.296-298]{Margulis91}.

\bigskip
In order to prove the theorem, we may simplify the framework a bit.
First, quotients by finite subgroups play no role, hence up to replacing $\bG$ by $\Ad \bG$ we may assume that $\bG$ is of adjoint type, in particular centerless. Thus, $\bG$ is the direct product of its $\C$-simple factors.
We may also suppose $\varphi$ is the identity map, and $\Lambda=G\cap (\bG_{\Z} \mod \bK_{\R})$.
In sum, we have reduced our working framework as follows.

\begin{notation}\label{not-sec-dimpos}
 $\bG\subseteq \SL_{d}$ is a connected adjoint semisimple $\Q$-subgroup, $\bK$ is the maximal $\R$-anisotropic $\R$-factor of $\bG$.
We fix a norm on $\C^d$ and endow $\M_{d}(\C)$ with the induced operator norm. 
 We define  $G=(\bG_{\R}/\bK_{\R})^\circ$, $\Lambda=G\cap (\bG_{\Z} \mod \bK_{\R})$, and set $X=G/\Lambda$ equipped with a quotient right $G$-invariant  Riemannian metric. We let $\mu$ be a  probability measure on $G\cap (\bG_{\QQ} \mod\bK_{\R})$  whose support is finite and generates  a Zariski-dense subgroup $\Gamma_{\mu}$ in $G$.  
\end{notation} 

The rest of the section is formulated within \Cref{not-sec-dimpos}, at the exception of \Cref{Sec-effective-recurrence}, which authorizes a more general framework to accommodate future use in \Cref{Sec-bootstrap}.
The symbols  $O(\cdot)$, $\ll$, $\simeq$, $\lll$ (see \ref{conventions-notations})  refer implicitely to constants possibly depending on \Cref{not-sec-dimpos} (with similar conventions  specified in \Cref{Sec-effective-recurrence}), additional dependences  are indicated in subscript. 

\subsection{Algebraic preliminaries}  \label{Sec-review-arith}

Given a subfield $\L\subseteq \C$ which is algebraic over $\Q$, we study the $\L$-rational points of $X$, and the subset of $\L$-rational points with bounded complexity (Mahler measure). This prepares the proof of the almost Diophantine property in \Cref{Sec-almost-diophantine}. We place ourselves in \Cref{not-sec-dimpos}.

\bigskip
Denote by $\cF$ the set of $\C$-simple factors of $\bG$.
In particular, $\{\bH_{\C} : \bH\in \cF\}$ are commuting subgroup of $\SL_{d}(\C)$ and $\prod_{\bH\in \cF}\bH_{\C} \rightarrow \bG_{\C}, (g_{\bH})_{\bH\in \cF}\mapsto g\defeq \prod_{\bH\in \cF}g_{\bH}$ is an isomorphism. 
For an element $g \in \bG$, its preimage by this isomorphism will be denoted by $(g_\bH)_{\bH \in \cF}$.

Note also that the Galois group $\Gal(\C/\Q)$ acts on $\cF$. 
The following is a standard fact.
\begin{lemma} \label{L-points-bG}
An element $g\in \bG_{\C}$ belongs to $\bG_{\L}$ if and only if for every $\sigma\in \Gal(\C/\L)$, $ \bH\in \cF$, 
\[\presigma{(g_{\bH})}=g_{(\presigma{\bH})}.\]
\end{lemma}

\begin{proof}
 The factor decomposition $g=\prod_{\bH\in \cF}g_{\bH}$ is unique and  
$\bG_{\L}$ is the set of fixed points of $\Gal(\C/\L) \acts \bG_{\C}$. 
\end{proof}

We define the set of \emph{$\L$\nobreakdash-rational points} of $G$ by $G_{\L} \defeq G\cap (\bG_{\L} \mod \bK_{\R})$. The \emph{$\L$\nobreakdash-rational points} of $X$ are then obtained by quotienting by $\Lambda$, namely $X_{\L} \defeq  \{\, g\Lambda \,:\,g\in G_{\L} \,\}$.  Here also we have a characterization in terms of the Galois action. Let $\cF_{nc}\subseteq \cF$ be the subset of $\C$-simple factors of $\bG$ that are not contained in $\bK$, so that $\bG = \bK \times \prod_{\bH \in \cF_{nc}} \bH$.
Note that $\cF_{nc}$ consists of the $\C$-simple factors that are either not defined over $\R$ or defined over $\R$ and $\R$-isotropic.
This is because if $\bH \in \cF$ is not defined over $\R$ then the group of $\R$-points of the product of $\bH$ with its complex conjugate is homeomorphic to $\bH_\C$ which is never compact.
Observe that an element $g$ in $G = (\bG/\bK)_\R^\circ$, being a coset of $\bK$, determines a vector $(g_\bH)_{\bH \in \cF_{nc}}$.

\begin{lemma}[$\L$-rational points of $G$] \label{L-points-G}
An element $g\in G$ belongs to $G_{\L}$ if and only if for every $\sigma\in \Gal(\C/\L)$, every $\bH\in \cF_{nc}$ such that $\presigma{\bH}\in \cF_{nc}$, we have 
\[\presigma{(}g_{\bH})=g_{(\presigma{\bH})}.\] 
\end{lemma}

\begin{remark}
From this, we see that if $(\prod_{\bH\in \cF_{nc}}\bH)$ is defined over $\L$, then our notion of $\L$-points on $G$ coincides with the usual one.
\end{remark}

\begin{proof} The ``only if'' part follows immediately from \Cref{L-points-bG}.
Let us show the ``if'' part. 
We complete the vector $(g_\bH)_{\bH \in \cF_{nc}}$ at $\bH \in \cF \setminus \cF_{nc}$ as follows.
If there is $\sigma\in \Gal(\C/\L)$ such that $\presigma{\bH}\in \cF_{nc}$ then set $g_{\bH} \defeq \prescript{\sigma^{-1}\!}{}(g_{(\presigma{\bH})})$.
The assumption on $g$ guarantees this $g_{\bH}$ is well defined independently of the choice for $\sigma$.
If there  is no such $\sigma$, we set $g_{\bH} = \Id$.
Then $\prod_{\bH \in \cF}g_{\bH}$ belongs to $\bG_{\L}$ by \Cref{L-points-bG}  and it projects to $g$ as desired. 
\end{proof}

\bigskip
We will be led to consider $\L$-rational points of small complexity. The complexity will be expressed using the Mahler measure which we now recall.

\begin{definition}[Mahler measure]
Let $\alpha\in \overline{\Q}$ and $\chi_{\alpha}=\sum_{i=0}^ma_{i}X^i\in \Z[X]$ its minimal polynomial with coefficients in $\Z$ (and $m=\deg(\alpha)$). 
The \emph{Mahler measure} of $\alpha$ is  
\[\Mah(\alpha) \defeq  \abs{a_{m}} \prod_{i=1}^m\max (1, \abs{\alpha_{i}})\]
where $\alpha_{1}, \dots, \alpha_{m}$ enumerate the roots of $\chi_{\alpha}$ (i.e. the Galois conjugates of $\alpha$).
\end{definition}

 For example, if $\alpha\in \Q\smallsetminus\{0\} $ then $\Mah(\alpha)=\max(\abs{p}, \abs{q})$ where $p,q\in \Z^2$ satisfy $\alpha=p/q$ with $\gcd(p, q)=1$. Moreover $\alpha\in \overline{\Q}\smallsetminus\{0\}$ satisfies $\Mah(\alpha)=1$ if and only if $\alpha$ is a root of unity (Kronecker's theorem). Another notion of algebraic complexity is given by the \emph{denominator} of $\alpha$, namley $\den(\alpha)$ is the smallest integer $Q\geq 1$ such that $Q\alpha$ is an algebraic integer ($Q\alpha\in \cO_{\QQ}$). This notion is looser than the Mahler measure, in the sense that $ \den(\alpha)\leq \Mah(\alpha)$ as we see below.
  The next lemma records the basic properties of the Mahler measure that we use in the paper.
  
\begin{lemma} \label{Facts-Mahler}
\begin{itemize} 
\item[(i)]  Let $n\geq 1$, $P\in \Z[X_{1}, \dots, X_{n}]$, $\alpha_{1}, \dots, \alpha_{n}\in \overline{\Q}$.
Let $\mathcal{L}(P)$ be the sum of the absolute values of the coefficients of $P$. 
Let $k_{i}$ be the degree of $P$ in the variable $X_{i}$, and $d=[\Q(P(\alpha_{1}, \dots, \alpha_{n})):\Q]$, $d_{i}=[\Q(\alpha_{i}):\Q]$. Then we have 
\[\Mah(P(\alpha_{1}, \dots, \alpha_{n}))\leq \mathcal{L}(P)^{d}\prod_{i=1}^n\Mah(\alpha_{i})^{k_{i} dd_{i}^{-1}}.\]
\item[(ii)] If $\alpha\neq 0$, then we have $\Mah(\alpha)=\Mah(\alpha^{-1})$,  as well as $ \den(\alpha)\leq \Mah(\alpha)$ and $\Mah(\alpha)^{-1}\leq \abs{\alpha}\leq \Mah(\alpha)$.
\end{itemize}
\end{lemma}

\begin{proof}
By \cite[Proposition 1.6.6]{BombieriGubler06}, for any number field $\mathbb{K}$, there exists  $(\abs{\cdot}_{v})_{v}$  a certain system of representatives for the places of $\mathbb{K}$, which  satisfies the product formula and that for all $\alpha \in \mathbb{K}$, we have 
$\Mah(\alpha) = \prod_{v}\max(1, \abs{\alpha}_{v})^{[\Q(\alpha) : \Q]}$.

Set $\mathbb{K} = \Q(\alpha_1,\dotsc,\alpha_n)$ and $\alpha=P(\alpha_{1}, \dotsc, \alpha_{n})$.
For any finite place $v$, the ultrametric property yields $\abs{\alpha}_{v} \leq \prod_{j} \max(1, \abs{\alpha_j}_{v})^{k_{j}}$. 
For an infinite place~$v$,  the triangle inequality gives
\[\abs{\alpha}_{v}\leq \mathcal{L}(P)\prod\nolimits_{j} \max(1, \abs{\alpha_j}_{v})^{k_j}.\]
Noting that $\abs{\alpha}_{v}$ can be replaced by $\max(1, \abs{\alpha}_{v})$ in the above inequalities, then taking the product over all places,  we obtain (i). 

For (ii), $\Mah(\alpha)=\Mah(\alpha^{-1})$ is straightforward by the product formula, and $\Mah(\alpha)^{-1}\leq \abs{\alpha}\leq \Mah(\alpha)$ follows. To check $ \den(\alpha)\leq \Mah(\alpha)$, write $\chi_{\alpha}=\sum_{i=0}^m a_{i}X^i\in \Z[X]$ ($a_{m}\neq 0$), then $a_{m}\alpha$ is a zero of $X^m+\sum_{i=0}^{m-1}a_{m}^{m-1-i}a_{i}X^{i}$, whence $a_{m}\alpha\in \cO_{\QQ}$, and the claim follows using $\Mah(\alpha)\geq \abs{a_{m}}$. 
\end{proof}
 
Given a matrix $A=(\alpha_{i,j})$ with  coefficients in $\overline{\Q}$, we set $\Mah(A) \defeq \max_{i,j}\Mah(\alpha_{i,j})$.  
Given $Q\geq 1$, we set 
\[
G_{\L, Q} \defeq \{\,g\in G_{\L}  \,:\, \max_{\bH \in \cF_{nc}} \Mah(g_{\bH})\leq Q \,\}\; \text{ and }\; X_{\L, Q} \defeq  \{g\Lambda \,:\,g\in G_{\L, Q} \}.
\]

\begin{lemma}[Polynomial separation 1] \label{separation-XK}
For $M \ggg_{\L} 1$, for every $Q \geq 2$, the set $X_{\L, Q}$ is $Q^{-M}$-separated and included in $\{\inj\geq Q^{-M}\}$.
\end{lemma}
  
\begin{proof}
Let $g_{1}, g_{2}\in G_{\L, Q}$, let  $\gamma_{1}, \gamma_{2} \in \Lambda$ such that  $g_{1}\gamma_{1}\neq g_{2}\gamma_{2}$, write $\omega=g_{1}\gamma_{1}\gamma_{2}^{-1}g_{2}^{-1}$.
We need to show $\dist(g_1\gamma_1, g_2 \gamma_2) = \dist(\omega, \Id)$ is bounded below by  $Q^{-M}$ provided $M\ggg_\L 1$.
Of course we may assume $\dist(\omega, \Id)\leq 1$. 

In this case  we have $\norm{\gamma_{1}\gamma_{2}^{-1}} \leq Q^{O(1)}$. This allows us to bound the absolute values of matrix entries of $(\gamma_{1}\gamma_{2}^{-1})_\bH$ for all $\bH \in \cF_{nc}$ and at all Archimedean places.
At non-Archimedean places, note that $(\gamma_{1}\gamma_{2}^{-1})_\bH$ is the projection of a $\Z$-point to a simple factor, hence the matrix entries have absolute values bounded by $O(1)$.
Put together, we find $\Mah\bigl( (\gamma_{1}\gamma_{2}^{-1})_\bH \bigr) \leq Q^{O(1)}$ for every $\bH \in \cF_{nc}$.

By \Cref{Facts-Mahler} (i), the previous paragraph implies   $\Mah(\omega_\bH)\leq Q^{O_\L(1)}$, in other terms $\omega\in G_{\L, Q^{O_{\L}(1)}}$.
Since $\omega \neq \Id$, there is some $\bH \in \cF_{nc}$ such that $\omega_\bH \neq \Id$.
Then \Cref{Facts-Mahler} (ii)  gives the lower bound $\norm{\omega -\Id} \gg \norm{\omega_\bH - \Id}\geq Q^{-O_{\L}(1)}$ and the result follows. 
\end{proof}

It will be useful to us to record a similar separation for the set $D_{Q}=D_{\Q, Q} \defeq \{g\Lambda \,:\,g\in G_{\Q} \,:\,\den(g)\leq Q\}$ of $\Q$-rational points with bounded denominator. Here $\den(g) \defeq \max_{\bH\in \cF_{nc}} \den(g_{\bH})$, with $\den(A) \defeq \max_{i,j}\den(\alpha_{i,j})$ for any matrix $A=(\alpha_{i,j})$.

 \begin{lemma}[Polynomial separation 2]  \label{DQ-sep}
For $M\ggg1$, for all $Q\geq 2$, the set $D_{Q}$ is $Q^{-M}$-separated and included in $\{\inj\geq Q^{-M}\}$.
 \end{lemma}
 
Note that $D_{Q}$ contains $X_{\Q, Q}$ (by \Cref{Facts-Mahler} (ii)) so the lemma implies polynomial separation for $X_{\Q, Q}$, whence a part of  \Cref{separation-XK}.
Here however, polynomial separation for $D_{\L, Q}$ only holds for $\L=\Q$ (for $\L\neq \Q$, we have $D_{\L, 1}$ is dense in $X$) whereas \Cref{separation-XK} concerns all $\L$-points of small Mahler measure.

\begin{proof}
Arguing as for \Cref{separation-XK}, we only need to check that for $g\in G_{\Q}\smallsetminus \{\Id\}$ with $\den(g)\leq Q$, we have $\dist(g, \Id)\geq Q^{-M}$ provided $M\ggg 1$.

Recall that the product of the conjugates of an algebraic integer is an integer, in particular for any  $\alpha\in Q^{-1}\cO_{\QQ}$, we have $\prod_{\alpha' \in \Gal(\C/\Q).\alpha} \abs{\alpha'} \geq Q^{-[\Q(\alpha):\Q]}$. Considering an element $\prod_{\bH'\in \cF}g_{\bH'}$  in $\bG_{\Q}$ projecting to $g$ modulo $\bK_{\R}$,  and letting $\bH\in \cF_{nc}$ such that $g_{\bH}\neq \Id$, we deduce from the condition on the denominator of $g$ that
\[\prod_{\bH'\in \Gal(\C/\Q).\bH}\|\Id- g_{\bH'} \| \gg Q^{-O(1)}.\]
For $\bH'\notin \cF_{nc}$ we have $\bH'$ defined over $\R$ and $\bH'_{\R}$ compact. The condition $g\in G_{\Q}$  implies $g_{\bH'}\in \bH'_{\R}$, so $\|\Id- g_{\bH'} \| $ is then bounded above by a constant depending on the data of \Cref{not-sec-dimpos}. It follows that 
\[
\| \Id-  g\| \gg Q^{-O(1)}
\]
 thus concluding the proof.
\end{proof}

\subsection{Effective recurrence} \label{Sec-effective-recurrence}
 
We recall the recurrence properties of the $\mu$-walk on $X$, namely \Cref{effective-recurrence} stating that the tail probabilities of the distribution $\mu^{*n}*\delta_{x}$ decay exponentially with a rate which is independent of the couple $(n,x)$  provided $n$ is large enough depending on $x$. 

As such result will be useful in \Cref{Sec-bootstrap} as well, we allow momentarily a \emph{more general setting}: $G$ is a noncompact connected semisimple real linear group, $\Lambda \subseteq G$ is a (non-necessarily arithmetic) lattice, $X=G/\Lambda$ equipped with a quotient right $G$-invariant  Riemannian metric and $\mu$ is a probability measure on $G$ with finite exponential moment generating a Zariski-dense subgroup of $G$. Recall that the finite exponential moment assumption means that for some $\eps>0$, we have $\int_{G}\norm{\Ad g}^\eps \dd \mu(g)<\infty$. In this subsection,  the symbols $\ggg, \gg, O(\,\cdot\,), \simeq$ refer to constants depending on this setting.

\begin{proposition}[Effective recurrence] \label{effective-recurrence}
Let $s_{0}>0$ be small enough. Then for every $x\in X$, $n\ggg \dist(x,x_{0})$, and $R>0$, we have 
\[(\mu^{*n}*\delta_{x})\{\, y\in X\,:\,\dist(y, x_{0})\geq R\,\}\ll e^{-s_{0}R}. \]
\end{proposition}

\begin{remark}
Here we may replace $\dist(\,\cdot\,,x_{0})$ by $\abs{\log \inj(\cdot)}$, the absolute value of the logarithm of the injectivity radius, see \Cref{comparison2} below. This observation will be useful in \Cref{Sec-bootstrap} where the injectivity radius will occur naturally in the context of robust measures introduced for the bootstrap argument in \Cref{Sec-bootstrap}.
\end{remark}

\Cref{effective-recurrence} is well-known to experts but it might be hard to find it explicitly in the literature. We give details about the proof, by explaining how it is connected to written results. We denote by $P_{\mu}$ the \emph{Markov operator} associated to $\mu$, namely for any  measurable non-negative function $f :X\rightarrow [0, +\infty]$, we have 
\[P_{\mu}f: x\mapsto \int_{G} f(gx) \dd \mu(g).\]
Recurrence properties of the $\mu$-walk are usually obtained using a drift function, i.e. a function that is proper and contracted by $P_{\mu}$ up to an additive constant.  
 This method originates in the work of Eskin-Margulis-Mozes~\cite{EskMarMoz98} on the quantitative Oppenheim conjecture. It was introduced in the context of random walks by Eskin-Margulis \cite{EskinMargulis}, then generalized by Benoist-Quint \cite{BQ_rec}, see also B\'enard-Saxc\'e \cite{BenDeS21}. A survey  on that topic is presented in \cite{benoist_rec}. 

\begin{lemma}[Drift away from infinity] \label{drift-infinity}
There exists a proper map $u_{0}: X\rightarrow [1, +\infty)$ such that for all small enough $\lambda, s>0$,  for all $n \geq 1$, 
\[P_{\mu}^n (u_{0}^{s}) \leq  e^{- s \lambda n}  u_{0}^s +  2.\]
\end{lemma}

\begin{proof} See  \cite[Equation (39)]{EskinMargulis} where it is established  for some $\lambda=\lambda_{0}$, $s=s_{0}$ and some additive error $A>1$ in the place of $2$. Clearly we may then allow $\lambda<\lambda_{0}$.
To allow for smaller $s$ and reduce the additive error, observe that for  $r\in (0, 1]$, the function $t\mapsto t^r$ ($t\in \R^+$) is concave,  whence $P_{\mu}^n (u_{0}^{ s_{0}r})\leq (P_{\mu}^n u_{0}^{ s_{0}})^r$, then apply the bound on $P_{\mu}^n (u_{0}^{ s_{0}})$ and  the inequality $(a+b)^r\leq a^r+b^r$ for all $a,b\in \R^+$.
\end{proof}

Moreover we may compare $\log u_{0}$  to the distance function to the basepoint $x_{0}=\Lambda/\Lambda\in X$.

\begin{lemma}[Comparison 1] \label{comparison} We have
\[ \dist(x, x_{0})  \ll \log u_{0}(x)+1 \quad \quad \quad \log u_{0}(x)\ll  \dist(x, x_{0}) + 1.\]
\end{lemma}

\begin{proof}
Writing $f=\log u_{0}$, it is shown in  \cite[Equation (36)]{EskinMargulis} that for large enough $A>1$, for all $g\in G$, $x\in X$, 
\[ |f(gx)-f(x)|\leq A \log \|\Ad g\|.\]
 As we have $\dist(x, x_{0}) + 1 \simeq \inf \set{\log \|\Ad g\|: g x_{0} = x } + 1$,  we deduce the upper bound $f(x)\leq A' (\dist(x, x_{0}) +1) + f(x_0)$ where $A'$ depends only on $A$ and $G,X, \|\cdot\|$.
 For the lower bound, observe that the drift property implies that for some $R, r > 0$ and for every $x$ such that $f(x)>R$, there exists some $y\in B^X_{R}(x)$ for which $f(y) < f(x) - r$.
 It follows that for all $x\in \{ f > R \}$,
\[\dist(x, \{f\leq R\})\leq \frac{f(x)-R}{r} R.\]
Note that $R$ can be chosen arbitrarily large, in particular, we may assume $ f(x_0)<R$.
Then
\[f(x)\geq \frac{r}{R}\dist(x, x_0) + R - \diam\{f\leq R\}\]
where $\diam$ stands for the diameter. 
\end{proof}

\begin{proof}[Proof of \Cref{effective-recurrence}]
It follows from \Cref{drift-infinity}, \Cref{comparison}, and the Markov inequality. 
\end{proof}

We record the following, which allows to compare $u_{0}$ to the inverse of the injectivity radius via \Cref{comparison}. 

\begin{lemma}[Comparison 2] \label{comparison2} For $x\in X$, we have
\[\abs{\log \inj(x)}  \ll \dist(x, x_{0})+1\quad \quad \quad \dist(x, x_{0}) \ll \abs{\log \inj(x)} + 1.\]
\end{lemma}

\begin{proof}
We  start with the left inequality.
Write $x=h \Lambda$ for some $h \in G$ with $\dist(h,\Id) = \dist(x,x_0)$ which results in $\log \norm{\Ad h } \ll \dist(x,x_{0})+1$.
We have   
\[\inj(x) \simeq \min \dist(\Id, h\Lambda h^{-1} \setminus \{\Id\} )\]
 and the result follows because conjugation map by $h$ on $G$ is bi-Lipschitz with bi-Lipschitz constant $\norm{\Ad h}^{O(1)}$.

For the right inequality, let $r=\inj(x)$.
We may assume $r\leq 1$.
By \Cref{effective-recurrence},  for any $x' \in X$, for large enough $n\geq 0$, we have 
\[\mu^{*n}*\delta_{x'}(B_{r}(x))\leq (\mu^{*n}*\delta_{x'})\{y\in X\,:\,\dist(y, x_{0})\geq \dist(x, x_{0})-1\}\ll e^{-s_{0} \dist(x, x_{0})}.\]
Choosing $x'$ in  a subset of full Haar measure, and letting $n \to +\infty$, we know that $ \mu^{*n}*\delta_{x'}$  converges to Haar in Ces\`aro average (by ergodicity of the Haar measure for the $\mu$-walk), so 
\[r^{\dim X}\ll e^{-s_{0} \dist(x, x_{0})}\]
which gives the desired inequality. 
\end{proof}

\subsection{Almost Diophantine property} \label{Sec-almost-diophantine}
Let us place ourselves again in \Cref{not-sec-dimpos}.
We will show \Cref{nonpos->rational} below claiming that if the $n$-step distribution of the $\mu$-walk on $X$ accumulates much on a ball, then the center of this ball must be very close to an algebraic point of small Mahler measure.
Our proof is inspired by the work of Bourgain-Furman-Lindenstrauss-Mozes \cite[Proposition 7.3]{BFLM} (or its generalization by He-Saxc\'e \cite[Proposition 5.2]{HS2022}) where a similar property is established to show positive dimension for walks on the \emph{torus}. A crucial difference though is that their proof makes significant use of the commutativity of the torus, while for us $X$ is covered by $G$ which is not abelian. 
\bigskip
 
We fix $\L\subseteq \C$ a number field such that \emph{$\Gamma_{\mu}\cap G_{\L}$ has finite index in $\Gamma_{\mu}$}.  Given $Q, r>0$, we also set $X_{\L,Q}^{(r)} \defeq \{\, x \,: \,\dist(x, X_{\L,Q})<r \,\}$  the $r$-neighborhood of $X_{\L,Q}$ in $X$ (see \Cref{Sec-review-arith} for the definition of $X_{\L,Q}$).

\begin{proposition}[Almost Diophantine property]\label{nonpos->rational}
Set $x_{0}=\Lambda/\Lambda\in X$.
Given $\eps \in (0, 1/2)$, there exist $C > 1$ and  $\kappa > 0$ such that for $x \in X$, for every ball $B_{\rho}(y)\subseteq X$ with $\rho\leq 1/C$, and every $n \geq \eps |\log \rho| + A \cdot \dist(x, x_{0})$ where $A \ggg 1$, one has:
\[\text{if $\mu^{*n} *\delta_{x}(B_{\rho}(y))\geq \rho^\kappa$ then $y$ belongs to $X_{\L,\rho^{-\eps}}^{(\rho^{1-\eps})}$.}\]
\end{proposition}

\begin{remark}
None of the parameters depend on the choice for $\L$. 
In fact, it is enough to show the statement in the case where $\L$ is the \emph{smallest} number field such that $\Gamma_{\mu}$ is included in $G_{\L}$ up to finite index (it exists because the collection of finite index subgroups of $\Gamma_{\mu}$ is stable by intersection). 
\end{remark}

We start by proving that a back-and-forth $\mu$-random trajectory on $G$ escapes proper algebraic subgroups. Given  a collection $E$ of elements in $G$, we write $\overline{\langle E \rangle}^{\mathrm{Zar}}$ the smallest real algebraic subgroup of $G$ containing $E$.

\begin{lemma}
\label{lm:rw-Zdense}
For every $m\ggg1\ggg\eta>0$, every proper real algebraic subgroup $H\subsetneq G$, we have
$$(\mu^{*m})^{\otimes 2}\left\{(g_{1},g_{2})\,:\, \dim \overline{\langle g_{1}g_{2}^{-1}, H \rangle}^{\mathrm{Zar}} =\dim H \right\} \leq e^{-\eta m}.$$
\end{lemma}

\begin{remark}By the same proof, we can see the estimate is also valid when the variable $(g_{1}g_{2}^{-1})_{(g_{1},g_{2})\sim(\mu^{*m})^{\otimes 2}}$ is replaced by $g\sim \mu^{*m}$, 
but we will not need this fact.
\end{remark}

\begin{remark}
For the example $G=\SL_{2}(\R)$, we can use the following short argument. Every proper algebraic subgroup is included in a conjugate $hBh^{-1}$ of the upper triangular subgroup $B$, so \Cref{lm:rw-Zdense} reduces  to showing that for all cosets $z_{1},z_{2}\in G/B$, the quantity $\mu^{*m}*\delta_{z_{1}}(z_{2})$ decays exponentially, uniformly in $z_{1},z_{2}$. This  result is standard, see \cite[Lemma 11.5, Proposition 11.10]{BQ_book}, showing that $\mu^{*m}*\delta_{z_{1}}$ converges to the Furstenberg measure on $G/B$  with exponential rate and uniformly in $z_{1}$. 
\end{remark}

\begin{proof} 
By \cite[Proposition 3.3]{BreBeck25}, there exists $k\in \N$ (bounded in terms in $d=\dim G$)  such that the  convolution power $\mu^{*k}$ is not supported on any coset modulo a proper algebraic subgroup of $G$. We may thus apply the anti-concentration result \cite[Theorem 1.3]{BreBeck25} to i.i.d. variables of law $\mu^{*k}$. We get  some $\eta'=\eta'(\mu^{*k})>0$ such that for all $m\geq 1$,  all  proper algebraic subgroups $\tilde H\subseteq G$, all  $g_{2}\in G$, 
\begin{equation}\label{eq-esc}
\mu^{*km}\left\{g_{1}\,:\, g_{1} \in \tilde H g_{2}\right\} \leq e^{-\eta' m}.
\end{equation}
If the Lie algebra $\kh:=\Lie(H)$ is not an ideal in $\kg:=\Lie(G)$, we choose $\tilde H=N_{G}(\kh)\subsetneq G$, i.e. $\tilde H$ is the normalizer of $\kh$ in $G$.  \Cref{lm:rw-Zdense} follows by applying \eqref{eq-esc} to this choice of $\tilde H$. 

It remains to deal with the case where $\kh$ is an ideal of $\kg$.  In this case the identity component $H_{e}$ of $H$ for the Zariski-topology is a normal subgroup and up to placing ourselves in $G/H_{e}$, we may assume $H=\{\Id\}$. We are reduced to showing 
$$(\mu^{*m})^{\otimes 2}\left\{(g_{1},g_{2})\,:\, g_{1}g_{2}^{-1} \text{ generates a finite group } \right\} \leq e^{-\eta m}.$$
By considering a proximal irreducible representation of $G$, this estimate follows from Benoist-Quint \cite[Equation (14.34)]{BQ_book} which states\footnote{To be precise here, \cite[Equation (14.34)]{BQ_book} is stated for the distribution $\mu^{*m}$ but the proof adapts to the distribution of $g_{1}g^{-1}_{2}$ where $(g_{1},g_{2})\sim(\mu^{*m})^{\otimes 2}$} that with probability exponentially close to $1$, the variable $g_{1}g_{2}^{-1}$ is proximal in the representation.

\end{proof}

Next, we convert the probabilistic assumption that $\mu^{*n} *\delta_{x}$ gives a lot of mass to a ball into a geometric property of the center of the ball. 
\begin{lemma} \label{almost-diophantine-v-1}
There exists $M\geq 1$ such that the following holds. Let $\nu$ be a probability measure on $X$, let $n\geq 1$, and $B_{\rho}(y)\subseteq X$ a ball with $\rho<1/2$. If $\mu^{*n} *\nu(B_{\rho}(y))\geq \rho^\kappa$ for some $\kappa>0$ and $n\geq  M \kappa |\log \rho|+M $, then there exists a finite set $S\subseteq G_{\L}$ such that 
\begin{itemize}
\item $\max_{g\in S}\Mah(g) \leq M\rho^{-M \kappa}$.
\item $S$ generates a Zariski-dense subgroup of $G$.
\item Writing $y=h\Lambda$ and $\sC_{h}$ the conjugation map by $h$, there exists a map $S\rightarrow \Lambda, g\mapsto \gamma_{g}$ such that 
\[\max_{g\in S} \| g - \sC_{h}(\gamma_{g})\| \leq M \rho^{1-M \kappa}.\]
\end{itemize}
\end{lemma} 

\begin{proof} 
We let $M>1$ be a parameter to specify below and $(\nu, n,B_{\rho}(y))$ as in the statement. The first half of the proof is dedicated to 
\bigskip

\noindent\emph{\underline{Step 1}: Provided $M\ggg 1$,  there exist $g_{1},\dots, g_{2s}\in G$ with $s = \dim G$ and such that
\begin{itemize}
\item $\max_{1 \leq i \leq 2s}\Mah(g_{i}) \leq M\rho^{-M \kappa}$.
\item The set $\{\, g_{2j-1}\,g_{2j}^{-1} : j=1, \dots s \,\}$ is included in $G_{\L}$ and generates a Zariski-dense subgroup of $G$.
\item $\bigcap_{i=1}^{2s} g_{i}^{-1} B_{\rho}(y) \neq \emptyset$.
\end{itemize}  }

Setting $m  \defeq  \lceil M \kappa |\log \rho|+M \rceil $, the idea is to select the $g_{i}$'s randomly with law $\mu^{*m}$ and show that each requirement has a good chance to be satisfied, whence their intersection as well.

The assumption that $\mu$ has finite support in $G_{\QQ}$ and \Cref{Facts-Mahler} (i) guarantee that there exists  $Q_{\mu}>1$, such that for all $m\geq 1$, every $g\in \supp \mu^{*m}$ satisfies $\Mah(g)\leq Q_{\mu}^m$.

Let  $N_{0}\geq 1$ be the index of $\Gamma_{\mu}\cap G_{\L}$ in $\Gamma_{\mu}$. Given a tuple $\ug\in G^{(N_{0}+1)s}$, we write its coordinates   as $\ug=(g_{j,k})$ where   $1\le j\leq s$ and $1\le k\le N_{0}+1$ (for a fixed enumeration of $\llbracket 1, s\rrbracket \times \llbracket 1, N_{0}+1\rrbracket$).
By an iterated application of \Cref{lm:rw-Zdense}, we may fix some $\eta>0$, $ m_{0}\geq1$,  such that for  $m\geq m_{0}$, for a $(\mu^{*m})^{\otimes (N_{0}+1)s}$-proportion at least $ 1-e^{-\eta m}$ of elements $\ug$ the following holds: for every $j \in \llbracket 1, s\rrbracket $,  all pairs of distinct elements $(k_{1}, k_{2}), \dots, (k_{2j-1},  k_{2j})\in \llbracket 1, N_{0}+1\rrbracket^2$, the group
\begin{equation} \label{eqHdiml}
H=\overline{\left\langle g_{1, k_{1}}g_{1,k_{2}}^{-1}, \,\dots, \,g_{j, k_{2j-1}} \,g_{j, k_{2j}}^{-1} \right\rangle}^{Zar} \text{ satisfies } \dim H\geq j.
\end{equation}
Note in particular that $H=G$ when $j=s$. 
Moreover, by the pigeonhole principle, we may always select a specific tuple of $2s$ elements $(k'_{1}, \dots, k'_{2s})$ as above and satisfying additionally that for all $j=1, \dots, s$, the elements $g_{j, k'_{2j-1}}$ and $g_{j, k'_{2j}}$ are in the same right $G_{\L}$\nobreakdash-coset or equivalently $g_{j, k'_{2j-1}}g^{-1}_{j, k'_{2j}}\in G_{\L}$.

Using Hölder's inequality, we have
\begin{align*}
\rho^{(N_{0}+1)s\kappa} &\leq \nu_{n}(B_{\rho}(y))^{(N_{0}+1)s}\\
&= \left(\int \1_{B_{\rho}(y)}(gx) \dd\mu^{* m}(g) \dd\nu_{n-m}(x) \right)^{(N_{0}+1)s}\\
&\leq \int \prod_{j=1}^{s}\prod_{k=1}^{N_{0}+1} \1_{B_{\rho}(y)}(g_{j,k}x)  \dd(\mu^{* m})^{\otimes (N_{0}+1)s}(\ug) \dd\nu_{n-m}(x)\\
&= \int \nu_{n-m}\bigl(\cap_{j,k}\,g_{j,k}^{-1} B_{\rho}(y)\bigr)  \dd(\mu^{* m})^{\otimes (N_{0}+1)s}(\ug)
\end{align*}
from which it follows that for a $(\mu^{* m})^{\otimes (N_{0}+1)s}$-proportion at least $\rho^{(N_0+1)s\kappa}/2$ of elements $\ug\in G^{(N_{0}+1)s}$, we have 
\begin{equation} \label{eqintnonemp}
\nu_{n-m}\bigl(\cap_{j,k} \,g_{j,k}^{-1} B_{\rho}(y)\bigr)\geq \rho^{(N_0+1)s\kappa}/2,
\end{equation}
in particular the intersection is non-empty.

Step $1$ now follows from the combination of the three previous paragraphs. More precisely, we ask that $m$ is large enough so that $m\geq m_{0}$ and $ \rho^{sN_0\kappa}/2>e^{-\eta m}$, or more simply $m \geq s\eta^{-1} N_0 \kappa |\log \rho|+\eta^{-1}\log 2+m_{0}$. This allows to select a tuple $\ug$ satisfying both \eqref{eqHdiml} and  \eqref{eqintnonemp}. The family $(g_{i})_{i=1\dots, 2s}$ defined by $g_{2j-1}:=g_{j, k'_{2j-1}}$, $g_{2j}:=g_{j, k'_{2j}}$ then satisfies the required conditions up to choosing $M$ large enough in terms of the data of \Cref{not-sec-dimpos}. 

\bigskip

\noindent\emph{\underline{Step 2}: Conclusion.} Set $S \defeq \{g_{2j-1}\,g_{2j}^{-1} : j=1, \dots s \}$ where the $g_{i}$'s come from Step 1.  The first item in Step 1 and \Cref{Facts-Mahler} (i) together give the desired bound on the Mahler measure of elements in $S$ (with a slightly bigger $M$).
In order to conclude, we need to interpret the last condition in Step 1 in terms of conjugation by $h$. 
This condition amounts to the existence of $\gamma_{1},\dots, \gamma_{s}$ in  $\Lambda$ such that in $G$, we have
\[\bigcap_{i=1}^{2s} g_{i}^{-1} B_{\rho}h \gamma_{i} \neq \emptyset\]
where $B_{\rho}$ denotes the ball of radius $\rho$ centered at the identity in $G$.  
The non-emptiness of $g_{2j-1}^{-1} B_{\rho}h \gamma_{2j-1}\cap g_{2j}^{-1} B_{\rho}h \gamma_{2j}$ gives
\[  \sC_{h} (\gamma_{2j}\gamma_{2j-1}^{-1})\in  B_{\rho}g_{2j}g_{2j-1}^{-1}B_{\rho} \subseteq g_{2j}g_{2j-1}^{-1} +B^{\M_{d}(\C)}_{\rho'}\]
where $\sC_{h}(x) \defeq hxh^{-1}$,  $ \rho'\ll\rho \| g_{2j}g_{2j-1}^{-1}\|$ and $B^{\M_{d}(\C)}_{\rho'}$ is the open ball of radius $\rho'$ centered at $0$ in $\M_{d}(\C)$. 
The desired bound on $\rho'$ follows using $\| g_{2j}g_{2j-1}^{-1}\|\ll \Mah(g_{2j}g_{2j-1}^{-1})$  (\Cref{Facts-Mahler} (ii)) and the bound on the Mahler measure of $g_{2j}g_{2j-1}^{-1}$ established above.
\end{proof}

The next lemma allows to extract information on $h$ from the condition on $\sC_{h}$ obtained in \Cref{almost-diophantine-v-1}.

\begin{lemma}[From $\sC_{h}$ to $h$] \label{Ch->h}
For   $M\ggg_{\L}1$, the following holds. Let $Q>2$, let $r \in (0,Q^{-M}]$, let $S\subseteq G_{\L}$ be a finite set generating a Zariski-dense subgroup of $G$ and such that $\max_{g\in S} \Mah(g)\leq Q$. Let $h\in G$ such that $\|h\|\leq Q$ and for all $g\in S$, there exists $\gamma_{g}\in G_{\L}$ with $\Mah(\gamma_{g})\leq Q$ and such that  
\[\|g- \sC_{h}(\gamma_{g}) \| \leq r.\]
Then there exists $h'\in G_{\L}$ such that $\dist(h,h')\leq r Q^M$ and $\Mah(h')\leq Q^M$. 
\end{lemma}

\begin{proof}
\noindent\emph{\underline{Step 1}: Construct a candidate $\varphi$ for $h'^{-1}$.}

In the next paragraphs, we define a finite dimensional vector space $V$ (over $\C$ or whichever big algebraically closed field we fix to talk about algebraic groups), an embedding $\sC \colon G\hookrightarrow \GL(V)$ and an element $\varphi\in \End(V)$ which will later turn out to be $\sC(h'^{-1})$.
 
Set $\bG_{nc} \defeq \prod_{\bH\in \cF_{nc}} \bH \subseteq \bG$ and let $\bbF\subseteq \R$ be the smallest real number field on which $\bG_{nc}$  is defined.  We let $V \subseteq \M_{d}$ be the vector space spanned by $\bG_{nc}$ (equivalently by $G$ that we identify with $\bG_{nc, \R}^\circ$). Note that $V$ is defined over $\bbF$.  

The morphism $\sC: \bG_{nc}\rightarrow \GL(V), g\mapsto \sC_{g}$ is defined over $\bbF$ and injective homomorphism of algebraic groups (recall that $\bG$ and hence $\bG_{nc}$ is centerless).
It follows that its image $Z$ is an $\bbF$-subvariety (see \cite[Corollary 1.4]{Borel}) 
and  $\sC$ is  an isomorphism of $\bbF$-varieties between $\bG_{nc}$ and $Z$ (see \cite[Corollary 3.1.11]{Zimmer84}. Also, $\sC(G)$ has finite index in $Z_{\R}$, the group of $\R$-points of $Z$ (see \cite[Corollary 3.19]{PlatonovRapinchuk}).

To define $\varphi$, note that $V$ is the unital associative algebra generated by $S$.
For dimension reasons, $V$ is the linear span of products of at most $\dim V$ elements of $S$.
Hence, we may extract a basis $(v_{1}, \dotsc, v_{\dim V})$ of $V$ consisting of $v_{i}=g_{i,1}\dots g_{i,m_{i}}$ where each $g_{i,j}\in S$ and $m_{i} \leq \dim V$.

We then define $\varphi\in \End(V)$ as the linear endomorphism of $V$ such that for every $1 \leq i \leq \dim V$ 
\[\varphi(v_{i}) =\gamma_{g_{i,1}}\dots \gamma_{g_{i, m_{i}}}.\]
where the $\gamma_{g_{i,j}} \in G_\L$ are given by the assumption in \Cref{Ch->h}.

\bigskip
\noindent\emph{\underline{Step 2}: Properties of $\varphi$.}

Clearly, $(v_i)_{1\leq i \leq \dim V}$ is a $\QQ$-basis of the $\QQ$-strucuture $V_{\QQ}$ of $V$.
It follows that $\varphi$ has algebraic coefficients (i.e. it preserves $V_{\QQ}$) because $\phi(v_i) \in G_\L \subset G_{\QQ}$. Moreover $\Mah(\varphi)\leq Q^{O_{\L}(1)}$,   thanks to \Cref{Facts-Mahler} and the assumptions that the elements $g, \gamma_{g}$ ($g\in S$) belong to $G_{\L}$, and satisfy
$\max_{g \in S} \Mah(g) \leq Q$ and $\max_{g \in S} \Mah(\gamma_g) \leq Q$.

Next we show that $\varphi$ is very close to $\sC_{h^{-1}}$. 
Indeed, starting from the assumption $\max_{g\in S}\|g- \sC_{h}(\gamma_{g}) \| \leq r$, by a simple induction, we see that for any integer $m \geq 1$ and all $g_1,\dotsc,g_m \in S$ we have 
\[\|g_1 \dotsm g_m - \sC_{h}(\gamma_{g_1} \dotsm \gamma_{g_m}) \| \leq r Q^{O(m)}.\] 
Remembering also $\norm{h} \leq Q$ and $\det h=1$, we obtain, for every $i \in \{1,\dotsc,\dim V\}$, 
\[\norm{\sC_{h^{-1}}v_{i}- \phi(v_{i})} \leq r Q^{O(1)}.\] 
Combined with  $\Mah(v_{i})\leq Q^{O_{\L}(1)}$ and \Cref{Facts-Mahler}, we get 
\begin{equation}
\label{eq-sC1}
\| \sC_{h^{-1}} -\varphi\| \leq rQ^{O_{\L}(1)} \leq Q^{-M+O_{\L}(1)}.
\end{equation}

Finally, up to imposing $M\ggg_{\L} 1$, we have
\[\text{(1) $\varphi (g)= \gamma_{g}$ for all $g \in S$  \,\,\,\,\,\,\,\,\,\,\,\,\,\,\,\,\,\,\,\,\,\,\,\,\,\,\,\,\,(2) $\varphi$ belongs to $Z$.}\]
Indeed, by \eqref{eq-sC1} we have for every $g \in S$, $\|\varphi (g)- \gamma_{g}\|\leq Q^{-M+O_{\L}(1)}$. Then (1) follows from the fact that $\Mah(\varphi)$, $\Mah(g)$, $\Mah(\gamma_{g})$ are bounded by $Q^{O_{\L}(1)}$ and  \Cref{Facts-Mahler} (ii) implying that points of small Mahler measure that are very close are actually equal. And (2) comes from a similar reasoning, stated in the general \Cref{rational-absorption} below. 

\bigskip
\noindent\emph{\underline{Step 3}: Conclusion.}

As $\varphi\in Z$, there exists $h'\in \bG_{nc}$ such that $\sC_{h'^{-1}}=\varphi$. Note that $h'\in \bG_{nc, \R}$ because $\varphi$ and $\sC$ are defined over $\R$, and up  to assuming $M \ggg_{\L} 1$, we deduce from \eqref{eq-sC1} that $h'\in G = \bG_{nc, \R}^\circ$.  
From \eqref{eq-sC1}, we also have  $\dist(h,h')\leq r Q^{O_{\L}(1)}$.

Recall that $\varphi$ is an $\L$-rational point and the inverse map $\sC^{-1}$ is defined over $\F$.
It follows that $h'\in \bG_{nc, \langle \F, \L\rangle}^\circ$ where $\langle \F, \L\rangle$ denotes the number field generated by $\F$ and $\L$.
Recalling also $\Mah(\varphi)\leq Q^{O_{\L}(1)}$, it follows from \Cref{Facts-Mahler} that $\Mah(h')\leq Q^{O_{\L}(1)}$.

It remains to check $h'\in G_{\L}$. Write  any element $x\in G$ as $x=\prod_{\bH\in \cF_{nc}}x_{\bH}$. Let $\sigma\in \Gal(\C/ \L)$, let $\bH\in \cF_{nc}$ such that $\presigma{\bH}\in \cF_{nc}$. 
For $g\in S$, using that $\varphi (g)= \gamma_{g}$ as well as the assumption $g, \gamma_{g}\in G_{\L}$ and  \Cref{L-points-G}, we have 
\[\sC_{\presigma{(h'_{\bH})} } (\presigma{(\gamma_{g\, \bH})}) = \presigma{g_{\bH}}= g_{\presigma{\bH}}= \sC_{(h'_{\presigma{\bH}})} (\gamma_{g\, \presigma{\bH}}))=\sC_{(h'_{\presigma{\bH}})} (\presigma{(\gamma_{g\, \bH})}) . \]
Hence $\presigma{(h'_{\bH})}=h'_{\presigma{\bH}}$ because $\{\presigma{\gamma}_{g\, \bH}\,: \,g\in S\}$ generates a Zariski-dense subgroup of $\presigma{\bH}$. In view of  \Cref{L-points-G}, we conclude  $h'\in G_{\L}$.
\bigskip
\end{proof}

We record the following general lemma which was used in the above proof.

\begin{lemma}[Absorption]\label{rational-absorption}
Let $\K$ be a number field.
Let $n\geq 1$ and $Z\subseteq \C^{n}$ be a $\K$-subvariety. There exists $C>0$  such that any  $a \in \K^n$ satisfying $\inf_{z\in Z} \|a- z\|\leq \frac{1}{C} \Mah(a)^{-C}$ must actually belong to $Z$.
 \end{lemma}
 
\begin{proof}
Let $P_{1}, \dots, P_{s}$ be polynomial maps on $\C^n$ with  coefficients in $\K$ and such that $Z=\bigcap_{i\leq s} \{P_{i}=0\}$.  
Let $C>0$ be a parameter to specify below, and $(a,z)\in \K^n\times Z$ such that $ \|a- z\|\leq \frac{1}{C} Q^{-C}$ where $Q \defeq \frac{1}{2}\Mah(a)$. For all $i\leq s$, one has 
\[|P_{i}(a)|= |P_{i}(a-z+z)-P_{i}(z) |\ll_{P_{i}} \frac{1}{C}Q^{-C} \|z\|^{\deg P_{i}-1}\ll_{P_{i}} \frac{1}{C}Q^{-C+\deg P_{i}-1} \]
where the last inequality uses that $\|z\|\ll Q$. 
On the other hand, \Cref{Facts-Mahler}(i) implies that $\Mah(P_{i}(a))\leq \bigl(C_{P_{i}}Q^{\deg P_{i}}\bigr)^{[\K : \Q]}$ where $C_{P_{i}}>0$ is a constant depending only on $P_{i}$.
Hence, by \Cref{Facts-Mahler}(ii), $|P_{i}(a)|$  is either $0$ or its modulus is greater than $\bigl(C_{P_{i}}Q^{\deg P_{i}}\bigr)^{-[\K : \Q]}$. Choosing $C$ large enough in terms of  $\K$ and  $P_{i}$, we then deduce 
\[P_{i}(a)=0.\qedhere\]
\end{proof}

We finally establish the almost Diophantine property. 

\begin{proof}[Proof of \Cref{nonpos->rational}] 
We may assume $\L$ to be the \emph{smallest} number field such that $\Gamma_{\mu}\cap G_{\L}$ has finite index in $\Gamma_{\mu}$. Thus $\L$ is determined by \Cref{not-sec-dimpos} and does not appear in subscript of $O(\cdot)$, $\lll$, etc. We fix  a constant $M>1$ depending on \Cref{not-sec-dimpos} and which is large enough to apply  Lemmas \ref{almost-diophantine-v-1}, \ref{Ch->h}.
We let $A,C,\kappa>0$ be parameters to specify along the proof. 
We may choose $C, \kappa^{-1}\ggg_{\eps}1 $ so that  whenever $\rho\leq 1/C$, we have $\eps |\log \rho| \geq M\kappa |\log \rho| +M$.
Hence any integer $n$ at least $\eps |\log \rho|+ A \cdot \dist(x, x_{0})$ satisfies the condition of application of \Cref{almost-diophantine-v-1}.
Let $S$ be as in \Cref{almost-diophantine-v-1}. By imposing $A\ggg 1$, we may assume $n$ also satisfies the condition of \Cref{effective-recurrence}.

Write $y=h\Lambda$ with $h$ of minimal norm. We observe that $\| h \|$ is bounded by a small power of $\rho^{-1}$.
Indeed, the assumption $\mu^{*n}*\delta_{x}(B_{\rho}(y))\geq \rho^\kappa$ implies via \Cref{effective-recurrence} that $e^{-s_{0}\dist(y,x_{0})}\gg \rho^{2\kappa}$, so noting that $\log\|h\| \simeq \dist(y, x_{0})$, we get 
\begin{align} \label{normh}
\|h\|\ll\rho^{-O(\kappa)}.
\end{align}

We have a similar bound for the Mahler measure of $\gamma_{g}$. Indeed for $g\in S$, we have by  \Cref{almost-diophantine-v-1} the inequalities $\|g\|\leq M \rho^{-M\kappa}$ and  $ \| g - \sC_{h}(\gamma_{g})\| \leq M  \rho^{1-M \kappa}$, hence using \eqref{normh} we get
\[ \| \gamma_{g}\| \ll \rho^{-O(\kappa)}.\]
Moreover, $\gamma_{g}=\wt{\gamma}_{g}\bK_{\R}$ where $\wt{\gamma}_{g}\in \bG_{\Z}$ and by compactness of $\bK_{\R}$, we have  $\| \wt{\gamma}_{g}\|\ll \rho^{-O(\kappa)}$.
This allows us to control the norms and hence the absolute values at all Archimedean places of the matrix entries of the projections $\wt\gamma_{g\, \bH}$ of $\wt\gamma_g$ into simple factors $\bH \in \cF$.
To control the absolute values at non-Archimedean places, note that $\wt\gamma_g \in \M_d(\Z)$ and the projections to simple factors are all regular maps defined over some number field.
We obtain
\[\Mah(\gamma_{g}) \defeq \max_{\bH\in \cF_{nc}}\Mah(\wt{\gamma}_{g \bH})\ll \rho^{-O(\kappa)}.\]

By taking $C\lll_{\kappa}1$, all these bounds allow to apply \Cref{Ch->h} and this concludes the proof (after choosing $\kappa\lll_{\eps}1$). 
\end{proof}

\subsection{Positive dimension when \texorpdfstring{$\Gamma_{\mu}$}{Γ\_μ} has finite orbits} \label{Sec-case-finite-orbits}

Keep \Cref{not-sec-dimpos}. We show \Cref{thm-positive-dimension} under the additional assumption that $\Gamma_{\mu}$ has a finite orbit on $X$.

\begin{lemma}
\label{reduction-commensurable}
We may assume that $\Gamma_{\mu}$ is commensurable to $\Gamma_{\mu}\cap \Lambda$, i.e. $\Gamma_\mu \cap \Lambda$ has finite index in $\Gamma_\mu$.
\end{lemma}

\begin{proof}
The assumption that $\Gamma_{\mu}$ has a finite orbit means that some conjugate $\Gamma'_{\mu} \defeq \sC_{g}(\Gamma_{\mu})$ of $\Gamma_{\mu}$ in $G$ is commensurable to $\Gamma'_{\mu}\cap \Lambda$. Such $\Gamma'_{\mu}$  has algebraic entries.
Indeed, it normalizes a finite index subgroup of $\Gamma'_{\mu}\cap \Lambda$, which implies that elements of $\{\sC_{\gamma} : \gamma\in \Gamma'_{\mu}\}$ seen as endomorphisms of $\Span_{\C} \bG_{nc}\subseteq \M_{d}(\C)$ have algebraic entries, so $\Gamma'_{\mu}$ does as well because $\bG_{nc}$ is centerless. 
Note also the bi-Lipschitz  map $G/\Lambda \to G/\Lambda$, $x \mapsto gx$ is $G$-equivariant if we let $G$ act on the second $G/\Lambda$ by $(h,x) \mapsto \sC_g(h)x$.
Thus, it suffices to show  \Cref{thm-positive-dimension} for the image measure $\sC_{g \star}\mu$ in order to get it for $\mu$. 
\end{proof}

\begin{remark} In the case  where $\Gamma_{\mu}$ has a finite orbit on $X$ but not all entries are algebraic, the above reduction still provides $\Gamma'_{\mu}$ with algebraic entries. Hence we may remove the algebraic condition on $\Gamma_{\mu}$ in this case. 
\end{remark}
\bigskip

For the remainder of \Cref{Sec-case-finite-orbits}, we assume that $\Gamma_{\mu}$ is commensurable to $\Gamma_{\mu}\cap \Lambda$. 
The next lemma describes the finite orbits of $\Gamma_{\mu}$ (equivalently of $\Lambda$). We set  
\[W_{\mu, Q} \defeq \{x\in X\,:\, \sharp \Gamma_{\mu}.x\leq Q\}\]
and recall the sets  $X_{\Q, Q}$, $D_{Q}$ have been defined in \Cref{Sec-review-arith}.

\begin{lemma}[finite orbits] \label{finite-orbits} Assume that $\Gamma_{\mu}$ is commensurable to $\Gamma_{\mu}\cap \Lambda$.
Let $M\ggg1$ and $Q\geq 2$. Then the sets $X_{\Q, Q}$, $D_{Q}$ and $W_{\mu, Q}$ are all included in one another up to replacing $Q$ by $Q^M$ (i.e. $W_{\mu, Q}\subseteq X_{\Q, Q^M}$ etc.).
\end{lemma}

\begin{proof} 
The inclusion $X_{\Q, Q}\subseteq D_{Q}$ follows from the inequality $\Mah(\alpha)\geq \den(\alpha)$ established in \Cref{Facts-Mahler}.

The polynomial separation of $D_{Q}$ from \Cref{DQ-sep} together with the fact that $X$ has finite volume implies $\sharp D_{Q}\leq Q^M$.
Combined with the observation that $\Lambda D_{Q}\subseteq D_{O(Q)}$, we get $D_{Q}\subseteq W_{\mu, Q^M}$ up to increasing $M$ depending on the index of $\Gamma_{\mu} \cap \Lambda$ in $\Gamma_\mu$.

It remains to check $W_{\mu, Q}\subseteq X_{\Q, Q^M}$.  Let $\eps\in (0, 1/2)$ be a parameter. Given $y\in W_{\mu, Q}$, we have 
\[\limsup_{n \to +\infty} \max_{x \in W_{\mu,Q}} \mu^{*n}*\delta_{x}(B_{\rho}(y))\geq 1/Q\]
for all $\rho>0$.
Applying \Cref{nonpos->rational} to $\L = \Q$, we get $y\in  \bigcap_{\rho < Q^{-A_{\eps}}} X_{\Q, \rho^{-\eps}}^{(\rho^{1-\eps})}$ where $A_{\eps}>1$ is a constant depending on the working framework \ref{not-sec-dimpos} and $\eps$.
But the polynomial separation of $X_{\Q, Q}$ from \Cref{separation-XK} implies that for every $\rho_{1}>\rho_{2}>0$ and $r_{1}>r_{2}$ in $(0, \rho_{2}^{(1-\eps)/2}]$, 
\[X_{\Q, \rho_{1}^{-\eps}}^{(r_{1})}\cap X_{\Q, \rho_{2}^{-\eps}}^{(r_{2})}=X_{\Q, \rho_{1}^{-\eps}}^{(r_{2})}.\]
as long as $\eps, \rho_{1} \lll 1$. Fixing such $\eps$ and using induction, we deduce
\[\bigcap_{\rho < Q^{-A_{\eps}}} X_{\Q, \rho^{-\eps}}^{(\rho^{1-\eps})}\subseteq X_{\Q, Q^M} \]
where $M$ only depends on the initial data of \Cref{not-sec-dimpos} (assuming $\eps$ does). Hence $y\in X_{\Q, Q^M}$ and this finishes the proof. 
\end{proof}

In view of \Cref{finite-orbits},  the almost Diophantine property \ref{nonpos->rational}  tells us that if a ball contradicts positive dimension, then it must be located very close to a finite $\Gamma_{\mu}$-orbit of small cardinality. We now prove that the only way this may happen is if the walk initially started extremely close from this small orbit. This is the classical argument that the walk is repelled by finite orbits due to the positivity of the Lyapunov exponent for the adjoint random walk.

In order to control the distance to finite orbits, given $Q\geq 1$, we set
\[u_{Q}(x) \defeq    \dist(x, W_{\mu, Q})^{-1} .\]
We recall that $P_{\mu}$ denotes the Markov operator associated to $\mu$, acting on functions $f : X \rightarrow[0, +\infty]$ via the formula $P_{\mu}f : x\mapsto \int_{G} f(gx) \dd\mu(g).$ 
Finally, we denote by $\lambda_{\mu}$  the minimum of the top Lyapunov exponents for the $\mu$-walk in the irreducible components of the adjoint representation. 
A theorem due to Furstenberg~\cite{Furstenberg} tells us that $\lambda_{\mu}>0$ because $G$ has no compact factor.

The next proposition claims that the functions $\{u_{Q}, \,Q\geq 1\}$ admit a small common power that is contracted under $P_{\mu}$ up to a polynomial additive constant.

\begin{proposition}[Drift away from finite orbits] \label{drift-finite-orbits}
Let $\lambda \in (0, \lambda_{\mu})$.  There exist $C_{0}, s_{0}>0$  such that for all $n,Q \geq 1$, $x\in X$, 
\[C_{0}^{-1} P_{\mu}^n u_{Q}^{ s_{0}}(x) \leq e^{- s_{0} \lambda n}  u_{Q}^{s_{0}}(x)+  Q^{C_{0}}.\]
\end{proposition}

We stress the fact that the contraction factor does not depend on $Q$ and the additive error is polynomial in $Q$.
This will be crucial for our arguments.
Proof of \Cref{drift-finite-orbits} is essentially already present in the literature \cite{BQ2, BenDeS21}, see also \cite[Lemma 6.6]{HS2023}.
We give the proof for the convenience of the reader.  
To this end, it is useful to consider a slight reformulation of the drift functions.  

The combination of \Cref{DQ-sep} and \Cref{finite-orbits} tells us that there exists $M>1$ depending only on the initial data  \ref{not-sec-dimpos} such that the subset $W_{\mu,Q}$ is $2Q^{-M}$-separated and included in $\{\inj \geq 2Q^{-M} \}$ for all $Q\geq 2$. In particular, writing $W_{\mu,Q}^{(r)}$ the $r$-neighborhood of $W_{\mu,Q}$, 
one has that  every $x\in W^{(Q^{-M})}_{\mu, Q}$ can be uniquely written as $x=\exp(v_{x})y_{x}$ where  $y_{x}\in W_{\mu, Q}$ minimizes the distance to $x$ and $v_{x}\in \kg$ has minimal norm. 
We let $\wt{u}_{Q}$ be the function defined by
\[
\wt{u}_{Q}(x)= \left\{
    \begin{array}{ll}
       \|v_{x}\|^{-1} & \mbox{if } x\in W^{(Q^{-M})}_{\mu, Q},\\
        0 & \mbox{otherwise.}
    \end{array}
\right.
\]
Up to assuming $M\ggg 1$,  we have $u_{Q}(x)- 2Q^M \ll \wt{u}_{Q}(x)\ll u_{Q}(x)$. 

\bigskip

\begin{proof}[Proof of \Cref{drift-finite-orbits}] 
Note that we may assume $Q\geq 2$ and replace $u_{Q}$ by $\wt{u}_{Q}$ without loss of generality.  
By compactness of the support of $\mu$, we may  fix $L_{\mu}>1$ such that for every $g\in \supp \mu$ is $L_{\mu}$\nobreakdash-bi-Lipschitz on $X$.
Let $k\geq 1$ be a parameter to specify below. Let  $\cV(k)$ be  the $L_{\mu}^{-k}Q^{-M}$-neighborhood of $W_{\mu, Q}$. Then for all $x\in \cV(k)$, $g\in \supp \mu^{*k}$, we have $gx\in W^{(Q^{-M})}_{\mu, Q}$ and 
\[ v_{gx} =\Ad(g)v_{x}.\]
It follows that  for  $s>0$, $x\in \cV(k)\smallsetminus W_{\mu, Q}$, writing $w_{x}=v_{x}/\|v_{x}\|$, we have
\begin{align*}
  \frac{(P_{\mu}^k \wt{u}^s_{Q})(x)}{\wt{u}^s_{Q}(x)}
& =\int_{G}      \|\Ad(g)w_{x}\|^{-s} \dd\mu^{*k}(g)\\
&=\int_{G}    \exp(-s \log \|\Ad(g) w_{x}\|) \dd \mu^{*k}(g)\\
&\leq  1 -s  \int_{G}  \log \|\Ad(g) w_{x}\| \dd\mu^{*k}(g) + s^2 \int_{G} \left(  \log \|\Ad(g)\|\right)^2   \dd\mu^{*k}(g)
\end{align*}
where the last line uses that $e^t\leq 1+t +t^2$ for all $t\in [-1, 1]$ and assumes that $s$ is small enough in terms of $k$ and $\supp \mu$ so that the term in the exponential is in $[-1, 1]$. We may choose $k=k_{0}$ depending only on $\mu$ and the norm $\|\cdot\|$ on $\kg$ such that $\int_{G}  \log \|\Ad(g) w_{x}\| \dd\mu^{*k_{0}}(g)\geq \frac{k_{0}}{2}(\lambda+\lambda_{\mu})$, thanks to the law of large number for the norm cocycle (\cite{Furstenberg} or \cite[Theorem 4.28 (d)]{BQ_book} applied to each irreducible component of $\kg$).
Then we can choose $s=s_{0}$ depending on $(\mu, \|\cdot\|, k_{0})$ such that for all $x\in \cV(k_{0})$, 
\begin{align*}
P_{\mu}^{k_{0}}  \wt{u}^{s_{0}}_{Q}(x) \leq e^{-{s_{0}} k_{0}\lambda}  \wt{u}^{{s_{0}}}_{Q}(x)
\end{align*}
Note that for $x\notin \cV(k_0)$, for any  $g\in \supp \mu^{*k_{0}}$, we  have $\dist(gx, W_{\mu, Q})\geq L_{\mu}^{-2k_{0}} Q^{-M}$, so  $ \wt{u}_{Q}(gx)\leq 2 L_{\mu}^{2k_{0}} Q^{M}$. 

Combining these estimates,  we get that for all $x\in X$, 
\begin{align*}
P_{\mu}^{k_{0}}  \wt{u}^{s_{0}}_{Q}(x) \leq e^{-{s_{0}} k_{0}\lambda}  \wt{u}^{{s_{0}}}_{Q}(x) +2^{s_{0}} L^{2k_{0}s_{0}}_{\mu} Q^{M{s_{0}}}
\end{align*}
It follows that for all $n\in k_{0}\N^*$, 
\begin{align*}
P_{\mu}^{n}  \wt{u}^{s_{0}}_{Q}(x) \leq e^{-{s_{0}} n\lambda}  \wt{u}^{{s_{0}}}_{Q}(x) +C Q^{M{s_{0}}}
\end{align*}
where $C=2^{s_{0}} L^{2k_{0}s_{0}}_{\mu}  (1-e^{-{s_{0}} k_{0}\lambda})^{-1}$. Now if $n\in \N^*$ is of the form $n=k_{0}p+q$ where $p, q\in \N$  and $0\leq q< k_{0}$, then 
\begin{align*}
P_{\mu}^{n}  \wt{u}^{s_{0}}_{Q}(x) \leq e^{-{s_{0}} k_{0}p\lambda} P_{\mu}^q  \wt{u}^{{s_{0}}}_{Q}(x) +C Q^{M{s_{0}}}
\end{align*}
and the result follows because $ \wt{u}(gx)/ \wt{u}(x)$ is bounded away from zero and infinity uniformly in $x\in X$, $g\in \supp \mu^{*q}$. 
\end{proof}

We are now able to conclude this subsection by proving the effective positive dimension in the case where $\Gamma_{\mu}$ has a finite orbit. 

\begin{lemma}  \label{dimpos-fo}
Assume that $\Gamma_{\mu}$ has a finite orbit in $X$. Then \Cref{thm-positive-dimension}  holds. 
\end{lemma}

\begin{proof} 
We may restrict to \Cref{not-sec-dimpos} and assume that $\Gamma_{\mu}$ is commensurable to $\Gamma_{\mu}\cap \Lambda$ (\Cref{reduction-commensurable}). 
Note also that the lower bound on $n$ in \Cref{thm-positive-dimension} is decreasing as a function of $\rho_0$. Hence, if this lower bound is satisfied for $\rho_0$, then it also holds for all $\rho \geq \rho_0$. In particular, it is enough to check \Cref{thm-positive-dimension}  for $\rho_0=\rho$. 

Let  $M>1$ be as in \Cref{finite-orbits},  let $\eps\in (0, 1/2)$, and $(A, C, \kappa)$ so that \Cref{nonpos->rational} holds with $\eps' \defeq \eps/M$. Let $\rho\leq 1/C$, $x,y\in X$, $n\geq \eps |\log \rho|+ A\,\dist(x, x_{0})$. To obtain a contradiction, suppose that 
\[\mu^{*n}*\delta_{x}(B_{\rho}(y))> \rho^\kappa. \]
By  \Cref{nonpos->rational} and \Cref{finite-orbits}, we must have 
\[y\in X^{(\rho^{1-\eps'})}_{\Q, \rho^{-\eps'}} \subseteq W^{(\rho^{1-\eps'})}_{\mu, \rho^{-\eps}}.\]
In particular, 
\[\mu^{*n}*\delta_{x}(W^{(\rho^{1-\eps})}_{\mu, \rho^{-\eps}})> \rho^\kappa. \]
Considering the drift function $u_{\rho^{-\eps}}$, we obtain for any $s>0$, 
\[P_\mu^n u^s_{\rho^{-\eps}}(x)= \mu^{*n}*\delta_{x}(u^{s}_{\rho^{-\eps}}) \geq \rho^{\kappa-s(1-\eps)}.\]
Fixing $\lambda\in (0, \lambda_{\mu})$ then choosing $s=s_{0}$ as in \Cref{drift-finite-orbits}, we may bound the left hand side to obtain 
\[C_{0}e^{-s_{0} n\lambda}u^{s_{0}}_{\rho^{-\eps}}(x)+\rho^{-\eps C_{0}} \geq \rho^{\kappa-s_{0}(1-\eps)}.\]
Imposing additionally $n\geq \frac{1}{\lambda} |\log \dist(x, W_{\mu, \rho^{-\eps}})|$, we get
\[C_{0}+ \rho^{-\eps C_{0}} \geq  \rho^{\kappa-s_{0}(1-\eps)}.\]
Taking $\eps$ and $\kappa$ such that $\eps C_{0}, \eps, \kappa<s_{0}/4$ and assuming $\rho$ small enough in terms of $C_{0}, s_{0}$, we get a contradiction. 
\end{proof}

\subsection{Positive dimension when \texorpdfstring{$\Gamma_{\mu}$}{Γ\_μ} has no finite orbit}

We show \Cref{thm-positive-dimension} in the case where $\Gamma_{\mu}$ has no finite orbit on $X$.  The first step is to show that the $\mu$-walk on $X$ gives exponentially small mass to algebraic points of small Mahler measure, see \Cref{Sec-pointwise-pos-dim}. Once this is established, we conclude in \Cref{Sec-pointwise-to-posdim}  using the almost Diophantine property (\ref{nonpos->rational}) and effective recurrence away from infinity  (\ref{effective-recurrence}).
Notations refer to \Cref{not-sec-dimpos}.

\subsubsection{Exponentially small probability to hit an algebraic point.} \label{Sec-pointwise-pos-dim} 

\Cref{Sec-pointwise-pos-dim} is dedicated to proving the following statement. 
We let $\L$ be a finite Galois extension of $\Q$  such that $\Gamma_{\mu}\subseteq G_{\L}$ and all factors of $\bG$ are defined over $\L$.   

\begin{proposition}[Pointwise decay] \label{pointwise-pos-dim}
Assume $\Gamma_{\mu}$ has no finite orbit on $X$. There exists $r>0$ such that for all  large enough $n$, all $x\in X$, $y\in X_{\L, e^{rn}}$,
\[\mu^{*n}*\delta_{x} (\{y\}) \leq e^{-rn}. \]
\end{proposition}

We first give some insight about the proof. Consider momentarily the case $G=\SL_{2}(\R)$, $\Lambda=\SL_{2}(\Z)$ with $\mu$ supported on $\SL_{2}(\R\cap \QQ)$. Essentially, the problem is to show that $\mu^{*n}(\Lambda)$ decays exponentially. In case $\mu$ is supported on $\SL_{2}(\Q)$, the absence of a finite $\Gamma_{\mu}$-orbit allows to find a suitable $p$-adic embedding in which $\Gamma_{\mu}$ is unbounded,  then we conclude observing that the $\mu$-walk has positive Lyapunov exponent in this embedding while $\Lambda$ is represented by a bounded set. In case that $\mu$ is supported on $\SL_{2}(\Z[\sqrt{2}])$ this argument does not work because $\Gamma_{\mu}$ is bounded for any non-Archimedean  norm. The idea is to use the restriction of scalars to identify $\Gamma_{\mu}$ with a Zariski-dense subgroup of a bigger group $\bR$ in which $\Lambda$ is represented by a proper algebraic subgroup. Then we conclude using the transience of proper algebraic subgroups. 

\bigskip

We now introduce the appropriate restriction of scalars, motivated by the previous paragraph and \Cref{L-points-G}.  Let $E$ be the set of couples $(\bH, \sigma)$ where $\bH\in \cF_{nc}$ and $\sigma\in \Gal(\L/ \Q)$ are such that $\presigma{\bH}\in \cF_{nc}$. Set 
\[
\Theta : G_{\L} \, \hookrightarrow \,\prod_{(\bH, \sigma)\in E}\presigma{\bH}_{\L},\,\,\,\, g\mapsto  (\presigma{g_{\bH}})_{(\bH, \sigma)\in E}.\]
Denote by $\bR$  the Zariski closure of $\Theta(\Gamma_{\mu})$ in $\prod_{(\bH, \sigma)\in E}\presigma{\bH}$.
 
\begin{lemma}\label{Hsemisimple}
$\bR$ is a semisimple algebraic group. 
\end{lemma}

\begin{proof}
We show that $\bR$ is reductive with trivial center. 

For each $(\bH, \sigma)\in E$, the standard representation $\presigma{\bH}_{\C} \acts \C^d$ is completely reducible by simplicity of $\presigma{\bH}_{\C}$. Observe that $\bR_{\C}$ acts faithfully on $\oplus_{(\bH, \sigma)\in E}\C^d$.
Moreover, each subspace $\C^d$ is preserved by $\bR_{\C}$, the restricted action factorizes through $\presigma{\bH}_{\C}$ and every element of $\presigma{\bH}_{\C}$ is represented by the Zariski-density of $\Gamma_{\mu}$. Hence $\bR_{\C}\acts \oplus_{E}\C^d$ is completely reducible and faithful, which forces $\bR_{\C}$ to be reductive.  

 Assume $h=(h_{(\bH, \sigma)})$ belongs to the  center of $\bR_{\C}$. Then for any $(\bH, \sigma)\in E$, we have for all $g\in \Gamma_{\mu}$ that $h_{(\bH, \sigma)}\presigma{g_{\bH}}=\presigma{g_{\bH}} h_{(\bH, \sigma)}$. As $\Gamma_{\mu}$ is Zariski-dense in $\bG_{\C}$, we get that $h_{(\bH, \sigma)}$ centralizes $\presigma{\bH}_{\C}$, hence it is trivial. 
\end{proof}

Let $\Delta=\{ a\in \prod_{(\bH, \sigma)\in E}\presigma{\bH}\,:\,  a_{(\bH, \sigma)}= a_{(\presigma{\bH}, \Id)}\}$, which, in view of \Cref{L-points-G},  corresponds to the Zariski-closure of $\Theta(G_{\Q})$.  
We will distinguish the two cases listed in \Cref{list-cases} below. 

\begin{lemma}[List of cases] \label{list-cases} At least one of the following holds. 

\begin{itemize}

\item[(i)] There exists $ g\in G_{\L}$ such that $\Gamma_{\mu} \cap g G_{\Q} g^{-1}$ has finite index in $\Gamma_{\mu}$.

\item[(ii)] For all $g\in G_{\L}$ the intersection  $\bR\cap \Theta(g)\Delta \Theta(g)^{-1}$ is an algebraic subgroup of dimension strictly less than $\dim \bR$. 

\end{itemize}
\end{lemma}

\begin{proof}
We assume that (ii) fails and we show (i). By assumption, there is $g\in G_{\L}$ such that $\Theta(g)\Delta \Theta(g)^{-1}$ contains the identity component of $\bR$.  Hence a finite index subgroup $\Gamma' \subseteq \Gamma_{\mu}$ satisfies 
\[ \Theta(\Gamma')\subseteq \Theta(g)\Delta \Theta(g)^{-1}\]
i.e.
\[ \Theta(g^{-1}\Gamma'g)\subseteq \Delta \]
and by \Cref{L-points-G}, this means
\[g^{-1}\Gamma'g\subseteq G_{\Q}.\qedhere\]
\end{proof}


We now establish \Cref{pointwise-pos-dim}, starting with case (ii) which is easier. 

\begin{lemma}[Pointwise decay - case (ii)] \label{pointwise-pos-dim-ii}
\Cref{pointwise-pos-dim} holds in case (ii) of \Cref{list-cases}.  
\end{lemma}

\begin{proof} First note that
\[\mu^{*n}*\delta_{x} (\{y\})^2 \leq (\mu^{*n}\otimes \mu^{*n})\set{(g_{1},g_{2}) : g_{1} g_{2}^{-1}y = y}.\]
Writing $y=h_{y}\Lambda$ with $h_y \in G_\L$, the equality $g_{1} g_{2}^{-1}y=y$ means $h_{y}^{-1}g_{1} g_{2}^{-1} h_{y} \in \Lambda$, which implies 
$\Theta(h_{y}^{-1}g_{1} g_{2}^{-1} h_{y}) \in\Delta$, 
i.e.
$\Theta(g_{1}g_{2}^{-1})  \in \Theta(h_{y}) \Delta \Theta(h_{y})^{-1}. $
But $\Theta(g_{1}g_{2}^{-1})$ also belongs to $\bR$, hence we get 
\[\Theta(g_{1}) \in \bR'_{y}\Theta(g_{2})\]
where $\bR'_{y}=\Theta(h_{y}) \Delta \Theta(h_{y})^{-1} \cap  \bR$ is an algebraic subgroup of $\bR$ with strictly smaller dimension. 
In the end, 
\[\mu^{*n}*\delta_{x} (\{y\})^2 \leq  \sup_{\dim \bR' <\dim \bR,  h\in \bR}(\Theta_{\star}\mu)^{*n}(\bR'h)\]
and then exponential decay follows from Breuillard-Becker \cite[Proposition 3.3, Theorem 1.3]{BreBeck25} (see also \cite[Proposition 1.2, Lemma 6.1, Lemma 6.2]{Breuillard-note}) which applies\footnote{Note that  \cite{BreBeck25} assumes the ambient group to be Zariski-connected. We may nevertheless apply their result by using the following observation. Let $G$ be a topological  group with finitely many connected components, $\eta$ a probability measure on $G$,  for $\ug\in G^{\N^*}$ distributed by $\eta^{\otimes \N^*}$ set $\tau_{0} \equiv 0$ and  $\tau_{n}(\ug) \defeq \inf\{i >\tau_{n-1}(\ug) : g_{i}\dots g_{1}\in G^{\circ}\}$ the $n$-th return time to the identity component $G^\circ$ of $G$, set $\eta_{\tau_{1}}$ the law of $g_{\tau_{1}}\dots g_{1}$, driving the induced random walk on $G^{\circ}$. Then for any  $\eps>0$, any $E\subseteq G^{\circ}$, one has $ \eta^{*n}(E)\leq  \mu^{\otimes \N^*} \{\tau_{\eps n}>n \} + \sum_{k\geq \eps n} \eta_{\tau_{1}}^{*k} (E)$. The first term of the sum has exponential decay (provided $\eps$ is small), hence exponential decay for $\eta^{*n}(E)$ boils down to that of $\eta_{\tau_{1}}^{*n} (E)$. Then in our context,  \cite{BreBeck25} applies to $\eta_{\tau_{1}}$.} 
thanks to  \Cref{Hsemisimple}. 
\end{proof}

We now address case (i).
Observe that our assumption on $\L$ allows to identify $G_{\L}$ with a  subgroup of $\prod_{\bH\in \cF_{nc}}\bH_{\L}$, in particular $G_{\L}\subseteq \SL_{d}(\L)$. Also $\Lambda\subseteq \SL_{d}(\frac{1}{M}\cO_{\L})$ for some $M \geq 1$ because each $\bH$-coordinate of a $\Z$-point of $\bG$ has bounded denominator.

\begin{lemma} \label{unbounded-embedding}
Let $\Gamma\subseteq G_{\L}$ be a finitely generated subgroup without finite orbit  on $X$.  If $\Gamma\cap G_{\Q}$ has   finite index in $\Gamma$, then there exists a prime number $p\geq 2$, a finite extension $\L'$ of $\Q_{p}$, and a field embedding $\psi \colon \L \hookrightarrow \L'$ such that $\psi(\Gamma)$ is unbounded in $\SL_{d}(\L')$. 
\end{lemma}

\begin{proof}
We just need to construct a non-Archimedean norm on $\L$ for which $\Gamma$ is unbounded. By Ostrowski's theorem, the prime ideals of $\cO_{\L}$ give a system of representatives for the non-Archimedean norms on $\L$ up to a power.
 More precisely, if $\kP\subseteq \cO_{\L}$ is a prime ideal, then we define for $x\in \L$ the corresponding norm $\|x\|_{\kP}=N(\kP)^{-v_{\kP}(x)}$ where $N(\kP)=\sharp \cO_{\L}/\kP$ and $v_{\kP}(x)\in \Z$ is the number of occurences of $\kP$ in the decomposition of the fractional ideal $x\cO_{\L}$ into prime ideals or their inverses. 

Now write $E_{\Gamma}\subseteq \L$ the set of entries of elements in $\Gamma$. As $\Gamma$ is finitely generated, only a finite collection of prime ideals $\kP_{1}, \dotsc, \kP_{s}$ may satisfy $\inf_{x\in E_{\Gamma}}v_{\kP_{i}}(x)<0$. Assume by contradiction that $\Gamma$ is bounded  with respect to all the norms $\|\cdot\|_{\kP_{i}}$. Then any $x\in E_{\Gamma}$ satisfies 
$x\cO_{\L}= \mathfrak{J}\kP^{i_{1}}_{1} \dots \kP^{i_{s}}_{s}$ where $ \mathfrak{J}\subseteq \cO_{\L}$ is an ideal and the $i_{k}$ are negative integers such that $\sup_{k} \abs{i_{k}}\leq M$ for a constant $M$ depending only on $(\L,\Gamma)$. Choosing an arbitrary element  $a_{i}\in\kP_{i}$  for each $i$, we then have 
\[x \in \frac{1}{a_{1}^{M}\dotsm a_{s}^{M}} \cO_{\L}.\]
In particular, $\sup_{x\in E_{\Gamma}} \den(x)<\infty$, so $\sup_{g\in \Gamma } \den(g)<\infty$. The orbit $(\Gamma\cap G_{\Q}) x_{0}$ (where $x_{0}=\Lambda /\Lambda\in X$) is therefore included in some set $D_{Q}$ where $Q\geq1$, and in particular  finite (\Cref{DQ-sep}).
Recalling the assumption $[\Gamma:\Gamma\cap G_{\Q}]<\infty$, we deduce that $x_{0}$ has a finite $\Gamma$-orbit. This contradicts our hypothesis.
\end{proof}

\begin{lemma}[Pointwise decay - case (i)] \label{pointwise-pos-dim-i} 
\Cref{pointwise-pos-dim} holds under the additional assumption that there exists $g\in G_{\L}$ such that $\Gamma_{\mu}\cap gG_{\Q}g^{-1}$ has finite index in $\Gamma_{\mu}$.
\end{lemma}

\begin{proof}  We may assume $\Gamma_{\mu}$ intersects $G_{\Q}$  in a finite index subgroup (without conjugation required).


 Let $r>0$ be a parameter to specify below. 
If $\mu^{*n}*\delta_{x} (\{y\}) > e^{-rn} $, then 
\[(\mu^{*n} \otimes \mu^{*n}) \set{(g_{1},g_{2})\,:\, g_{1}g_{2}^{-1}y=y }> e^{-2rn}.\]
The equality $g_{1}g_{2}^{-1}y=y$ can be rewritten $g_{1}g_{2}^{-1}\in h_{y}\Lambda h_{y}^{-1}$ where $h_{y}\in G_{\L}$ is such that $y=h_{y}\Lambda$ and $\Mah(h_{y})\leq e^{rn}$.  Letting $(\L', \psi)$ as in \Cref{unbounded-embedding} applied to $\Gamma_{\mu}$ and observing that $\sup_{\gamma \in \Lambda}\norm{\gamma} \leq 1$ for the non-Archimedean norm $\norm{\cdot}$, we obtain that
\begin{equation}\label{anti-ldp}
(\mu^{*n} \otimes \mu^{*n})\setbig{ (g_{1},g_{2}) : \norm{\psi(g_{1}g_{2}^{-1})} \leq C_{\psi} e^{C_{\psi} rn}  } > e^{-2rn} 
\end{equation}
for some constant $C_{\psi}>1$ depending on $\psi$.
But $\psi(\Gamma_{\mu})$ is unbounded in $\SL_{d}(\L')$ with semisimple Zariski-closure. 
The positivity of the top Lyapunov exponent~\cite[Theorem 10.9 (e)]{BQ_book} and the large deviation estimate~\cite[Theorem 13.11 (iii)]{BQ_book} for the associated $\mu$-walk  yields a contradiction with \eqref{anti-ldp} for small enough $r$  and large enough $n$. 
\end{proof}

\begin{proof}[Proof of \Cref{pointwise-pos-dim}]
 Combine Lemmas \ref{list-cases}, \ref{pointwise-pos-dim-i}, \ref{pointwise-pos-dim-ii}.   
\end{proof}

\subsubsection{From pointwise decay to positive dimension} \label{Sec-pointwise-to-posdim}

We combine the pointwise exponential decay  (\Cref{pointwise-pos-dim}) and the almost Diophantine property (\Cref{nonpos->rational}) to prove that the $\mu$-walk generates effective positive dimension at a single exponentially small scale. 

\begin{lemma}[Positive dimension at single scale] \label{single-scale-posdim} 
Assume $\Gamma_{\mu}$ has no finite orbit on $X$,  set $x_{0}=\Lambda/\Lambda\in X$. There exist $A, M, r>0$ such that for all $x,y\in X$, $n\geq A\cdot\dist(x, x_{0}) + M$. 
\[\mu^{*n}*\delta_{x} (B_{e^{-Mn}}(y)) \leq e^{-rn}. \]
\end{lemma}

The idea of the proof is that if the walk starting at $x$ gives a lot of mass to a ball of exponentially small radius, then by \Cref{nonpos->rational} the center $y$ of the ball is close to an algebraic point $y'$ with small Mahler measure, hence walking backward the same holds  for $x$ with respect to some algebraic point $x'$, and by effective separation of algebraic points with small complexity, the walk starting at $x'$ has a large probability to reach $y'$ in $n$ steps, which contradicts \Cref{pointwise-pos-dim}.

\begin{proof}
We let $A>0$ be as in \Cref{nonpos->rational}.
Let $M, n_{0}\geq 1$, $\eps, r\in (0, 1/2)$ be  parameters to specify below.
Let $(C,\kappa)$ be the constants associated to $\eps$ by  \Cref{nonpos->rational}.  We  take $n_{0}$ large enough so that $e^{-n_{0}}\leq 1/C$ and we also assume $r<\kappa$. We consider  $x,y\in X$, $n\geq 2A\cdot\dist(x, x_{0}) + n_{0}$ and we note that, up to assuming $\eps M<1/4$,  the condition on $n$ is sufficient to apply  \Cref{nonpos->rational} with $\rho=e^{-Mn}$. 
Suppose that 
\[\mu^{*n}*\delta_{x} (B_{e^{-Mn}}(y)) > e^{-rn}.\]
Then by \Cref{nonpos->rational}, denoting by $\L$  the smallest (finite) Galois extension of $\Q$ on which the factors of $\bG$ are defined and such that $\Gamma_{\mu}\subseteq G_{\L}$, we have 
\[y\in X^{(e^{-(1-\eps)Mn})}_{\L, e^{\eps M n}}.\]
Walking backward in time and recalling \Cref{Facts-Mahler} (i), we deduce that 
\[x\in X^{(e^{-(1-\eps)Mn}L_{\mu}^n)}_{\L, e^{O(\eps M n)}Q_{\mu}^n} \]
where  $L_{\mu},Q_{\mu}>1$ are constants depending only on $X$ and the support of $\mu$. 
Let $x'\in X_{\L, e^{O(\eps M n)}Q_{\mu}^n}$ and $y' \in X_{\L, e^{\eps M n}} $ be the algebraic points minimizing $\dist(x,x')$ and $\dist(y,y')$.
We deduce that for a $\mu^{*n}$-proportion at least $e^{-rn}$ of $g\in G$, we have 
\[\dist(gx', y')\leq \dist(gx', gx)+\dist(gx,y)+\dist(y,y')\leq 3e^{-(1-\eps)Mn}L_{\mu}^{2n}. \]
 But $gx', y'\in X_{\L, e^{O(\eps M n)}Q_{\mu}^{2n}}$ and this set $(e^{-\eps M n}Q_{\mu}^{-2n})^{O(1)}$-separated by \Cref{separation-XK}.  Hence if $\eps^{-1},M$ are large enough in terms of $(X, \mu)$, and $n_{0}$ large enough depending on $(X, \mu, M)$,  we must have $gx'=y'$. This leads to
\[\mu^{*n}*\delta_{x'} (\{y'\}) > e^{-rn}.\]
Up to assuming that $r$ is small enough in terms of $(X, \mu)$, this contradicts \Cref{pointwise-pos-dim}. This concludes the proof (up to replacing $M$ by $\max(n_{0}, M)$). 
\end{proof}

By cutting the walk trajectory into two suitable pieces and applying the single scale estimate of \Cref{single-scale-posdim}, we finally obtain an effective positive dimension at all scales above an exponentially decreasing threshold. 

\begin{lemma}
Assume $\Gamma_{\mu}$ has no finite orbit on $X$. Then \Cref{thm-positive-dimension} holds.
\end{lemma}

\begin{proof} Note that the lower bound on $n$ is decreasing as a function of $\rho$, in particular it is enough to check \eqref{eq-posdim-c} for $\rho_0=\rho$. 
We let $A, M, r$ be as in \Cref{single-scale-posdim}. 
Given $\rho\in(0,1)$ we write $n_{\rho}\defeq \lfloor \frac{1}{M} \abs{\log \rho} \rfloor$ the largest integer such that $\rho\leq e^{-M n_{\rho}}$. For $x,y\in X$, $n\geq n_{\rho}$, we have 
\begin{align} \label{eq-tpd-2}
\mu^{*n}*\delta_{x} (B_{\rho}(y)) &=\int_{G} \mu^{*n_{\rho}}*\delta_{gx} (B_{\rho}(y)) \dd\mu^{*(n-n_{\rho})}(g) \nonumber\\
&\leq 3\rho^{r/M} + \mu^{*(n-n_{\rho})} \{\, g \,:\, n_{\rho} < A\cdot\dist(gx, x_{0})+M \,\}
\end{align}
where the upper bound uses \Cref{single-scale-posdim}. Assuming $n-n_{\rho}\ggg \dist(x, x_{0})$ we may apply \Cref{effective-recurrence} to bound the second term 
\[ \mu^{*(n-n_{\rho})} \{ g \,:\, n_{\rho} < A\cdot\dist(gx, x_{0})+M\} \ll e^{-s_{0}(n_{\rho}-M)/A} \ll \rho^{\frac{s_{0}}{2MA}}\]
and this concludes the proof.
\end{proof}

\section{From small dimension to equidistribution}
\label{Sec-bootstrap}

The goal of the section is to show that a random walk on a (possibly non-arithmetic) finite-volume homogeneous space modeled over $\so(2,1)$ or $\so(3,1)$ equidistributes exponentially fast toward the Haar probability measure, provided the initial distribution of the walk has positive dimension and is not too concentrated near infinity. 

\bigskip


The rate of convergence is estimated on a space of regular functions. Given a metric space $X$ and $\b \in {(0,1]}$, we let $C^{0,\b}(X)$ denote the space of bounded $\beta$-Hölder continuous functions on $X$, endowed with its usual norm $\norm{\,\cdot\,}_{C^{0,\b}}$ :
\begin{align} \label{def-holder}
\forall f \in C^{0,\b}(X), \quad \norm{f}_{C^{0,\b}} \defeq \norm{f}_\infty + \sup_{x\neq y \in X} \frac{\abs{f(x)-f(y)}}{\dist(x,y)^\b}.
\end{align}
The corresponding Wasserstein distance between two probability measures $\nu,\nu'$ on 
X is defined as 
\begin{align} \label{def-wassertstein}
\cW_\b(\nu,\nu') \defeq \sup_{f \in C^{0,\b}(X),\norm{f}_{C^{0,\b}} \leq 1} \abse{\int_X f \dd \nu - \int_X f \dd \nu'}.
\end{align}



\begin{thm}[From small dimension to equidistribution] \label{small-dim-equid}
Let $G$ be connected real Lie group with Lie algebra $\so(2,1)$ or $\so(3,1)$. Let $\Lambda\subseteq G$ be a lattice, $X=G/\Lambda$ equipped with a quotient right-invariant Riemannian metric. Let $\mu$ be a Borel probability measure on $G$ having a finite exponential moment and whose support generates a Zariski-dense subgroup of $G$. 

Given $\beta \in {(0,1]}$ and $\kappa \in (0, 1]$, there exist $\eps = \eps(X,\mu,\beta,\kappa) > 0$  such that for small  enough $\delta>0$, the following holds.

Let $\nu$ be a probability measure on $X$ satisfying 
\[\nu (B_\rho(x)) \leq \rho^\kappa \text{ for all } x \in X, \rho \in [\delta, \delta^\eps].\]
Then for all $n\geq  |\log \delta|$, one has
\[\cW_{\beta}(\mu^{*n}*\nu, m_{X})\leq \delta^\eps+2\nu\{\inj\leq \delta^\eps\}\]
where  $m_{X}$ denotes the Haar probability measure on $X$. 
\end{thm}

\begin{proof}[Proof admitting the key steps]
Note that for parameters $0 < \beta \leq \beta' \leq 1$, any function $f \in C^{0,\beta'}(X)$ belongs also to $C^{0,\beta}(X)$ and
\(
\norm{f}_{C^{0,\beta}} \leq 2 \norm{f}_{C^{0,\beta'}}.
\)
Consequently, $\cW_{\beta'}$ is dominated by $\cW_{\beta}$.
Therefore, it suffices to show the theorem for small $\beta > 0$.
If $\beta$ is small enough so that $\int_{G}\norm{\Ad(g)}^\b \dd \mu(g) < \infty$, then we can
combine \Cref{pr:persist}, \Cref{dimension-bootstrap} and \Cref{endgame} established in the next subsections, and this gives the theorem under the condition $n\geq C \abs{\log \delta}$ where  $C>1$ is a constant depending on  $(X, \mu)$. But writing $|\log \delta|=C |\log \delta^{1/C}|$, we see that $C$ can be assumed to be equal to $1$ up to changing $\eps$ and taking a smaller upper bound for $\delta$. 
\end{proof}

Let us describe the outline of the actual proof of this theorem.
We will be working with a notion of dimension of measures at discretized scales.
One could use the notion of $(C,\alpha)$-regular measure of \cite[Definition 5.1]{BFLM}.
But since in our proof, we remove small mass from the measure while we do convolution with $\mu$, the notion of robust measures of \cite[Definition 2.5]{Shmerkin} is more suited to our need.
In \Cref{ss:robust-measure} we will recall this notion and adapt it to our problem.

Then the proof of \Cref{small-dim-equid} will be divided into three relatively independent steps.
In the first step, \Cref{ss:persist}, we show that the condition that $\nu$ has a small positive dimension at a wide range of scales is preserved under taking convolution with $\mu^{*n}$.
This is achieved by using drift functions (a.k.a. Margulis functions) similar to those used for proving recurrence in Eskin-Margulis~\cite{EskinMargulis} or Benoist-Quint~\cite{BQ_rec}.
In the second step, \Cref{ss:bootstrap}, we show that the dimension of $\nu$ increases under convolution by $\mu^{*n}$,
and this can continue until the dimension reaches any prescribed number less than $\dim X$.
In this step, the key is to apply the multislicing theorem established in \Cref{sec:slicing}. 
Finally, in \Cref{ss:high-to-equid}, we show that a probability measure whose dimension is close enough to $\dim X$ equidistributes exponentially fast under convolution by  $\mu$.
This is done by exploiting the presence of a spectral gap.
The use of spectral gap to obtain effective equidistribution is not new. For example, in the context of random walks on homogeneous spaces, this already appeared in \cite{HLL_Nil} and in \cite{KimKogler}.
\bigskip

The level of generality will depend on the subsection. It will be handy to refer to the following framework, which we will occasionally specify. 

\begin{notation} \label{notation-bootstrap}
 $G$ is a semisimple connected real Lie group with finite center and no compact factor, $\Lambda\subseteq G$ is a lattice, $X=G/\Lambda$ equipped with a quotient right-invariant Riemannian metric, $\mu$ is a Borel probability measure on $G$ having a finite exponential moment and whose support generates a Zariski-dense subgroup $\Gamma_{\mu}$ of $G$. The distance on $X$ comes from a Euclidean norm $\|\cdot\|$ on $\kg\defeq \Lie(G)$. Fixing $K\subseteq G$ a maximal compact subgroup and $\ka\subseteq \kg$ a Cartan subspace orthogonal to $\Lie(K)$ for the Killing form, we assume (as we may \cite[Lemma 6.33]{BQ_book}) that $\|\cdot\|$ is $\Ad(K)$-invariant and $\ad(\ka)$ is made of self-adjoint matrices. We choose a Weyl chamber $\ka^+$ for $\ka$ (so that the Cartan decomposition is well defined).
\end{notation}

 The symbols  $O(\cdot)$, $\ll$, $\simeq$, $\lll$ (see \ref{conventions-notations}) will refer implicitely to constants possibly depending on \Cref{notation-bootstrap}, additional dependences will be indicated in subscript.

\subsection{Robust measures} \label{ss:robust-measure}
Note that for the problem of quantitative equidistribution, in the course of the random walk, we can toss away measures of exponentially small mass which will go into the error term.
The notion of robust measure introduced by Shmerkin in~\cite{Shmerkin} turns out to be convenient.
We will adapt this notion to our problem, taking into account the possible presence of cusps in $X$. The setting refers to \ref{notation-bootstrap}. Recall that  $\inj \colon X \to \R_{>0}$ denotes the injectivity radius.


\begin{definition} \label{def-robust}
Let $\alpha \in {[0,1]}$ be a parameter and $\rho > 0$ a radius.
We say a Borel measure $\nu$ on $X$ is \emph{$(\alpha, \cB_\rho, 0)$-robust} if 
\[
\nu\{\inj < \rho \}=0 \quad \text{and} \quad \forall x \in X,\; \nu(B_\rho(x)) \leq \rho^{\alpha \dim X}.
\]

Let $I\subseteq \R_{>0}$ be a subset.
We say that $\nu$ is \emph{$(\alpha, \cB_I, 0)$-robust} if for every $\rho \in I$, 
$\nu$ is $(\alpha,\cB_\rho,0)$-robust.
This amounts to say that $\nu$ is supported on the compact $\{\inj \geq \sup I\}$ and
\[
\forall \rho \in I,\, \forall x \in X,\quad \nu(B_\rho(x)) \leq \rho^{\alpha \dim X}.
\]

Let $\tau > 0$ be another parameter, 
We say that  $\nu$ is \emph{$(\alpha, \cB_I, \tau)$-robust} if one can decompose $\nu$ into a sum of  Borel measures $\nu = \nu'+\nu''$ such that $\nu'$ is $(\alpha, \cB_{I}, 0)$-robust and $\nu''(X)\leq \tau$.
In this case, we call $\nu = \nu' + \nu''$ a \emph{robust decomposition} of $\nu$ and $\nu'$ is the \emph{regular part} of $\nu$. 
\end{definition}

Here, $\cB_I$ stands for the sets of balls of radius $\rho \in I$. 
In practice, $I$ will always be an interval, and measures $\nu$ will be finite measures of total mass at most $1$.
Roughly speaking, $\nu$ is $(\alpha, \cB_{I}, \tau)$-robust if, up to ignoring a part of $\nu$ of mass $\tau$, it is $(\alpha \dim X)$-Frostman with respect to scales appearing in $I$ and its support does not go to high in the cusp. 

\bigskip
We record some elementary combinatorial properties concerning robustness.

\begin{lemma} \label{lm:robust elem}
Let $\alpha \in {[0,1]}$ and  $\rho, \tau, \tau_1, \tau_2>0$, $I, I_1, I_2\subseteq \R_{>0}$.
\begin{itemize}
\item If $\nu$ is $(\alpha, \cB_{\rho}, \tau)$-robust then for every $r\in (0, 1)$ it is also  $(\alpha r,  \cB_{[\rho^{1/r}, \rho]}, \tau)$-robust. 
\item If $\nu$ is $(\alpha, \cB_{I_1}, \tau_1)$-robust and $(\alpha, \cB_{I_2}, \tau_2)$-robust then $\nu$ is $(\alpha, \cB_{I_1\cup I_2}, \tau_1+\tau_2)$-robust.  
\end{itemize}
\end{lemma}

\begin{proof}
The first claim is immediate. For the second, consider $\nu=\nu_{i}'+\nu_{i}''$ a robust decomposition corresponding to $I_{i}$. Write $\nu''_{i}=f_{i}\nu$ where $f_{i}$ takes values in $[0, 1]$, set $f=\max(f_{1}, f_{2})$, $\nu''=f\nu$. Then $\nu=\nu'+\nu''$ where $\nu'$ is $(\alpha, \cB_{I_1\cup I_2}, 0)$-robust, and $\nu''(X)\leq \tau_1+\tau_2$. 
\end{proof}

Note that if $\nu$ is $(\alpha, \cB_{I}, \tau)$-robust then $\nu$ is automatically $(\alpha, \cB_{\rho}, \tau)$-robust for every $\rho\in I$. We will need the following partial converse.

\begin{lemma} \label{single-vs-multiple}
Let  $\alpha, s, \delta \in (0, 1]$, $\tau \in \R^+$. If $\nu$ is $(\alpha, \cB_{\rho}, \tau)$-robust for all $\rho  \in [\delta, \delta^s]$, then for any $\eps \in (0, \alpha)$, the measure $\nu$ is $(\alpha -\eps, \cB_{[\delta, \delta^s]},  \lceil\frac{\log s}{\log (1-\eps)}\rceil \tau)$-robust.
\end{lemma}

\begin{proof} 
Let $r \defeq 1-\eps$. 
Let $k$ be the smallest integer such that $\delta^{s/r^k}\leq \delta$, namely $k=\lceil \log s/ \log r \rceil$. By \Cref{lm:robust elem}, for  each $i\in \{0, \dots, k-1 \}$, we know that $\nu$ is $(\alpha r, \cB_{[\delta^{s/r^{i+1}},\delta^{s/r^{i}}]}, \tau)$-robust. By \Cref{lm:robust elem} again, those estimates add up and the claim follows using that $\alpha r\geq \alpha -\eps$. 
\end{proof}




\subsection{Persistence of small dimension} \label{ss:persist}
In this subsection, we prove  \Cref{pr:persist}, stating that the $\alpha$-robustness of a measure on $X$ for small $\alpha> 0$  is preserved under convolution by $\mu$.
This fact will be important to initiate the bootstrap argument developed in \Cref{ss:bootstrap}. Notations refer to \ref{notation-bootstrap}.

\begin{proposition}[Persistence of small positive dimension]  \label{pr:persist} 
 Let $M \ggg 1$. 
Let $\nu$ be a $(\kappa, \cB_{I}, \tau)$-robust measure on $X$ for some $\kappa \in {(0, 1]}$, $I \subseteq \R_{>0}$ and $\tau \geq 0$, and assume $\nu(X)\leq 1$. 
Then, provided that $\sup I \lll_{\kappa} 1$, for all $n \geq 0$, the measure $\mu^{*n}*\nu$ is $(\frac{\kappa}{4M}, \cB_{I'}, \tau + \sup I )$-robust where $I'\defeq\{\rho^{M}, \rho \in I\}$.
\end{proposition}

Note that we only prove the persistence of \emph{small} positive dimension. Persistence of \emph{large} positive dimension will be a consequence of our bootstrap argument. 

\bigskip

Consider the space $X \times X$ on which $G$ acts diagonally.
By an abuse of notation, we still write $P_\mu$ to denote the Markov operator on this space associated to $\mu$. We let $u_{0}: X\rightarrow [1, +\infty)$ be the drift function from \Cref{Sec-effective-recurrence}. Up to replacing $u_{0}$ by some $R u_{0}^R$ where $R$ is large, we may assume $u_{0}^{-1}\leq \inj$ (see Lemmas \ref{comparison}, \ref{comparison2}). 
We introduce a new drift function $\omega^s_{C} : X\times X \to [0, +\infty]$ given by
\[\omega^s_C(x,y) \defeq \frac{1}{\dist(x,y)^s} + Cu_{0}^s(x)\]
where $C \geq 0$ and $s > 0$ are parameters to be adjusted.

\begin{lemma}[Drift away from the diagonal] \label{drift-function-omega} For $C\ggg 1$ and $\lambda, s \lll 1$, we have for all $n\geq 0$,
\begin{equation}\label{eq-dfo}
P_{\mu}^n\omega_{C}^s  \ll_{C} e^{-s \lambda n} \omega_{C}^s + 1.
\end{equation}
\end{lemma}


\begin{proof} First recall from \Cref{drift-infinity} that small powers of $u_{0}$ are contracted by the Markov operator up to an additive constant. For small enough $\lambda', s>0$ for all $k\geq 0$,  $x\in X$, 
\begin{equation}
\label{eq:HC inj}
P_\mu^k u_{0}^{s}(x) \leq  e^{-s\lambda' k} u_{0}^{s}(x) + 2.
\end{equation}

We claim that  for all large enough $k$, for small enough $\lambda, s>0$,  the function $\omega^s_{0}(x,y)\defeq \dist(x,y)^{-s}$ satisfies for all $x,y\in X$,  
\begin{equation} \label{claim-omega0}
P_{\mu}^{k}\omega^{s}_{0}(x,y)\leq e^{-s \lambda k} \omega^{s}_{0}(x,y)+ 2 \,u_{0}^s(x).
\end{equation}

Assuming the claim, we conclude as follows. Fix $ k, s, \lambda', \lambda$ such that \eqref{eq:HC inj}, \eqref{claim-omega0} hold and satisfying $\lambda<\lambda'$.  We have for any $C\geq 0$, 
\begin{align*}
P_{\mu}^{k}\omega^s_C(x,y) &= P_\mu^k \omega_0^s(x,y) + C P_\mu^k u_{0}^s(x)\\
&\leq e^{-s \lambda k} \omega^s_0(x,y) + (2 + C e^{-s\lambda' k}) u_{0}^{s}(x) + 2C \\
&\leq e^{-s \lambda k} \omega^s_C(x,y) + 2C,
\end{align*}
up to taking  $C$ large enough so that $2 + C e^{-s \lambda' k} \leq C e^{-s \lambda k}$.
This shows that $P_{\mu}^{k}$ contracts $\omega^s = \omega^s_C$ up to an additive constant, and we may upgrade this to the desired property \eqref{eq-dfo} (arguing as for \Cref{drift-finite-orbits} to go from $P_{\mu}^k$ to $P_{\mu}^n$ for any $n$, and as for \Cref{drift-infinity} to allow smaller $s$ or $\lambda$).

It remains to establish \eqref{claim-omega0}.
For any $x,y\in X$, $g \in G$, distinguish two cases according to whether 
\begin{equation}
\label{eq:x close to y}
\dist(x,y) < \frac{u_{0}^{-1}(x)}{2\norm{\Ad(g)}\norm{\Ad(g^{-1})}}	
\end{equation}
holds or not.
If it holds then there is $v \in \kg$ such that $y = \exp(v)x$ and $\dist(x,y) \leq \norm{v} < 2 \dist(x,y)$.
Then, recalling $u_{0}^{-1}\leq \inj$, 
\begin{align*}
\norm{\Ad(g)v} & < 2\norm{\Ad(g)} \dist(x,y)\\
& < \norm{\Ad(g^{-1})}^{-1}\inj(x)\\
&\leq \inj(gx),
\end{align*}
which, together with $gy = \exp(\Ad(g)v)gx$, implies that $\dist(gx,gy) \gg \norm{\Ad(g)v}$.
On the other hand, if \eqref{eq:x close to y} does not hold, then
\[
\dist(gx,gy) \gg \norm{\Ad(g^{-1})}^{-1} d(x,y) \gg \norm{\Ad(g)}^{-O(1)} u_{0}^{-1}(x).
\]

Put together, we have shown that for all $x,y\in X$, there is some $v\in \kg\smallsetminus \{0\}$ such that for all $g\in G$ and all $s\in (0, 1)$,
\begin{equation}\label{eq-omega0}
\omega_0^s(gx,gy) \ll \frac{\norm{v}^{s}}{\norm{\Ad(g)v}^{s}}\omega_0^s(x,y) + \norm{\Ad(g)}^{O(s)}u_{0}^{s}(x).
\end{equation}
Using the positivity of the Lyapunov exponent and finite exponential moment (see for example by \cite[Lemma 4.2]{EskinMargulis}), we may fix (arbitrarily small) $s, \lambda>0$ depending only on $(\mu, \|\cdot\|)$ such that for all large enough $k$, we have 
\begin{equation}
\label{eq-drift>0}
\int_{G}\frac{\norm{v}^{s}}{\norm{\Ad(g)v}^{s}}\,d\mu^{*k}(g)\ll e^{-s\lambda k}.
\end{equation}
Plugging \eqref{eq-drift>0} into \eqref{eq-omega0}, using the finite exponential moment condition, and up to replacing  $\lambda$ by $\lambda/2$  and taking  $k$ larger to absorb the implicit constant in the $\ll$, we get
\[
P_{\mu}^k\omega_0^{s}(x,y) \leq e^{-s\lambda k}  \omega_0^s(x,y) + O(1)u_{0}^{s}(x).
\]
This clearly holds for smaller $\lambda$, and as for \Cref{drift-infinity},  we may use a convexity argument to allow smaller $s$   and replace $O(1)$ by $O(1)^s$. In the end, we have obtained~\eqref{claim-omega0}.
\end{proof}

We deduce the following lemma which shows how the regularity properties of a measure $\nu$ can be transferred to $\mu^{*n}*\nu$.

\begin{lemma} \label{persist-ineq}
There are $s > 0$ and $\lambda > 0$ depending only on $\mu$ and $X$ such that the following holds.
Let $\nu$ be a finite Borel measure on $X$ of total mass at most $1$. 
Then for any $n \geq 0$ and any $\rho, r > 0$, we have
\begin{equation}
\label{eq:fuir cusp}
(\mu^{*n}*\nu)\{\inj\leq r\} \ll  r^s (e^{-s\lambda n} \rho^{-1} + 1) + \nu\{\inj \leq \rho\},
\end{equation}
and 
\begin{equation}
\label{eq:fuir diag}
\sup_{x\in X}\mu^{*n}*\nu(B_{r}(x))^2 \ll  r^s (e^{-s\lambda n} \rho^{-1} + 1) + \sup_{x\in X}\nu(B_\rho(x)) + \nu\{\inj \leq \rho\}.
\end{equation}
\end{lemma}


\begin{proof}
We fix $s > 0$, $\lambda > 0$ small enough and $C > 1$ large enough so that they satisfy \Cref{drift-infinity}, \Cref{drift-function-omega} and moreover $u_{0}^{-1}\leq \inj \ll u_{0}^{-s}$ (recall \ref{comparison}, \ref{comparison2}).

Let $\nu_{|\rho}$ denote the restriction of $\nu$ to the compact subset $\{ \inj \geq \rho \} \subset X$. 
Then 
\[
(\mu^{*n}*\nu) \{\inj \leq r\} \leq (\mu^{*n}*\nu_{|\rho}) \{\inj \leq r\} + \nu \{ \inj \leq \rho \}
\]
By $u^{-1}_{0}\leq \inj$ and Markov's inequality,
\[
(\mu^{*n}*\nu_{|\rho}) \{\inj \leq r\} \leq (\mu^{*n}*\nu_{|\rho}) \{u_{0} \geq r^{-1} \} \leq r^s (\mu^{*n}*\nu_{|\rho})(u_{0}^{s}).
\]
Note that $(\mu^{*n}*\nu_{|\rho})(u_{0}^{s}) = \nu_{|\rho} (P_\mu^n u_0^s)$. Thus, by \Cref{drift-infinity}, 
\[
(\mu^{*n}*\nu_{|\rho})(u_{0}^{s}) \leq  e^{-s\lambda n}\nu_{|\rho}(u_0^{s}) + 2
\]
By the choice of $\nu_{|\rho}$,
\[
\nu_{|\rho}(u_0^{s}) \ll \nu_{|\rho}(\inj^{-1}) \leq \rho^{-1}.
\]
Putting these together we find \eqref{eq:fuir cusp}.

The proof of \eqref{eq:fuir diag} is similar but uses \Cref{drift-function-omega} instead.
Introduce the measure $\nu'_{|\rho}$ on $X\times X$ obtained by restricting $\nu_{|\rho} \otimes \nu_{|\rho}$ further to the set $\set{(x,y) \in X \times X : y \notin B_\rho(x)}$.
Note that $\nu \otimes \nu - \nu'_{|\rho}$ is a non-negative measure of total mass at most 
\[M_\rho \defeq \sup_{x\in X} \nu(B_\rho(x)) + 2 \nu\{ \inj \leq \rho \}.\]

By the Cauchy-Schwarz inequality,
\begin{align*}
\mu^{*n}*\nu(B_{r}(x))^2 
  &\leq \int_{G} g_{\star}\nu(B_{r}(x))^2 \dd\mu^{*n}(g) \\
  &=    \int_{G} g_{\star}(\nu \otimes \nu) \bigl(B_{r}(x) \times B_{r}(x) \bigr)  \dd\mu^{*n}(g) \\
  &= \mu^{*n} * (\nu \otimes \nu) \bigl(B_{r}(x) \times B_{r}(x) \bigr).
\end{align*}
Remark that $y,z \in B_r(x)$ implies $\dist(y,z) \leq 2r$ and hence $\omega_C^s(y,z) \geq (2r)^{-s}$, the last term is at most
\[
\mu^{*n}*(\nu \otimes \nu) \{\omega_{C}^s \geq 2^{-s}r^{-s} \} \leq (\mu^{*n}*\nu'_{|\rho})\{\omega_{C}^s \geq 2^{-s}r^{-s} \}+  M_\rho.
\]
We then proceed as above, using \Cref{drift-function-omega},
\begin{align*}
(\mu^{*n}*\nu'_{|\rho})\{\omega_{C}^s \geq 2^{-s}r^{-s} \}
  &\leq 2^s r^s\,  \nu'_{|\rho}(P_\mu^n \omega_C^s) \\
  &\ll r^s \bigl(e^{-s\lambda n}\,\nu'_{|\rho}(\omega_C^s) + 1\bigr)\\
  &\ll r^s \bigl(e^{-s\lambda n} \rho^{-1} + 1\bigr).
\end{align*}
Putting these together, we obtain \eqref{eq:fuir diag}.
\end{proof}

\begin{proof}[Proof of \Cref{pr:persist}]
Note we may assume $\tau=0$ and $I\subseteq (0, 1)$.
Let $s > 0$ be a constant depending only on $(X, \mu)$ as in \Cref{persist-ineq}.
Let $M>1$ be large enough such that $M s-1 > \dim X$.
We apply \eqref{eq:fuir diag} to get for all  $n\geq 0$, $\rho\in I$, $x\in X$, 
\[
\mu^{*n}*\nu (B_{\rho^M}(x))^2 \ll \rho^{M s-1}+ \rho^{\kappa \dim X} \leq \rho^{\frac{\kappa}{2} \dim X}
\]
where the last inequality assumes $\sup I$ small enough in terms of $(X,\mu, \kappa)$.
Moreover, \eqref{eq:fuir cusp} gives
\[
\mu^{*n}*\nu \{\inj \leq \rho^M\}  \ll \rho^{M s-1} \leq \sup I.
\]
Hence $\mu^{*n}*\nu$ is $(\frac{\kappa}{4M}, I', \sup I)$-robust with $I'\defeq \{\rho^M\,:\,\rho\in I\}$.
\end{proof}

\subsection{From small dimension to high dimension}  \label{ss:bootstrap}

The goal of this subsection is to show dimension increment results.
Roughly speaking if a distribution on $X$ has some dimension away from both $0$ and $\dim X$, then the random walk will increase this dimension. 
The precise statement is \Cref{global-single-scale-bootstrap} below.
Then it can be iterated to get to a dimension arbitrarily close to be full, yielding the statement of \Cref{dimension-bootstrap}.

\begin{proposition}[Dimension increment] \label{global-single-scale-bootstrap}
Consider the setting \ref{notation-bootstrap} and assume $\kg=\so(2,1)$ or $\kg=\so(3,1)$.
Let $\kappa, \eps, \delta > 0$, $\tau \geq 0$ and $\alpha \in {[\kappa, 1 - \kappa]}$ be some parameters. 
Let $\nu$ be a Borel measure on $X$ that is $(\alpha, \cB_{[\delta, \delta^{\eps}]}, \tau)$-robust. Set $n_{\delta} \geq 0$ to be the integer part of $\frac{1}{2\lambda_{\mu}}\abs{\log \delta}$.

If $\eps,  \delta \lll_{ \kappa} 1$,  then $\mu^{*n_{\delta}}*\nu$ is $(\alpha+\eps, \cB_{\delta^{1/2}}, \tau+\delta^{\eps})$-robust. 
\end{proposition}

The next three paragraphs are devoted to the proof of this proposition.
Then in paragraph \ref{sss:dimension-bootstrap}, we show how to iterate it to get the dimension bootstrap.

Now, we give a heuristic of the argument.
The aim is to bound for every $x \in X$, after throwing away exponentially small mass, the measure $\mu^{*n}*\nu(B_\rho(x))$ of the ball of radius $\rho = \delta^{1/2}$ centered at $x$.
It is the average of $\nu(g^{-1}B_\rho(x))$ with $g$ distributed according to $\mu^{*n}$.
If we understand the typical shape of $g^{-1}B_\rho(x)$ then we can hope to apply the multislicing theorem\footnote{Needless to say, it is by analyzing this problem of dimension increment that we realized that \Cref{coro:slicing} is the relevant slicing theorem.}.

In our setting, the random walk on $G$ induced by $\mu$ has three Lyapunov exponents: $\lambda_\mu$, $0$ and $-\lambda_\mu$.
Hence a typical $g$ will dilate in its unstable direction by a factor of $e^{n \lambda_\mu}$ and shrink in its stable direction by a factor of $e^{- n \lambda_\mu}$.
Thus, $g^{-1}B_\rho(x)$ viewed in an appropriate local chart is roughly a rectangle of length $e^{n \lambda_\mu}\rho$ in the $g$-stable direction and $\rho$ in the $g$-central direction and $e^{-n \lambda_\mu}$ in the $g$-unstable direction.
With our choice of $n = n_\delta$, these lengthes are roughly $1$, $\delta^{1/2}$ and $\delta$.
The flag associated to this rectangle depends on the first $K$-element in the $KAK$-decomposition of $g^{-1}$.
Denote this element by $\theta =\theta_g$.

We will define in paragraph~\ref{sss:charts}, for this $\theta \in K$ associated to $g$, a map $\psi_\theta \colon \kg \to G$ whose local inverse $\phi_\theta \colon B^G_{\rho_0} \to \kg$ is a chart in which every ball $g^{-1}B_\rho(x)$ is roughly a rectangle of side lengths $1$, $\delta^{1/2}$ and $\delta$, i.e. $D_\delta^{(0,1/2,1)}$ in the notation of \Cref{sec:slicing}.
Actually, $\phi_\theta$ is a chart that trivializes simultaneously the $g$-stable and $g$-stable-central foliations.
Next, in paragraph~\ref{sss:noncen-foliation}, we will show that the distribution of these foliations, or equivalently, the distribution of the corresponding flags satisfies the non-concentration condition~\ref{it:NCsigma} of \Cref{thm:slicing}.
Thus \Cref{coro:slicing} can be applied to get an upper bound for $\nu(g^{-1}B_\rho(x))$ for a typical $g \in G$ and $x \in X$.

\subsubsection{Charts straightening images of balls under group action}\label{sss:charts}
In this paragraph, notations refer to \ref{notation-bootstrap}  with the additional assumption that the Lie algebra of $G$ is  $\mathfrak{so}(n,1)$ ($n\geq 2$); i.e. $G$ is a simple connected real Lie group whose restricted root system is of type $A_1$.

Recall that $\ka$ denotes a Cartan subspace (chosen to be in the orthogonal of $\Lie(K)$ for the Killing form) and set $A = \exp(\ka)$.
Thus the subgroup $A$ is isomorphic to $\R$ and there is a unique parametrization $t \in \R \mapsto a^t \in A$ so that $\exp(\ka^+) = \set{a^t: t\geq 0}$ and the restricted root space decomposition of $\kg$ is
\[
\kg = \kg_+ \oplus \kg_0 \oplus \kg_-,
\]
where for $\bullet \in \{+,-\}$
\[
\kg_\bullet = \set{ v \in \kg : \forall t \in \R,\, \Ad(a^t) v = e^{\bullet t} v}
\]
and
\[
\kg_0 = \set{ v \in \kg : \forall t \in \R,\, \Ad(a^t) v =  v}
\]
is the centralizer of $\ka$. 

As a substitute to the exponential map, define the map 
\[\psi \colon \kg \to G\]
such that for all $(v_+, v_0, v_-)  \in\kg_+ \times \kg_0 \times  \kg_-$,
\[
\psi(v_+ + v_0 + v_-) = \exp(v_+)\exp(v_0)\exp(v_-).
\]

\begin{lemma}
\label{lm:param-loc} 
There is a constant $\rho_0=\rho_{0}(G, \|\cdot\|, \ka^+)  > 0$ such that 
for any $t \geq 0$, any $\rho > 0$ with $e^t\rho \leq \rho_0$ and any $h \in G$, there is $w \in \kg$ such that
\begin{equation}\label{eq:param-loc}
\set{v \in B^\kg_{\rho_0} : \psi(v) \in a^t B^G_\rho h} \subset \Ad(a^t) B^\kg_{10^6\rho} + w.
\end{equation}
\end{lemma}
One could also prove a converse inclusion, with $ B^\kg_{10^6\rho}$ replaced by $ B^\kg_{10^{-6}\rho}$ but we do not dwell on this as it is not useful for our bootstrap argument. 

Note that $\Ad(a^t) B^\kg_{10^6\rho}$ is roughly the rectangle $B^{\kg_+}_{10^6e^t\rho} \times B^{\kg_0}_{10^6\rho} \times B^{\kg_-}_{10^6e^{-t}\rho}$.
So the local chart $\psi$ straightens all the images of balls under $a^{t}$ \emph{simultaneously} into rectangles. 
It is not difficult to see that the exponential map $\exp \colon \kg \to G$ does not fulfill this job.

\begin{proof}
In this proof, each appearance of the notation $\bar\rho$ stands for a quantity of the form $2^C\rho$ where $C$ is a number that one can make explicit by following carefully the proof.

Let $v, w$ be arbitrary elements of the set on the left-hand side of \eqref{eq:param-loc}.
The goal is to show that $v - w \in \Ad(a^t) B^\kg_{\bar\rho}$.

Note that both $a^{-t}\psi(v)$, $a^{-t}\psi(w) \in B^G_\rho h$.
Hence by the triangle inequality,
\[
\psi(v) \in a^t B^G_{2\rho}a^{-t} \psi(w).
\]
We choose $\rho_0> 0$ small so that the local diffeomorphism $\psi \colon \kg \to G$ is injective and $2$-bi-Lipschitz on $B_{2\rho_0}$.
It follows that there is $u = u_+ + u_0 + u_- \in \kg$ with 
\[
u_+ \in B^{\kg_+}_{e^t\bar\rho},\; u_0 \in B^{\kg_0}_{\bar\rho}\;\text{and}\; u_- \in B^{\kg_-}_{e^{-t}\bar\rho} 
\]
and such that
\[
\psi(v) = \psi(u)\psi(w) = \exp(u_+)\exp(u_0)\exp(u_-) \exp(w_+)\exp(w_0)\exp(w_-),
\]
where, obviously, $w = w_+ + w_0 + w_-$ is the decomposition of $w$ in the weight spaces.

Then, we permute the the order of $\exp(u_-)$ with $\exp(w_+)\exp(w_0)$.
Note that the map 
\[
(w_+,w_0,u_-) \in \kg_+ \times \kg_0 \times \kg_- \mapsto (w_+',w_0',u_-') \in \kg_+ \times \kg_0 \times \kg_-, 
\]
defined by the relation  
\begin{equation}
\label{eq:permute-+0}
\exp(w_+')\exp(w_0')\exp(u_-') = \exp(u_-) \exp(w_+)\exp(w_0)
\end{equation}
is a local diffeomorphism at $0$ with differential at $0$ being the identity map. 
Hence it is $2$-Lipschitz on $B_{2\rho_0}$ if we choose $\rho_0 \lll 1$.
For our vectors, note that $(w_+,w_0,0)$ is fixed by this map. Hence \eqref{eq:permute-+0} holds
for some 
\[
w_+' \in B^{\kg_+}_{e^{-t}\bar\rho}(w_+),\; w_0' \in B^{\kg_0}_{e^{-t}\bar\rho}(w_0) \text{ and } u_-' \in B^{\kg_-}_{e^{-t}\bar\rho}.
\]

Next, we permute the order of $\exp(u_0)$ with $\exp(w_+')$.
Note that $\kg_+ \oplus \kg_0$ is a Lie-subalgebra, hence we may perform the same operation as above.
The map 
\[
(w'_+,u_0) \in \kg_+ \times \kg_0 \mapsto (w_+'',u_0') \in \kg_+ \times \kg_0, 
\]
defined by the relation  
\begin{equation}
\label{eq:permute0+}
\exp(w_+'') \exp(u_0') = \exp(u_0)\exp(w_+')
\end{equation}
is $2$-Lipschitz on $B_{2\rho_0}$ and fixes the point $(w_+',0)$.
We conclude that for our vectors, \eqref{eq:permute0+} holds for some
\[
w_+'' \in B^{\kg_+}_{\bar\rho}(w_+') \subset B^{\kg_+}_{\bar\rho}(w_+) \text{ and } u_0' \in B^{\kg_0}_{\bar\rho}
\]

After the two permutations, we find 
\[
\psi(v) = \exp(u_+) \exp(w_+'') \exp(u_0') \exp(w_0') \exp(u_-') \exp(w_-).
\]
Then, using that $\exp$ is $2$-bi-Lipschitz near $0$ on each $\kg_\bullet$, $\bullet \in \{+,0,-\}$,
we have (for $\rho_0 \lll 1$)
\[
\psi(v) \in \exp(B^{\kg_+}_{e^t\bar\rho} + w_+) \exp(B^{\kg_0}_{\bar\rho} + w_0) \exp(B^{\kg_-}_{e^{-t}\bar\rho} + w_-).
\]
Using that $\psi$ is injective on $B^\kg_{2 \rho_0}$ and $e^t\rho \leq \rho_0$, we conclude that $v - w \in \Ad(a^t) B^\kg_{\bar\rho}$.

\end{proof}

\begin{remark}
The reason we restrict to simple real Lie groups of type (restricted root system) $A_1$ is that they enjoy the following property:
for any affine half space $D$ of the dual $\ka^*$, the sum of rootspaces parametrized by roots in $D$ is a Lie subalgebra.
Simple real Lie groups of other types, e.g. $\sl_{3}(\R)$, do not have this property.
\end{remark}


\bigskip

For $\theta \in K$, define $\psi_\theta \colon \kg \to G$ by 
\[
\forall v \in \kg,\quad \psi_\theta(v) = \theta \psi(v) \theta^{-1}.
\]
Then \Cref{lm:param-loc} gives
for any $t \geq 0$, $h \in G$, the existence of some $w \in \kg$ such that
\begin{equation}\label{eq:param-loc-theta}
\set{v \in B^\kg_{\rho_0} : \psi_\theta(v) \in \theta a^t B^G_\rho h} \subset \Ad(a^t) B^\kg_{10^6\rho} + w.
\end{equation}
Hence for any element $g\in G$ of Cartan decomposition $g=\theta a^{t} \theta'$, the $g$-translate of a small ball looks like a Euclidean rectangle in the chart  $\psi_{\theta}$.

\subsubsection{Non concentration estimates for unstable and central unstable leaves}\label{sss:noncen-foliation}

Given $g \in G$, we consider the Cartan decomposition of $g^{-1}$ : 
\begin{equation*}
g^{-1} = \theta_g a^{t_g} \theta'_g  \,\,\,\,\,\,\,\,\,\,\text{ where $\,\,\,\,\,\,\,\,\,\,\theta_{g}, \theta'_{g}\in K$, $t_{g}\geq 0$.}
\end{equation*}
We plan to apply the multislicing theorem (more precisely \Cref{coro:slicing}) to rectangles of the form $\psi_{\theta_{g}}(\Ad(a^{t_{g}}) B^\kg_{\rho} + w)$ where $g$ is random element in $G$ chosen with law $\mu^{*n}$. Here we check the relevant non-concentration properties for these rectangles.  We let $\rho_{0}>0$ be small enough in terms of $G, \|\cdot\|, \ka^+$ so that every $\psi_{\theta}$ restricts as a diffeomorphism between a neighborhood of $0$ in $\kg$ and the ball $B_{\rho_{0}}^G$, we write $\varphi_{\theta}: B_{\rho_{0}}^G\rightarrow \kg$ the inverse map. 


\begin{lemma}\label{lm:noncon-foliation}
Consider the setting \ref{notation-bootstrap} and assume moreover  $\kg = \so(2,1)$ or $\kg=\so(3,1)$.
There exist constants $\kappa > 0$ and $C = C> 1$ such that the following holds for every integer $n \geq C$.

Writing $\sigma$ for the image measure of $\mu^{*n}$ by the map $g \mapsto \theta_g$, we have for any $h \in B^G_{\rho_{0}}$, any $\rho \geq e^{-n}$, 
\[
\forall W \in \Gr(T_h G,\dim \kg_0 \oplus \kg_-), \quad \sigma\set{ \theta \in K : \dang( (D_h \varphi_{\theta})^{-1}\kg_{+} ,W) \leq \rho} \leq C\rho^\kappa,
\]
and 
\[
\forall W \in \Gr(T_h G,\dim \kg_-),\quad \sigma\set{ \theta \in K : \dang((D_h \varphi_{\theta})^{-1}(\kg_{+}\oplus \kg_{0}),W) \leq \rho} \leq C\rho^\kappa.
\]
\end{lemma}

\begin{proof}
We see directly from the definition of $\psi$ that $D_{v}\psi(\kg_+)$ is given by the set of derivatives  $\frac{d}{dt}_{|t=0}\exp(t Z)\psi(v)$ where $Z$ runs through $\kg_{+}$. It follows that $\set{(D_h \varphi_{\theta})^{-1}\kg_{+} : h\in B^G_{\rho_{0}}}$ coincides with the right $G$-invariant sub-bundle of $TG$ determined by $\Ad(\theta)\kg_{+}$ at $\Id$. 
This also holds  with $ \kg_{+}\oplus \kg_{0}$ in place of $\kg_{+}$, by the same argument. 


Thus, in order to prove the lemma, it suffices to show the claim for $h = {\Id}$, i.e. we need to check that for $n \ggg 1$ and for all $\rho \geq e^{-n}$,
\begin{equation}
\label{eq:noncon-Adtheta2}
\forall W \in \Gr(\kg,\dim  \kg_0 \oplus \kg_-),\quad \mu^{*n}\setbig{ g \in G : \dang(\Ad(\theta_g)\kg_+, W) \leq \rho} \ll \rho^\kappa,
\end{equation}
and
\begin{equation}
\label{eq:noncon-Adtheta1}
\forall W \in \Gr(\kg,\dim \kg_-),\quad \mu^{*n}\setbig{ g \in G : \dang\bigl(\Ad(\theta_g)(\kg_+ \oplus \kg_0),W\bigr) \leq \rho} \ll \rho^\kappa.
\end{equation}
Both these estimates are known and are manifestations of the Hölder regularity of Furstenberg measures first shown by Guivarc'h~\cite{Guivarch}.

In the case  $\kg = \so(2,1)$, the linear random walk on $\kg$ defined by $\Ad_\star \check\mu$ (where $\check\mu$ is the image of $\mu$ by $g\mapsto g^{-1}$) is strongly irreducible and proximal.
Thus \eqref{eq:noncon-Adtheta2} is a special case of \cite[Lemma 4.5]{BFLM}.
For \eqref{eq:noncon-Adtheta1}, we use the relation, see for example \cite[Equation (14)]{He2020JFG},
\[
\dang\bigl(\Ad(\theta_g)(\kg_+ \oplus \kg_0),W\bigr) = \dang\bigl(\Ad(\theta_g)(\kg_+ \oplus \kg_0)^\perp,W^\perp\bigr),
\]
and, by the assumption that $\Ad(a^t)$ is self-adjoint, we have 
\[\Ad(\theta_g)(\kg_+ \oplus \kg_0)^\perp = \Ad(\theta_g) \kg_-\]
is the image of the expanded subspace in the singular decomposition of the adjoint $\Ad(g)^*$ of $\Ad(g)$.
Hence \eqref{eq:noncon-Adtheta1} follows from \cite[Lemma 4.5]{BFLM} applied to the random walk defined by the pushforward of $\mu$ by the map $g \mapsto \Ad(g)^*$.

In the case  $\kg = \so(3,1) = \sl_2(\C)$, the linear random walk on $\kg$ defined by $\Ad_\star\check\mu$ is strongly irreducible but not proximal.
However, it satisfies the property $(S)$ of \cite[Definition 1.1]{He2020IJM} by \cite[Proposition 2.5]{He2020IJM}, so that \cite[Proposition 3.1(iii)]{He2020IJM} applies.
Applied with $\omega = 1$ and $l = \log \rho$ it gives immediately~\eqref{eq:noncon-Adtheta2}.
In a dual manner like above, it also gives \eqref{eq:noncon-Adtheta1}.
\end{proof}

One may ask if this argument works for $G = \SO(n,1)$ with $n \geq 4$, which also has type $A_1$ restricted root system.
Note that for $n \geq 4$, $G$ does not have the property $(S)$ of \cite{He2020IJM}.
Hence \Cref{lm:noncon-foliation} fails.
So \Cref{coro:slicing} in its current form does not apply to these groups.

\subsubsection{Proof of dimension increment}
Note that since the restricted root system on $\kg$ is of type $A_1$, we know that 
$\lambda_{\mu}$ is also the top Lyapunov exponent of the random walk on $\kg$ induced by $\Ad_\star\check\mu$ where $\check\mu$ denotes the pushforward measure of $\mu$ by the map $g \mapsto g^{-1}$.

\begin{proof}[Proof of \Cref{global-single-scale-bootstrap}]
Clearly, this reduces to the regular part of $\nu$. Hence, we may assume that $\tau = 0$.

Then by the robustness assumption, $\nu$ is supported on $\{\inj \geq \d^\eps\}$.
We can cover the latter by balls $B_i = B^X_{\d^{2\eps}}(x_i)$ of radius $\d^{2\eps}$:
\[
X_{\d^\eps} \subset \bigcup_{i\in I} B_i
\]
The number of balls we need can be bounded: $\sharp I \leq \d^{-O(\eps)}$.

Restrict $\nu$ to each $B_i$ and pull back to $B^G_{\d^{2\eps}}$ by the isometries $g \in B^G_{\d^{2\eps}} \mapsto g x_i \in B_i$.
Thus, we can write 
\[
\nu = \sum_{i\in I} \nu_i * x_i
\]
as the sum of measures with $\supp(\nu_i) \subset B^G_{\d^{2\eps}}$. Here to avoid confusion between different meanings for $\delta$ we just write $x_{i}$ for the Dirac mass at $x_{i}$ (written  $\delta_{x_{i}}$ elsewhere in the text). 
Note that each $\nu_i$ is $(\alpha,\cB_{[\delta,\delta^\eps]},0)$-robust.

For $g \in G$, consider its Cartan decomposition
\[
g^{-1} = \theta_g a^{t_g} \theta'_g
\]
with $\theta_g$, $\theta'_g \in K$ and $t_g \geq 0$.

Consider the family of diffeomorphisms $\phi_\theta \colon B^G_{\rho_0} \to \kg$, indexed by $\theta \in \Theta \defeq K$, and introduced in \S\ref{sss:noncen-foliation}.
Let $\sigma$ be the pushforward measure of $\mu^{*n}$ by the map $g \mapsto \theta_g$.
In view of \Cref{lm:noncon-foliation}, we can apply \Cref{coro:slicing} in this setting, with the flag $\{0\} \subset \kg_+ \subset \kg_+\oplus \kg_0 \subset \kg$, the parameters $(r_1, r_2) = (0,\frac{1}{2})$ and at scale $\delta$.
This gives some $\eps_1>0$ depending only on Setting \ref{notation-bootstrap}  such that,  provided  $\eps, \delta  \leq \eps_1$ and $n\geq n_{\delta}$, for each $i \in I$,
there is a set $\sD_i \subset G$ with $\mu^{*n}(\sD_i) \geq 1 - \delta^{\eps_1}$ and that 
for every $g \in \sD_i$, $\nu_i$ contains a component measure $\nu_{i,g}$ of total mass $\nu_{i,g}(G) \geq \nu_i(G)-\d^{\eps_1}$ such that 
\begin{equation}\label{eq:psithetaw+R}
\forall w \in \kg,\quad \nu_{i,g}(\varphi^{-1}_{\theta_g}(w + R_\d)) \leq \d^{\frac{\a}{2}\dim G+ \eps_1}
\end{equation}
where $R_\d = B^{\kg_+}_1 + B^{\kg_0}_{\d^{1/2}} + B^{\kg_-}_{\d}$ is the rectangle associated to the data of the flag and $(r_1,r_2)$ at scale $\d$.

Observe that $t_g = \log \norm{\Ad(g^{-1})}$. Given a parameter $\eps_{2}>0$, 
the large deviation estimate (we can use \cite[Theorem V.6.2]{BougerolLacroix}, see also \cite[Theorem 13.17(iii)]{BQ_book}) for the random walk on $\kg$ defined by $\Ad_\star\check\mu$ asserts that there is $c = c(\mu, \eps_2) > 0$ such that for all $n \ggg_{\eps_{2}} 1$, we have $\mu^{*n}(\sD') \geq 1 - e^{-c n}$, where
\[\sD' \defeq  \setbig{g \in G : \abs{t_g - n \lambda_{\mu}} \leq \eps_2 n \lambda_{\mu}  }.\]
Assuming $\delta\lll_{\eps_{2}}1$, we may specialize to $n = n_\delta = \lfloor\frac{1}{2\lambda_{\mu}} \abs{\log \delta} \rfloor$ so that $e^{-cn} \leq \delta^{\eps_3}$ for some $\eps_3 = \eps_3(\mu,\eps_2) > 0$ and
\begin{equation}
\label{eq:LD tg}
\forall g \in \sD', \quad \delta^{-\frac{1}{2} + 2 \eps_2} \leq e^{t_g} \leq \delta^{-\frac{1}{2} - 2 \eps_2}.
\end{equation}

Let $\rho = \delta^{\frac{1}{2} +  2 \eps_2 + 2 \eps}$.
We claim that if $\delta\lll_{\eps}1$, then for $g \in \sD_i \cap \sD'$, 
\begin{equation}\label{eq:gnuixi}
\forall x \in X,\quad (g * \nu_{i,g} * x_i) (B^X_\rho(x)) \ll \d^{\frac{\a}{2}\dim G +\eps_{1}- O(\eps_2) - O(\eps)}.
\end{equation}
This is the $\nu_{i,g}$-measure of the set of $h \in B^G_{\delta^{2\eps}}$ such that
\[
ghx_i \in B^G_\rho x.
\]
If $h$ and $h_0$ are both in this set,
then by the triangle inequality,
\[
gh x_i \in B^G_{2\rho}\, gh_0 x_i,
\]
and then 
\[
h x_i \in g^{-1}B^G_{2\rho}\, gh_0 x_i.
\]
Note that $h \in B^G_{\d^{2\eps}}$ and $g^{-1}B^G_{2\rho}\, gh_0\subset B^G_{4 e^{t_g} \rho} B^G_{\d^{2\eps}} = B^G_{4 e^{t_g} \rho + \d^{2\eps}}$.
By the choice of $\rho$ and \eqref{eq:LD tg}, we have $4 e^{t_g} \rho + \d^{2\eps} \leq \d^\eps$.
Remembering that $\inj(x_i) \geq \d^{\eps}$, we obtain
\[
h \in g^{-1}B^G_{2\rho}\,gh_0 = \theta_g a^{t_g} B^G_{2\rho}\, a^{-t_g} \theta_g^{-1} h_0.
\]
In view of \eqref{eq:param-loc-theta} after \Cref{lm:param-loc} and provided $\delta \lll_{\eps}1$, such $h$ is contained in a set of the form
\[
\varphi^{-1}_{\theta_g}(w + \Ad(a^{t_g}) B^\kg_{2 \cdot 10^6 \rho})
\]
Note that $\Ad(a^{t_g}) B^\kg_{2 \cdot 10^6 \rho}$ is contained in a constant dilation of the rectangle of side length $e^{t_g}\rho$, $\rho$ and $e^{-t_g \rho}$ in respectively the $\kg_+$, $\kg_0$, $\kg_-$ direction.
In view of \eqref{eq:LD tg}, this is covered by at most $\d^{-O(\eps_2) - O(\eps)}$ translates of $R_\d$.
Now the claim \eqref{eq:gnuixi}  follows from \eqref{eq:psithetaw+R}.

To conclude, consider inside $\mu^{*n} * \nu$ the component
\[
\nu' = \sum_{i \in I} \int_{\sD_i \cap \sD'} g * \nu_{i,g} * x_i \dd \mu^{*n}(g).
\]
On the one hand, for every $x \in X$, $B^X_{\delta^{1/2}}(x)$ can be covered by at most $\bigl(\frac{\d^{1/2}}{\rho}\bigr)^{O(1)} \leq \d^{- O(\eps_2) - O(\eps)}$
balls of radius $\rho$, hence \eqref{eq:gnuixi} implies
\[
\nu'\bigl(B^X_{\delta^{1/2}}(x)\bigr) \leq \abs{I} \d^{- O(\eps_2) - O(\eps)} \d^{\frac{\a}{2}\dim G +\eps_1 - O(\eps_2) - O(\eps)}  \leq \d^{\frac{\a}{2}\dim G + \eps_1- O(\eps_2) - O(\eps)}.
\]
Hence 
\[
\nu'\bigl(B^X_{\delta^{1/2}}(x)\bigr) \leq \d^{\frac{\alpha+\eps}{2}\dim X},
\]
provided $\eps \leq \eps_2 \lll \eps_1$. We may suppose $\eps_{2}$ is fixed from the start, depending on the initial setting \ref{notation-bootstrap} and $\eps_{1}$, so that this holds.

On the other hand, the mass missing in $\nu'$ as compared to $\mu^{*n}*\nu$ is at most
\begin{align*}
& \sum_{i \in I} \mu^{*n}(G \setminus \sD_i \cup G \setminus \sD') \nu_i(G) + \sum_{i \in I} \int_{G_2} (\nu_i(G) - \nu_{i,g}(G)) \dd \mu^{*n}(g)\\
 \leq & \,\d^{-O(\eps)}(\d^{\eps_1} + \d^{\eps_3} + \d^{\eps_1})  \leq \d^{2\eps},
\end{align*}
provided  $\eps \lll_{\kappa} 1$  and $\delta \lll_{\eps}1$.

We also need to control the probability of falling into a cusp. By \Cref{effective-recurrence},  there exist $A, s>0$ depending on $(X,\mu)$ such that for every $x\in X$, $k\geq A |\log \inj(x)|$,  $r>0$, 
 \[
\mu^{*k}*x \,\{\inj \leq r\}\leq A \,r^{s}.
\]
Up to assuming $\eps$ small enough so that $\frac{1}{2\lambda_{\mu}}>A\eps$ and $\delta\lll_{\eps}1$, we may apply this to each point point $x$ in the support of $\nu$, to $k=n(=n_{\delta})$, and $r=\delta^{1/2}$ to get 
 \[
\mu^{*n}*\nu \{\inj \leq \delta^{1/2}\}\leq A\delta^{s/2} \nu(X)\leq A\delta^{s/2-O(\eps)}\leq \delta^{2\eps}
\]
where the last inequality holds up to assuming $\eps\lll1$ and $\delta$ small enough. 

We conclude that the measure $\mu^{*n}*\nu$ is $(\alpha+\eps, \cB_{\delta^{1/2}}, \delta^{\eps})$-robust.

\end{proof}

\subsubsection{Iterate the dimension increment}\label{sss:dimension-bootstrap}

We conclude the subsection with \Cref{dimension-bootstrap}, stating that the convolution with $\mu^{*n}$ bootstraps the $\kappa$-robustness of a measure on $X$ to a $(1-\kappa)$-robustness for any given $\kappa > 0$. 

\begin{corollary}[From small dimension to high dimension] \label{dimension-bootstrap}
Consider the setting \ref{notation-bootstrap} and assume $\kg=\so(2,1)$ or $\kg=\so(3,1)$.
Let $\kappa,c,\eta, \delta> 0$, $\tau \geq 0$ be some parameters and let $\nu$ be a Borel measure on $X$ such that for all $n \geq 0$, we have
\[\text{$\mu^{*n}*\nu\,$ is  $\,(\kappa, \cB_{[\delta, \delta^\eta]}, \tau)$-robust}.\]
If  $c \lll_{\kappa} 1$ and $\eta, \delta \lll_{\kappa, c}1$ then 
 for all $n \geq  \abs{\log \delta}$, we have
\[\text{$\mu^{*n}*\nu\,$ is $\,(1 - \kappa,  \cB_{\delta^c},  \tau')$-robust},\]
where $\tau' \ll_{\kappa,c} \tau + \delta^\eta$.
\end{corollary}



The proof is rather straightforward.
Note however that \Cref{global-single-scale-bootstrap} yields better dimensional properties only at a single scale, while it requires robustness on a wide range of scales as input.
So in order to iterate the dimension increment, we first upgrade the output of \Cref{global-single-scale-bootstrap} to a robustness on a (arbitrarily wide) range of scales. 

\begin{lemma}[Dimension increment at multiple scales] \label{global-multiscale-bootstrap}
Let $\kappa > 0$,  $\alpha \in {[\kappa,1-\kappa]}$, $\eps, s \in (0,1/10)$, $m \geq 1$ and $\tau \geq 0$ be some parameters. 
The following holds for all $\eps\lll_{\kappa}1$ and $\delta  \lll_{\kappa, s} 1$.

If $\nu$ is a Borel measure on $X$ such that for all $n \geq m$, 
\[\mu^{*n}*\nu \,\text{ is }\, (\alpha, \cB_{[\delta, \delta^{s\eps}]}, \tau)\text{-robust},\]
then  for all $n \geq m + \frac{1}{2\lambda_{\mu}} \abs{\log \delta}$,
\[\text{$\mu^{*n}*\nu\,$ is $\,\bigr(\alpha+\eps/2,  \cB_{[\delta^{1/2}, \delta^{s/2}]}, 10\eps^{-1}\abs{\log s}(\tau+\delta^{s\eps})\bigr)$-robust.}\]
\end{lemma}

\begin{proof}
We assume $\eps$ small enough so that \Cref{global-single-scale-bootstrap} holds.
Let $\rho \in [\delta,\delta^s]$.
Recall $n_\rho \defeq \bigl\lfloor\frac{\abs{\log \rho}}{2\lambda_\mu}\bigr\rfloor$ as in the statement of \Cref{global-single-scale-bootstrap} and note that $n \geq n_\rho$.
Taking $\delta \lll_{\kappa, s} 1$, we may assume $\rho$ arbitrarily in terms of $X, \mu, \kappa$.
Thus, \Cref{global-single-scale-bootstrap} can be applied to the measure $\mu^{*(n-n_\rho)}*\nu$ at the scale $\rho$. This yields that 
$\mu^{*n}*\nu$ is  $(\alpha+\eps,  \cB_{\rho^{1/2}}, \tau+\rho^{\eps})$-robust.

Using \Cref{single-vs-multiple}, 
we combine those single-scale estimates to obtain that $\mu^{*n}*\nu$ is  $((\alpha+\eps) - \eps/2,  \cB_{[\delta^{1/2}, \delta^{s/2}]}, \bigl\lceil\frac{\abs{\log s}}{\log (1-\eps/2)}\bigr\rceil (\tau+\delta^{s\eps}))$-robust.
This proves the lemma.
\end{proof}

\begin{proof}[Proof of  \Cref{dimension-bootstrap}]
We prove the statement under the assumption  $n\geq \frac{1}{\lambda_{\mu}} \abs{\log \delta}$ instead of $n\geq \abs{\log \delta}$. Note that the new statement is equivalent (if $\lambda_\mu<1$, replace $\delta$ by $\delta^{\lambda_\mu}$). 

Fix $\eps=\eps(X, \mu, \kappa)>0$  as in  \Cref{global-multiscale-bootstrap}.
Let $k \geq 0$ be the smallest integer such that $\kappa + k\eps/2 > 1-\kappa$.
Up to assuming $c \leq 2^{-k}$, $\eta \leq \eps^k c$ and $\delta \lll_{\kappa, c} 1$, 
we may apply $k$ times \Cref{global-multiscale-bootstrap}.

Indeed, at step $j \in \{0,\dotsc, k-1\}$, we apply \Cref{global-multiscale-bootstrap} at scale $\delta = \delta_j \defeq \delta^{2^{-j}}$ with the parameters $\alpha = \alpha_j \defeq \kappa + j\eps/2$ and $s = s_{j+1} \defeq c \,\eps^{k-(j+1)}$ and to the integer $m = m_j$ defined recursively by $m_0 = 0$ and $m_j = m_{j-1} + \bigl\lceil\frac{\abs{\log \delta_{j-1}}}{2\lambda_\mu}\bigr\rceil$.
This gives  for all $n \geq m_{j+1}$ that $\mu^{*n}*\nu$ is $(\alpha_{j+1}, \cB_{[\delta_{j+1},\delta_{j+1}^{s_{j+1}}]},\tau_{j+1})$-robust, where
$\tau_j$ defined recursively by $\tau_0 \defeq \tau$ and $\tau_{j+1} \defeq 10\eps^{-1}\abs{\log s_{j+1}}(\tau_j + \delta_j^{s_j})$.
This finishes the proof since $m_k \leq \frac{\abs{\log\delta}}{\lambda_\mu}$ for $\delta \lll_{k} 1$, the interval $[\delta_k,\delta_k^{s_k}] = [\delta^{2^{-k}},\delta^c]$ contains $\delta^c$ and a simple induction shows that for each $j = 1,\dotsc,k$,
\(
\tau_{j+1} \ll_{\eps,c,j} \tau + \delta^\eta.
\)
In particular, $\tau' \defeq \tau_k \ll_{\kappa,c} \tau + \delta^\eta$.
\end{proof}

\subsection{From high dimension to equidistribution} \label{ss:high-to-equid}

To conclude the section, we explain why a measure $\nu$ with high dimension at scale $\delta$ reaches $\delta$-equidistribution exponentially fast under convolution by $\mu$. The argument can be formulated in the following more general framework.  

Let $X$ be a locally compact separable metric space. Let $P$ be a \emph{Markov-Feller operator} on $X$. Here we mean that $P$ is an  operator on $C^{0, 0}(X)$  the space of bounded continuous  functions on $X$, and that $P$ is non-negative and satisfies $P {\mathbf 1}_{X}={\mathbf 1}_{X}$. In particular, $P$ has an operator norm equal to $1$ when  $C^{0, 0}(X)$ is endowed with the supremum norm. 

Let $m_{X}$ be a \emph{$P$-invariant} Borel probability measure on $X$.  The $P$-invariance of $m_{X}$ means that  $m_{X}P=m_{X}$ where $m_{X}P$ is the probability measure characterized by the relation $m_{X}P(\varphi)=m_{X}(P\varphi)$ for all $\varphi\in C^{0, 0}(X)$. By Jensen's inequality, $P$ extends to a norm $1$ operator of $L^2(X, m_{X})$ (see for instance \cite[Lemma 2.1]{BQ_book}), and of course $P$ preserves  the closed subspace of zero-mean  functions $L^2_{0}(X, m_{X})$.

Given another  Borel probability measure $\nu$ on $X$ whose support is included in that of $m_{X}$, and given $\delta>0$, we define the \emph{mollification of $\nu$ at  scale $\delta$} as the probability measure such that
for any measurable $f \colon X\rightarrow {[0, +\infty]}$,
\begin{align} \label{defnudelta}
\int_X f \dd \nu_\delta = \int_X  \frac{1}{m_{X}(B_\delta(y))} \int_{B_\delta(y)} f \dd m_{X}\dd \nu (y).
\end{align}
Note that $\nu_{\delta}$ is absolutely continuous with respect to $m_{X}$. By abuse of notation, we still write $\nu_\delta$ to denote its Radon-Nikodym derivative:
\[
\nu_\delta(x) =  \int_{B_\delta(x)}\frac{1}{m_{X}(B_\delta(y))}\dd \nu(y).
\]
\begin{proposition}[Endgame] \label{endgame}
Assume $P$ has spectral radius strictly less than $1$ on $(L^2_{0}(X, m_{X}), \|\cdot\|_{L^2})$ and that there exists $\b \in (0, 1]$ such that $P$ restricts as a bounded operator on $(C^{0, \b}(X), \|\,\cdot\,\|_{C^{0, \beta}})$.  
Given $r_{0}\in (0, 1)$ there are constants $\kappa , r, \eps > 0$ with $r<r_{0}$ such that the following holds for all $\d > 0$ small enough.

If $\nu$ is a Borel probability measure of the form $\nu=\nu'+\nu''$ with $\nu', \nu''$ non-negative measures such that $\|\nu'_{\delta}\|_{L^\infty}\leq \delta^{-\kappa}$,  
then for any integer $n \in [r\abs{\log \d},2r\abs{\log \d}]$, we have 
\[
\cW_\beta(\nu P^n, m_X) \leq  \delta^\eps +2\nu''(X).
\]
\end{proposition}

Recall that $(C^{0, \b}(X), \norm{\,\cdot\,}_{C^{0, \beta}})$ denotes the space of bounded $\b$-Hölder functions on $X$ and $\cW_{\b}$ the associated Wasserstein distance \eqref{def-wassertstein}. 

The spectral gap condition on $P$ is a way to say that $P$ determines a Markov chain that is exponentially mixing for the stationary measure $m_{X}$, while the condition that $P$ preserves  $\b$-Hölder functions reflects that $P$ does not distort too much the metric, e.g. $P$ satisfies a suitable exponential moment condition. The requirement for $\nu'$ expresses that $\nu'$ is not too far from $m_{X}$, in a homogeneous context it is equivalent to a high dimension condition at scale $\delta$. 

For the purpose of the paper, the important case for us lies in the following example. 

\medskip
\noindent {\bf Example}. 
Consider $X=G/\Lambda$ where $G$ is a connected semisimple Lie group with finite center and no compact factor, and $\Lambda$ is a lattice. 
Assume that $X$  is equipped with a quotient right-invariant metric, write $m_{X}$ the Haar probability measure. Let $\mu$ be a probability measure on $G$ whose support generates a Zariski-dense subgroup, and set $P=P_{\mu}$  the Markov operator associated to $\mu$.

In this setting, we know that the spectral radius of $P_\mu$ on $L^2_0(X,m_X)$ is strictly less than $1$.
This can be obtained by combining \cite[Lemma 3]{Bekka} of Bekka and \cite[Theorem C]{Shalom} of Shalom as explained in \cite[Proposition 3.2.5]{BenardPhD}.

The requirement that $P_\mu$ is a bounded operator on $(C^{0, \b}(X), \|\,\cdot\,\|_{C^{0, \beta}})$ is equivalent to say that $\mu$ has a finite exponential moment with $\int_{G} \norm{\Ad(g)}^\beta  \dd \mu(g) < \infty$.

Finally, let $\nu$ be a Borel probability measure on $X$. 
The requirement on $\nu$ in \Cref{endgame} holds if  $\nu$ is $(1-\frac{\kappa}{2\dim G}, \cB_{\delta}, \tau)$-robust and $\nu=\nu'+\nu''$ is an associated robust decomposition. 
Indeed, the regular part $\nu'$ of $\nu$ satisfies that for every $x \in X$ with $\inj(x)\geq \delta$,
\[ \nu'_\delta(x) = \frac{\nu'(B_\delta(x))}{m_X(B_\delta(x))} \ll_{G, \|\cdot\|} \delta^{-\kappa/2},\]
where we used that $m_X(B^X_\delta(x)) = m_G(B^G_\delta) \gg_{G, \|\cdot\|} \delta^{\dim G}$.
Thus $\norm{\nu'_\delta}_{L^\infty} \leq \delta^{-\kappa}$ for $\delta \lll_{G, \|\cdot\|, \kappa} 1$.

\bigskip

\begin{proof}[Proof of \Cref{endgame}] 
We  fix $\k'>0$ such that the operator norm satisfies
\begin{equation}
\label{eq:spectral gap}
\norm{P^n}_{L^2_0} \leq e^{-\kappa' n}
\end{equation}
for all large enough $n$. We also let $M$ be a parameter such that $M\geq\|P\|_{C^{0, \b}}$.  

As the function $\sigma\mapsto \cW_{\beta}(\sigma, m_{X})$ is convex on the set of probability measures, we may assume $\nu=\nu'$.  Let $f \in C^{0,\b}(X)$.
After shifting by a constant, we assume $\int f \dd m_X = 0$.
The goal is to bound from above
$\abse{\int_X P^n f \,\dd \nu}$.
We have
\begin{align*}
 \abse{\int_X P^n f \dd \nu} 
\leq \abse{\int_X P^n f \dd \nu_\d}  + \abse{\int_X P^n f \dd \nu_\d - \int_X P^n f \dd \nu}.
\end{align*}
The first term is bounded for large $n$ by:
\begin{align*}
\abse{\int_X P^n f \dd \nu_\d} & \leq \norm{P^n f}_{L^1} \norm{\nu_\d}_{L^\infty}\\
&\leq \norm{P^n f}_{L^2} \norm{\nu_\d}_{L^\infty}\\
&\leq \norm{P^n}_{L^2_0} \norm{f}_{L^2} \norm{\nu_\d}_{L^\infty}\\
&\leq e^{-\kappa' n} \delta^{-\kappa}\norm{f}_{L^2}.
\end{align*}
Using \eqref{defnudelta}, the second term is bounded by:
\begin{equation*}
\abse{\int_X P^n f \dd \nu_\d - \int_X P^n f \dd \nu} \leq \d^\b \norm{P^n f}_{C^{0,\b}}\leq \d^\b M^n \norm{f}_{C^{0,\b}}.
\end{equation*}

Put together, we obtain for large $n$:
\[
\abse{\int_X P^n f \dd \nu} \leq \bigl(e^{-\kappa' n} \delta^{-\kappa} +  \d^\b M^n\bigr)\norm{f}_{C^{0,\b}} \leq \delta^\eps \norm{f}_{C^{0,\b}} 
\]
where the final upper bound holds provided $n \in [\frac{\b}{4}\frac{\abs{\log\d}}{\log M},\frac{\b}{2}\frac{\abs{\log\d}}{\log M}]$,  $\kappa = \frac{\kappa' \b}{8\log M}$, $\eps = \frac{1}{2}\min\{\frac{\b}{2},\frac{\kappa' \b}{8\log M}\}$, and $\delta$ is small enough depending on $\eps$. Note that we can always choose $M$ large enough from the start so that the lower bound on $n$ is an arbitrary small multiple of $\abs{\log \delta}$, so this concludes the proof. 
\end{proof}

\section{Proof of the main statements} \label{Sec-proof-main-statements}
In this section, we prove the statements presented in the introduction.

\begin{proof}[Proof of  \Cref{thm:small-dim-equid-0} ]  
It is \Cref{small-dim-equid} specified to the case $\beta=1$.
\end{proof}

\begin{proof}[Proof of \Cref{main-thm}] 
Let $\kappa, A>0$ be constants depending on $X,\mu$ and for which \Cref{thm-positive-dimension} holds. Let $\eps>0$ be associated to $\kappa/2$ by \Cref{thm:small-dim-equid-0}. Then for every sufficiently large $R$, and every $m \geq \log R + A \max \{|\log \dist(x,W_{\!\mu, R})|,\, \dist(x,x_{0})\}$, 
\Cref{thm-positive-dimension} tells us that the measure $\mu^{*m}*\delta_{x}$ satisfies the non-concentration assumption in \Cref{thm:small-dim-equid-0} with parameters $\kappa/2$, $\eps$, and $\delta=R^{-1}$. Applying \Cref{thm:small-dim-equid-0}, we get for all $n\geq \log R $, and all $f \in \Lip(X)$ with $\norm{f}_{\Lip} \leq 1$, 
\[\abs{\mu^{*(n+m)}*\delta_{x}(f) -  m_{X}(f)} \leq R^{-\eps}+ 2\mu^{*m}*\delta_{x}\{\inj\leq R^{-\eps}\}.\]
The claim follows after applying \Cref{effective-recurrence} (with $\dist(.,x_{0})$ replaced by $\inj^{-1}$, see \Cref{comparison2}) to bound $\mu^{*m}*\delta_{x}\{\inj\leq R^{-\eps}\}$ by a small power of $R^{-1}$. 
\end{proof}

It remains to show \Cref{thm:equidis} and the corollaries presented in the introduction.

\subsection{Proof of Theorem \ref{thm:equidis}}

As explained in \Cref{Sec-ideas-more}, with \Cref{thm:small-dim-equid-0} at hand, we only need to complete phase I. 
That is, for a random walk starting at a point $x$  which is not trapped in a finite invariant set, we want to see a positive dimension at some range of scales $[\rho_{0},\rho_{0}^\eps]$ after some time.
This was achieved in \Cref{thm-positive-dimension} and with rate, at the cost of extra arithmetic assumptions. Here we make no such assumptions, but the rate becomes unspecified, it is allowed to depend arbitrarily on $x$ and $\rho_{0}$.

The proof of this result allows for homogeneous spaces modeled on more general groups. Namely, we show 

\begin{proposition}[Positive dimension without arithmetic condition] \label{thm-positive-dimension-without-arithm}
Let   $G$ be a connected semisimple real linear group with no compact factor, $\Lambda$ a lattice in $G$, set $X=G/\Lambda$ equipped with a quotient right $G$-invariant  Riemannian metric. Let $\mu$ be a Zariski-dense probability measure on $G$ with finite exponential moment. 

 There exists $C,\kappa>0$ such that, for every  $x\in X$ with infinite $\Gamma_{\mu}$-orbit,  $\rho_{0}>0$, we have for all large enough $n\geq 1$
\begin{equation}\label{eq-posdim-without-arithm}
\forall \rho\geq \rho_0,\, \forall y\in X,\quad \mu^{*n}*\delta_{x}(B_{\rho}(y)) \leq C\rho^\kappa.
\end{equation}
\end{proposition}

Given $x$ and  $\rho_{0}$, we argue by dividing the time into two times $n_1 + n_2$.
In the first time, atoms of $\nu_1 = \mu^{*n_1} * \delta_x$ all reach a mass smaller than $\rho_{0}$.
Then in the second time, using drift function arguments, namely \Cref{persist-ineq}, we show that
$\nu_2 = \mu^{* n_2} * \nu_1$ satisfies  \eqref{eq-posdim-without-arithm} for some $\kappa > 0$ depending only on $X$ and $\mu$.

The following lemma will allow us to find $n_1$.
\begin{lemma}\label{lm:atom orbit}
For every $x \in X$, either 
\[
\limsup_{n\to +\infty} \max_{y \in X} (\mu^{*n} *\delta_x)\{y\} = 0
\]
or $x$ is in a finite $\Gamma_{\!\mu}$-orbit.
\end{lemma}

\begin{proof}
For $n \in \N$, consider the function $f_n : X \to \R$ defined by
\[
f_n(x) = \max_{y \in X} (\mu^{*n} * \delta_x) \{y\} \quad \text{for} \quad  x \in X.
\]
Note that for every fixed $x \in X$, $f_n(x)$ is non-increasing in $n \in \N$.
This is because, for $n \geq m \geq 0$ and $y \in X$, we have
\[
(\mu^{*n} * \delta_x) \{y\} = \int_{\Gamma_{\!\mu}} (\mu^{*m} * \delta_x)\{g^{-1}y\} \dd \mu^{*(n-m)}(g).
\]

To prove the lemma, we let $x \in X$ be an element such that $f_n(x)$ does not converge to $0$ and we show that $\Gamma_{\!\mu} x$ is finite.
By monotonicity, $\inf_{n \in \N} f_n(x) > 0$, so we may fix $\eps > 0$  such that  $f_n(x) \geq 2 \eps$ for all $n \in \N$.

For $m \in \N$, consider the subset 
\[
A_m = \set{ z \in X : f_m (z) \geq \eps }.
\]
We claim that for large enough $m$, the set $A_m$ is discrete.

To show the claim, we use the following fact: for any $\eps' > 0$, for all $m$ sufficiently large,
\begin{equation}
\label{eq:Adtransient}
\forall v,w \in \kg \setminus\{0\},\quad  \mu^{*m} \set{g \in \Gamma_{\!\mu} : \Ad(g) v = w} < \eps'.
\end{equation}
This fact is well known, it follows for instance from \cite[Proposition 1.2]{Breuillard-note}, or \cite[Proposition 3.3, Theorem 1.3]{BreBeck25}.

We set $\eps' = \eps^2 /20$ and let $m \in \N$ be such that \eqref{eq:Adtransient} holds.
Then let $M \geq 1$ be large enough so that
\[
\mu^{*m} (B_M) \geq 1 - \eps',
\]
where
\[
B_M = \setbig{g \in \Gamma_{\!\mu} : \max\{ \norm{\Ad(g)}, \norm{\Ad(g^{-1})} \} \leq M}.
\]

Let $B \subset X$ be a ball whose radius, say $r> 0$, is small enough so that every point in $B$ has an injectivity radius at least $2M^2r$.
Since  $X$ is a locally finite union of such balls, it suffices to show that $A_m \cap B$ is finite.
More precisely, we show that $A_m \cap B$ has at most $2 / \eps$ points.
Otherwise, for each $z \in A_m \cap B$, using the definition of $A_m$, pick $y_z \in X$ such that $(\mu^{*m}*\delta_z)\{y_z\} \geq \eps$.
Set
\[
G_z = \set{g\in \Gamma_{\!\mu} : g z = y_z}
\]
so that $\mu^{*m}(G_z) \geq \eps$.
Using the inclusion-exclusion principle 
\[
1 \geq \mu^{*m}\bigl(\cup_z G_z\bigr) \geq \sum_{z} \mu^{*m}(G_z) - \sum_{z \neq z'} \mu^{*m}(G_z \cap G_{z'}),
\]
we see that if $\sharp (A_m \cap B) \geq 2 /\eps$, then there are two points $z\neq z' \in A_m \cap B$ such that
\[
\mu^{*m}(G_z \cap G_{z'}) \geq 2 \eps'.
\]
We can write $z' = \exp(v) z$ with some $v \in B_r^\kg \setminus \{0\}$.
Then for every $g,h \in G_z \cap G_{z'} \cap B_M$, we have
\begin{align*}
\exp(v) z = z' = g^{-1}hz'&= g^{-1}h \exp(v) z \\
&=\exp\bigl(\Ad(g^{-1}h)v\bigr) g^{-1}h z = \exp\bigl(\Ad(g^{-1}h)v\bigr) z
\end{align*}
Moreover,
\[
\max\{\, \norm{v}, \, \norm{\Ad(g^{-1}h)v} \, \}\leq M^2 r < \inj(z).
\]
Hence, $v = \Ad(g^{-1}h)v$ or equivalently,
\[
\Ad(g)v = \Ad(h) v.
\]
This being true for all $g \in G_z \cap G_{z'} \cap B_M$ which has measure $\mu^{*m}( G_z \cap G_{z'} \cap B_M) \geq 2\eps' - \eps'$, it contradicts \eqref{eq:Adtransient}.

This finishes the proof of the claim that $A_m$ is discrete.
Now, note that for every $n \in \N$ and any $y \in Y$, writing $X = A_m \cup (X \setminus A_m)$, we have
\begin{align*}
(\mu^{*(n+m)}*\delta_x)\{y\} &= \int_X (\mu^{*m} * \delta_z)\{y\} \dd (\mu^{*n} * \delta_x)(z)\\
&\leq (\mu^{*n} * \delta_x) (A_m) + \eps.
\end{align*}
Taking the maximum over $y \in Y$, we find
\[
\forall n \in \N, \quad 2\eps \leq f_{n+m}(x) \leq (\mu^{*n} * \delta_x) (A_m) + \eps.
\]

By the recurrence property \Cref{effective-recurrence}, there is a large compact set $K \subset X$ such that for all $n$ large enough, we have $(\mu^{*n} * \delta_x) (K) \leq \eps/2$.
It follows that, for $A \defeq A_m \cap K$ we have
\begin{equation*}
(\mu^{*n} * \delta_x) (A) \geq \eps /2 \quad \text{ for all  large enough $n$.}
\end{equation*}
But $A$ is finite, hence the existence of a point $z \in A$ such that
\begin{equation}
\label{eq:z is recurrent}
\liminf_{N \to +\infty} \Bigl( \frac{1}{N} \sum_{n=1}^N \mu^{*n} * \delta_x \Bigr) \{z\} > 0.
\end{equation}

Using the compactness of the space of probability measure on the one-point compactification of $X$, we can find a limit $\nu$ of a subsequence of $\bigl(\frac{1}{N} \sum_{n=1}^N \mu^{*n} * \delta_x\bigr)_{N \geq 1}$ in the weak-$*$ topology.
It is automatic that $\nu$ is a $\mu$-stationary measure.
By construction, the point $z$ is an atom for $\nu$.
By the maximum principle, $\Gamma_{\!\mu} z$ is a finite orbit.
But \eqref{eq:z is recurrent} also implies that $z$ is in the $\Gamma_{\!\mu}$ orbit of $x$.
Therefore $\Gamma_{\!\mu} x$ is a finite orbit, finishing the proof of the lemma.
\end{proof}

\begin{proof}[Proof of \Cref{thm-positive-dimension-without-arithm}]
Let $x \in X$ be a point whose orbit is infinite, and let $\rho_{0}>0$.
By  \Cref{lm:atom orbit}, there exists $n_1 \in \N$, such that the measure  $\nu_1 \defeq \mu^{*n_1} * \delta_x$ satisfies
\begin{equation}
\label{eq:y not atom}
\sup_{y \in X} \nu_{1} \{y\} < \rho_{0}.
\end{equation}
Then for all $r > 0$ small enough, we have 
\[
\nu_1\{ \inj \leq r \} \leq \rho_{0}, \,\,\,\,\,\,\,\,\,\,\,\,\,\,\,\sup_{y \in X} \nu_1 (B_\rho(y)) \leq \rho_{0}.
\]
Indeed, the first inequality is clear. If the second did not hold, then there would be points $y_{m} \in X$ such that $\nu_1 (B_{1/m}(y_{m})) > \rho_{0}$ for every $m \geq 1$.
Because $\nu_{1}$ has finite mass, the sequence $y_{m}$ must remain in a compact set. Let $y \in X$ be an accumulation point.
Then for every $t > 0$, $\nu_1 (B_t(y)) \geq \rho_{0}$ and hence $\nu_1\{y\} \geq \rho_{0}$, contradicting \eqref{eq:y not atom}.

Let $s > 0$ and $\lambda > 0$ be the constants given by \Cref{persist-ineq}.
Applying \Cref{persist-ineq} to the measure $\nu_1$, choosing $n_2 \geq - \frac{\log \rho}{s\lambda}$ and writing $\nu_2 = \mu^{*n_2}*\nu_1$, we have for every $r > 0$,
\[
\sup_{y \in X} \nu_2(B_r(y))^2 \ll r^s + \rho_{0},
\]
which concludes the proof. 
\end{proof}

We may finally conclude with the

\begin{proof}[Proof of \Cref{thm:equidis}] 
It is identical to the proof of \Cref{main-thm}, except the effective estimate \Cref{thm-positive-dimension} is replaced by \Cref{thm-positive-dimension-without-arithm}.
\end{proof}

\subsection{Proof of the corollaries}



\bigskip

\begin{proof}[Proof of \Cref{cr:dioph}]
$\ref{it:exprate1} \implies \ref{it:exprate2}$. Assume that $x$  is $(\mu,D)$-Diophantine  ($D>1$).  Let $\lambda, A$ as in \Cref{main-thm}. For  $n> A \log R +A^2D\log R + A \log D + A\cdot\dist(x, x_{0})$,  \Cref{main-thm} yields,
\[|\mu^{*n}*\delta_{x}(f)-m_{X}(f)| \leq \|f\|_{\Lip} R^{-1}.\]
which justifies the first implication, and the ``moreover'' part.

\noindent $\ref{it:exprate2}  \implies \ref{it:exprate1} $.  Assume $x\in X$  is not $(\mu,D)$-Diophantine for any $D$. In other words, there exist sequences $D_{i}\to +\infty$ and $R_{i}>1$ such that for all $i$, one has $W_{\!\mu, R_{i}}\neq \emptyset$ and
\[\dist(x, W_{\!\mu, R_{i}}) \leq \frac{1}{D_{i}}R_{i}^{-D_{i}}.\]
Let us show that \eqref{exp-eq} fails.
One may assume that $x$ is not trapped in a finite orbit. Note that necessarily $R_{i}\rightarrow +\infty$. 
Let $L_{\mu}>1$ such that for all $g\in \supp \mu$, $y,z\in X$, one has $\dist(gy, gz)\leq L_{\mu}\dist(y,z)$. Let $n_{i}=\frac{D_{i}}{2\log L_{\mu}}\log R_{i}$. Then $\mu^{*n_{i}}*\delta_{x}$ is supported on the $R_{i}^{-D_{i}/2}$-neighborhood of $W_{\!\mu, R_{i}}$.
We saw in \Cref{finite-orbits} that $W_{\!\mu, R_{i}}$ is $R_{i}^{-M}$-separated for some  $M=M(X, \mu)>0$. Hence for large $i$, we have that $\mu^{*n_{i}}*\delta_{x}$ is at most $R_{i}^{-M/2}$-equidistributed. As $D_{i}\to +\infty$, this forbids exponential equidistribution. 
\end{proof}

\bigskip

\begin{proof}[Proof of \Cref{cor-orbites}]
Let $\lambda, A$ as in \Cref{main-thm}, write $Y=\Gamma_{\!\mu}x$.  To prove the claim, one may assume $R$ to be large enough so that $R'\defeq (\frac{1}{2}R)^{1/A}\geq 2$.  Note that $\dist(x, W_{\!\mu, R'^A})>0$. We deduce from \Cref{main-thm}
 that for large enough $n$, 
for all $f\in C_{c}^\infty(X)$, 
\[|\mu^{*n}*\delta_{x}(f)-m_{X}(f)| \leq \|f\|_{\Lip}R'^{-1}.\]
The result follows by letting $n$ go to infinity and noting that $\mu^{*n}*\delta_{x}$ converges toward $m_{Y}$ in Ces\`aro average. 
\end{proof}

\appendix
\section{Non-linear subcritical projection theorem}
\label{ss:Appendix}
This appendix is dedicated to the proof of \Cref{pr:subcritical} in the case of $m = 1$ and $r_1 = 0$,
which is a nonlinear discretized projection theorem in the subcritical regime.

Let us restate it here. 
\begin{proposition}[Nonlinear subcritical projection theorem]\label{pr:CoarseNonProj}
Given an integer $d\geq 2$ and $\kappa > 0$,
there exists $C=C(d, \kappa) > 1$  such that the following holds for all $\eps \in (0, 1/2]$ and all $\delta \lll_{d,\kappa,\eps} 1$.

Let $(F_\theta)_{\theta \in \Theta}$ be a  family of differentiable maps $F_\theta \colon B^{\R^d}_1 \to \R^k$ where $0< k <d$ and such that $\theta \mapsto D_{x}F_{\theta}$ is measurable for every  $x\in B^{\R^d}_1$.
Let $\sigma$ be a probability measure on $\Theta$ and $A\subseteq B_1^{\R^d}$ a subset satisfying
\begin{enumerate}[resume*=cond]
\item \label{it:Fdistort} for $\sigma$-almost every $\theta \in \Theta$, every $x \in A$, all the singular values of $D_xF_\theta$ are between $\delta^\eps$ and $\delta^{-\eps}$ and moreover 
\[
\forall x, y \in A,\; \norm{F_\theta(x) - F_\theta(y) - D_x F_\theta (x-y)} \leq \delta^{-\eps} \norm{x-y}^2.
\]
\item \label{it:NCF} $\forall x\in A$, $\forall \rho \geq \delta$, $\forall W \in \Gr(\R^d, k)$,
\[
\sigma \set{\theta \in \Theta : \dang(\ker D_x F_\theta, W)\leq \rho} \leq \delta^{-\eps} \rho^\kappa.
\]
\end{enumerate}

Then the exceptional set  
\begin{equation}
\label{eq:excep-proj}
\begin{split} 
\cE \defeq \bigl\{\, \theta  \in \Theta : \exists A' \subseteq A \,\,&\text{ with }\,\,\cN_{\delta}(A')\geq \delta^{\eps} \cN_{\delta}(A)\\
& \text{ and } \,\,\cN_\delta(F_\theta A') < \delta^{C \eps \abs{\log \eps}} \cN_\delta(A)^{k/d} \,\bigr\}
\end{split}
\end{equation}
has measure $\sigma(\cE) \leq \delta^{\eps}$.
\end{proposition} 
Recall that the singular values of a $k \times d$ matrix $M$ are the square roots of the eigenvalues of the $k \times k$ matrix $M M^{\mathrm{T}}$. We can decompose $M \colon \R^d \to \R^k$ as the composition of a rotation of $\R^d$, the $k \times d$ matrix with the singular values on the main diagonal and zeros elsewhere, and a rotation of $\R^k$. 
 
Observe that the case $m=1$ and $\ur=(0,1)$ of \Cref{pr:subcritical} reduces to \Cref{pr:CoarseNonProj}. Indeed,  recall that  we are given a family $(\phi_\theta)_{\theta \in \Theta}$ of $C^1$-diffeomorphisms and a subspace $V_1 \subset \R^d$ of dimension $0 < \dim V_1 < d$.
We can then set $k = d - \dim V_1$ and $F_\theta = \pi_{\parallel V_1} \circ \phi_\theta$ for every $\theta \in \Theta$, where $\pi_{\parallel V_1} : \R^d \to \R^k$ denote the orthogonal projector of kernel $V_1$. 
Condition~\ref{it:distortion} implies immediately condition~\ref{it:Fdistort}.
Moreover, for any $\theta \in \Theta$ and any subset $A'$, we have $\cN_\delta(F_\theta A') \simeq \cN_\delta^{(0,1)}(\phi_\theta A')$ where $\cN_\delta^{\ur}$ is the notation used in \Cref{sec:slicing} with respect to the flag $0 \subsetneq V_1 \subsetneq \R^d$. Hence the reduction. On the other hand, \Cref{pr:CoarseNonProj} essentially boils down to \Cref{pr:subcritical} using the local normal form of submersions. 

\medskip
The proof is based on the linear case, which has a slightly better conclusion.
\begin{proposition} \label{A2}
\Cref{pr:CoarseNonProj} holds if the maps $F_\theta$ are all linear and $\cE$ is replaced by 
\begin{equation*}
\begin{split}
\cE \defeq \Bigl\{\, \theta  \in \Theta : \exists A' \subseteq A \,\,&\text{ with }\,\,\cN_{\delta}(A')\geq \delta^{\eps} \cN_{\delta}(A)\\
& \text{ and } \,\,\cN_\delta(F_\theta A') < \delta^{C \eps} \cN_\delta(A)^{k/d} \,\Bigr\}.
\end{split}
\end{equation*}
\end{proposition}
\begin{proof}
When $F_\theta$ is linear, condition~\ref{it:Fdistort} implies that for $\sigma$-almost every $\theta \in \Theta$, $F_\theta$ is the composition of a rotation in $\R^d$ with the projection to the first $k$ coordinates and then with a $\delta^{-\eps}$-bi-Lipschitz map of $\R^k$. 
Thus, we reduce the the setting where $F_\theta$ are orthogonal projections.
Then the statement follows from the proof of~\cite[Proposition 29]{He2020JFG}. 
\end{proof}

We upgrade \Cref{A2}  to its nonlinear counterpart \Cref{pr:CoarseNonProj}  in the same way as the supercritical projection of Bourgain has been upgraded to the nonlinear setting by Shmerkin in \cite{Shmerkin}. 
The extra $\abs{\log \eps}$ in \eqref{eq:excep-proj}  comes from the linearization procedure. The same phenomenon appears in Shmerkin \cite{Shmerkin} but there it can be ignored because in the supercritical statement, it is harmless to replace any occurrence of $\eps$ by a smaller number (provided $\delta$ is small enough accordingly). 

To reduce \Cref{pr:CoarseNonProj} to the linear case, it is more convenient to argue in terms of the Shannon entropy rather than covering numbers. Recall that, given a measurable space endowed with a probability measure $\nu$ and a countable partition $\cP$ into measurable subsets, the \emph{Shannon entropy} is 
\[H(\nu,\cP) \defeq - \sum_{P\in \cP} \nu(P) \log \nu(P).\]
It is essentially the logarithm of the covering number by $\cP$-cells of a set of large $\nu$-mass.

\begin{lemma}[Entropy vs covering number] \label{entropy-vs-cn}
For any $c\in (0, 1)$, 
\[
\log \cN_{\cP}(\supp \nu)\geq H(\nu,\cP)\geq (1-c) \inf \{\,\log \cN_{\cP}(E)\,:\,\nu(E)\geq c\,\}.
\]
\end{lemma}

\begin{proof}
The left inequality comes from the concavity of the logarithm. For the inequality on the right, set $h \defeq H(\nu,\cP)$, let $R>0$ be a parameter. By Markov's inequality, 
$\nu (\cup\{P\,:\, -\log \nu(P) > R h \})  \leq 1/R$, so witing $A' \defeq \cup\{P\,:\, -\log \nu(P)\leq R h \}$ the complementary set,  we get $\nu(A') \geq 1-1/R$. For $P\in \cP(A')$,  we have $\nu(P)\geq e^{-R h}$ so $\cN_{\cP}(A')\leq e^{Rh}$, i.e. $\log \cN_{\cP}(A') \leq R \,H(\nu,\cP)$. The result follows choosing $R=(1-c)^{-1}$. 
\end{proof}

The next proposition is taken from \cite{Shmerkin}.
It shows how to bound from below the entropy of a pushforward measure $F_{\star} \nu$ in terms of entropies of local pieces, where the map $F$ has been replaced at each scale and local neighborhood by linear projector approximations. 

We use the notation introduced in \Cref{sec:slicing}. Abusing slightly, $\cD_\delta$ can denote either the cubic tiling in $\R^d$ of side length $\delta$ rounded to a power of $2$, or that in $\R^k$.
For a measure $\nu$ on $\R^d$ and a dyadic cube $Q \in \cD_{\rho}(\nu)$, we let 
$\nu^Q \defeq (\Delta_Q)_\star (\nu_Q)$ where $\nu_Q$ is the normalized restriction of $\nu$ to $Q$ and $\Delta_Q$ is the affine dilation that sends bijectively $Q$ to the unit cube.

\begin{proposition}[Linearization] 
\label{pr:linearize}
Let $1 \leq k < d$.
Let $\nu$ be a probability measure supported on $B^{\R^d}_1$, let $F \colon U \to \R^k$ be a differentiable map defined on a neighborhood of $\supp \nu$. 
Assume that there is a constant $L \geq 1$ such that for every point $x \in \supp \nu$, all $k$ singular values of $D_xF$ are between $L^{-1}$ and $L$ and moreover
\[
\forall x, y \in \supp \nu,\quad \norm{F(x) - F(y) - D_xF(x - y)} \leq L \norm{x - y}^2.
\]

Let $\delta, \delta_1,\dotsc, \delta_q, \rho_1, \dotsc, \rho_q$ be powers of $2$ such that 
\[
\forall i \in \{1,\dotsc,q\},\quad \rho_i \leq \delta_i \leq 1
\]
and 
\[
\delta \leq \delta_1 \rho_1 \leq \rho_1 \leq \delta_2 \rho_2 \leq \rho_2 \leq \dots \leq \delta_q \rho_q \leq \rho_q \leq 1.
\]
For each $Q \in \cD_{\rho_i}(\nu)$, $1 \leq i\leq q$, pick an arbitrary point $x_Q \in Q \cap \supp\nu$.
As an approximation of $F$ on $Q$, consider the orthogonal projection parallel to $\ker D_{x_Q} F \in \Gr(\R^d,d-k)$ and denote it by $\pi_Q$. 
Then,
\begin{equation*}
H(\nu, F^{-1}\cD_\delta) \geq -3k q \log L - O_d(q) + \sum_{i=1}^{q} \sum_{Q \in \cD_{\rho_i}} \nu(Q) H \bigl( \nu^Q, \pi_{Q}^{-1} \cD_{\delta_i}\bigr).
\end{equation*}
\end{proposition}

\begin{proof}
This is a restatement of \cite[Proposition A.1]{Shmerkin} with explicit dependence  on $L$ of the deficit term.
We summarize the proof to trace the effect of $L$.

By basic properties of the conditional entropy, we have
\[
H(F_\star \nu, \cD_\delta) \geq \sum_{i=1}^q H(F_\star \nu, \cD_{\delta_i\rho_i} | \cD_{\rho_i}).
\]
For each $i = 1,\dotsc, q$, write $\nu = \sum_{Q \in \cD_{\rho_i}} \nu(Q) \nu_Q$.
By the concavity of the conditional entropy,
\[
H(F_\star \nu, \cD_{\delta_i\rho_i} | \cD_{\rho_i}) \geq  \sum_{Q \in \cD_{\rho_i}} \nu(Q) H(F_\star \nu_Q,\cD_{\delta_i\rho_i} | \cD_{\rho_i}).
\]
Note that by the assumptions on $F$, for each $Q \in \cD_{\rho_i}(\nu)$, $F (\supp \nu_Q)$ has diameter at most $O_{k}(L \rho_i)$.
It follows that $\sharp \cD_{\rho_i}(F_\star \nu_Q) \ll_k L^k$ and hence using \Cref{entropy-vs-cn},
\begin{align*}
H(F_\star \nu_Q,\cD_{\delta_i\rho_i} | \cD_{\rho_i}) &=  H(F_\star \nu_Q,\cD_{\delta_i\rho_i}) - H(F_\star \nu_Q, \cD_{\rho_i})\\
&\geq H(\nu_Q, F^{-1}\cD_{\delta_i\rho_i}) - k \log L - O_k(1).
\end{align*}
We claim that we have the following rough refinement relation (recall notation from \Cref{Sec-prepa}),
\[
(\pi_Q^{-1} \cD_{\delta_i\rho_i})_{|\supp \nu_Q} \overprec{O(L^{2k})} ( F^{-1}\cD_{\delta_i\rho_i})_{|\supp \nu_Q}.
\]
It follows that
\begin{align*}
H(\nu_Q, F^{-1}\cD_{\delta_i\rho_i}) &= H(\nu_Q, \pi_Q^{-1}\cD_{\delta_i\rho_i}) - H(\nu_Q, F^{-1}\cD_{\delta_i\rho_i} | \pi_Q^{-1}\cD_{\delta_i\rho_i}) \\
& \geq H(\nu^Q, \pi_Q^{-1}\cD_{\delta_i}) - 2 k \log L - O_d(1).
\end{align*}
Combining all the inequalities above proves the desired estimate.

It remains to show the claim.
If $x,y \in \supp\nu_Q$ are in the same cell for the partition $F^{-1}\cD_{\delta_i \rho_i}$, then
\[
\norm{F(x) - F(y)} \ll_k \delta_i \rho_i.
\]
By the assumption on $F$, 
\[
\norm{F(x) - F(y) - D_{x_Q}F(x-y)} \leq L (\norm{x - x_Q}^2 + \norm{y - x_Q}^2) \ll_d L \rho_i^2 \leq L \delta_i \rho_i
\]
and by the assumption on the singular values of $D_{x_Q}F$, $\pi_Q$ is the composition of $D_{x_Q} F$ with a $L$-bi-Lipschitz map, whence
\[
\norm{ \pi_Q (x-y)} \leq L \norm{D_{x_Q}F(x -y)} \ll_d L^2 \delta_i \rho_i.
\]
This shows the claim.
\end{proof}

\medskip
We are all set to prove \Cref{pr:CoarseNonProj}.
\begin{proof}[Proof of \Cref{pr:CoarseNonProj}]
Let $\eps, \delta \in (0,1/2]$ and let $q \geq 1$ be an integer to be determined in terms of $\eps$.
For $0 \leq i \leq q$, set $\delta_i = \delta^{2^{-i}}$. We may assume throughout the proof that \ref{it:Fdistort} holds for \emph{every} $\theta\in \Theta$. 


Similarly to the proof of \Cref{thm:slicing}, using the regularization procedure from \Cref{tree-structure} and the exhaustion technique of \cite[Proposition 25]{He2020JFG}), we may assume without loss of generality that $A$ is regular with respect to the filtration $\cD_{\delta_q} \prec \dots \prec \cD_{\delta_1} \prec \cD_\delta$.
Without loss of generality assume further that $A \cap P$ is a singleton for every $P \in \cD_\delta(A)$.
Let $\nu$ be the normalized counting measure on $A$.

Let $A' \subset A$ be such that $\sharp A' \geq \delta^\eps \sharp A$, 
 denote by $\nu'$  the normalized counting measure on $A'$. Let $\theta\in \Theta$.  
Using condition~\ref{it:Fdistort} and \Cref{pr:linearize} applied with $\rho_i = \delta_i = \delta^{2^{-i}}$, $1 \leq i \leq q$, we have for any $q \geq 1$, 
\begin{equation}
\label{eq1-subc}
H(\nu', F_{\theta}^{-1}\cD_{\delta}) + O_{d}(q \eps \abs{\log \delta}) \geq \sum_{i=1}^{q} \sum_{Q \in \cD_{\delta_i}} \nu'(Q) H \bigl( \nu'^Q, \pi_{Q, \theta}^{-1} \cD_{\delta_i}\bigr).
\end{equation}
where $\pi_{Q, \theta}$ denotes the orthognal projection parallel to $\ker D_{x_Q}F_\theta$ with $x_Q$ being an arbitrary point chosen in $Q \cap A$, just as in \Cref{pr:linearize}.

For every $i \in \{1,\dotsc,q\}$ and $Q \in \cD_{\delta_i}$, by \Cref{entropy-vs-cn}, 
\begin{equation}\label{eq2-subc}
H \bigl( \nu'^Q, \pi_{Q,\theta}^{-1} \cD_{\delta_i}) \geq  \inf_{E\subseteq [0, 1)^d \,:\, \nu'^Q(E)\geq \delta^{\eps}} \log \cN_{\delta_i}(\pi_{Q, \theta} E) - \delta^{\eps} d \abs{\log \delta_i}
\end{equation}

We now see that for most $Q$, such $E$ is actually large inside $\Delta_{Q}(A\cap Q)$. This will allow us to apply the linear subcritical projection theorem (\Cref{A2}) to bound below the infimum in \eqref{eq2-subc}.
For each $i$, consider 
\[\cQ^{(i)}_{\mathrm{large}}(A')\defeq \{\, Q\in \cD_{\delta_i} \,:\, \sharp(A'\cap Q)\geq \delta^{3\eps}\sharp(A\cap Q)\,\}.\]
Observe that for every $Q\in \cQ^{(i)}_{\mathrm{large}}(A')$, the lower bound $\nu'^Q(E)\geq \delta^{\eps}$ implies $\nu^Q(E)\geq \delta^{4\eps}$, which in turns implies $ \cN_{\delta_i}(E) \geq \delta^{4\eps} \cN_{\delta_i}(\Delta_{Q}(A\cap Q))$ thanks to \Cref{lm:reg and subset}.
Moreover,  $\cQ^{(i)}_{\mathrm{large}}(A')$ has almost full $\nu'$-mass, indeed
\begin{align*}
\sharp A'  &\leq  \delta^{3\eps}\sharp A + \sum_{\cQ^{(i)}_{\mathrm{large}}(A')} \sharp (A'\cap Q) \leq  \delta^{3\eps}\sharp A+ \nu'(\cup\cQ^{(i)}_{\mathrm{large}}(A')) \sharp A'
\end{align*} 
yielding, by the condition $\sharp A' \geq \delta^\eps \sharp A$, that 
\begin{equation}
\nu'(\cup \cQ^{(i)}_{\mathrm{large}}(A')) \geq 1-\delta^{2\eps}. \label{eq3-subc}
\end{equation}

For each $Q\in \cD_{\delta_i}(A')$, we now apply the \emph{linear} subcritical projection theorem  (\Cref{A2}) with exponent $\eps_{i}\defeq 2^{i+2}\eps$ and at scale $\delta_i$ to the set $\Delta_{Q}(A\cap Q)$ and random projectors $(\pi_{Q, \theta})_{\theta\sim \sigma}$.
This gives a constant $C = C(d, \kappa)>1$ and an event $\cE_{Q}\subseteq \Theta$  of mass $\sigma(\cE_{Q})\leq \delta_i^{\eps_i}=\delta^{4\eps}$ and such that for all $\theta \notin \cE_{Q}$, provided that $\delta \lll_{d, \kappa, \eps, i} 1$,  
\[
\inf_{\cN_{\delta_i}(E) \geq \delta^{4\eps}\cN_{\delta_i}(\Delta_{Q}(A\cap Q))} \log \cN_{\delta_i}(\pi_{Q, \theta} E)    \geq \frac{k}{d}\log \cN_{\delta_{i}^2}(A\cap Q) - C\eps \abs{\log \delta}
\]
and in particular, if $Q\in \cQ^{(i)}_{\mathrm{large}}(A')$, 
\[
\inf_{\nu'^Q(E)\geq \delta^{\eps}} \log \cN_{\delta_i}(\pi_{Q, \theta} E) \geq \frac{k}{d}\log \frac{\cN_{\delta_{i-1}}(A)}{\cN_{\delta_i}(A)} - C\eps \abs{\log \delta},
\]
where we have also used the regularity of $A$ between $\cD_{\delta_i} \prec \cD_{\delta_{i-1}}$.

Moreover, Fubini's theorem implies that this estimate holds for most $\theta$ and most $Q$ simultaneously:
\[\cF_{i}\defeq \setbig{\theta \in \Theta \,:\, \nu \bigl( \cup \set{ Q \in \cD_{\delta_i} : \theta \notin \cE_Q}  \bigr) \geq 1-\delta^{2\eps} } \,\text{ satisfies } \,\sigma(\cF_{i})\geq 1-\delta^{2\eps}.\]
Note that $\cF_{i}$ does not depend on $A'$. On the other hand for $\theta \in \cF_{i}$, we have $\theta \notin \cE_{Q}$ for most $Q$ chosen with $\nu'$ as well; indeed $\nu'\leq \delta^{-\eps} \nu$ implies 
\[\cF_{i}\subseteq \bigl\{\, \theta \,:\, \nu'  \bigl( \cup \set{ Q \in \cD_{\delta_i} : \theta \notin \cE_Q} \bigr) \geq 1-\delta^{\eps} \,\bigr\}.
\]

To summarize, for fixed $i$, if $\theta \in \cF_i$, we can restrict the sum over $Q$ in \eqref{eq1-subc} to $\set{Q \in \cQ_{\mathrm{large}}^{(i)} : \theta \notin \cE_Q}$, recalling \eqref{eq2-subc} and \eqref{eq3-subc},
\[
\sum_{Q \in \cD_{\delta_i}} \nu'(Q) H \bigl( \nu'^Q, \pi_{Q, \theta}^{-1} \cD_{\delta_i}\bigr) \geq (1 - 2 \delta^\eps) \frac{k}{d} \log \frac{\cN_{\delta_{i-1}}(A)}{\cN_{\delta_i}(A)} - (C\eps + 2^{-i}\delta^{\eps}d)\abs{\log \delta}.
\]
Summing over $i$ and using \eqref{eq1-subc}, we obtain for $\theta\in \cap_{i\leq q}\cF_{i}$, and under the condition $\delta \lll_{d, \kappa, \eps, q} 1$, that
\begin{align*}
&\log \cN_{\delta}(F_{\theta}A') + O_{d}(q \eps \abs{\log \delta}) \\
\geq &H(\nu',F_{\theta}^{-1} \cD_\delta) + O_{d}(q \eps \abs{\log \delta}) \\
\geq & (1-2\delta^{\eps}) \frac{k}{d} \log \cN_{\delta}(A) - \log \cN_{\delta_q}(A) - (q C\eps + \delta^{\eps}d)\abs{\log \delta}\\
\geq & \frac{k}{d}\log \cN_{\delta}(A) -O_d(q C \eps+ \delta^\eps + 2^{-q}) \abs{\log \delta}. 
\end{align*}
Choosing $q$ to be the smallest integer such that $2^{-q}\leq \eps$, i.e. $q \defeq \lceil \log_{2}\eps\rceil$, we obtain the relevant lower bound. The formula for $q$ also justifies that the set $\cap_{i\leq q}\cF_{i}$ of good parameters $\theta$ satisfies $\sigma(\cap_{i\leq q}\cF_{i})\geq 1-\delta^{\eps}$ provided $\delta\lll_{\eps}1$. This concludes the proof.
\end{proof}

\bibliographystyle{abbrv} 
\bibliography{EffEq}

\end{document}